\documentclass[a4paper,11pt,oneside]{article}
\usepackage[top=2cm, bottom=2cm, left=2cm, right=2cm]{geometry}

\usepackage{amsmath,amsthm,amsfonts,amssymb,amscd}
\usepackage{subfig}
\usepackage{booktabs}
\usepackage{mathtools}
\usepackage[dvipsnames]{xcolor}
\usepackage{graphicx,bm}
\usepackage{algorithm}
\usepackage{algpseudocode}
\usepackage{url}
\usepackage{siunitx}
\usepackage{hyperref}

\usepackage{authblk}

\newtheorem{theorem}{Theorem}[section]
\newtheorem{proposition}[theorem]{Proposition}

\newtheorem{remark}[theorem]{Remark}
\newtheorem{definition}[theorem]{Definition}

\newtheorem{problem}[theorem]{Problem}

\theoremstyle{definition} %

\theoremstyle{remark} %

	\title{Numerical modeling and open-source implementation of variational partition-of-unity localizations of space-time dual-weighted residual estimators for parabolic problems}
	\author[1,2]{J.P. Thiele}
	\author[1,3,4]{T. Wick}
	
        \affil[1]{Leibniz University Hannover,
	Institute of Applied Mathematics,
	Scientific Computing
	Welfengarten 1, 30167 Hannover, Germany}

        \affil[2]{Weierstrass Institut f\"ur Angewandte Analysis 
        und Stochastik, Mohrenstr. 39, 10117 Berlin, Germany}
        
        \affil[3]{Cluster of Excellence PhoenixD (Photonics, Optics, and
        Engineering - Innovation Across Disciplines), Leibniz University Hannover, Germany}
	
        \affil[4]{Universit\'e Paris-Saclay, LMPS - Laboratoire de Mecanique Paris-Saclay, 91190 Gif-sur-Yvette, France}

	\date{}

\begin{document}\maketitle
\begin{abstract}
In this work, we consider space-time goal-oriented a posteriori 
error estimation for parabolic problems. Temporal and spatial discretizations 
are based on Galerkin finite elements of continuous and discontinuous type.
The main objectives are the development and analysis of space-time 
estimators, in which the localization
is based on a weak form employing a partition-of-unity. The resulting 
error indicators are used for temporal and spatial adaptivity. 
Our developments are substantiated with several numerical examples.
\end{abstract}
	
\section{Introduction}
\label{Section: Introduction}
Space-time methods for solving differential equations with Galerkin-type finite 
element methods go back to \cite{Od69b} 
and a recent state-of-the-art summary was compiled in \cite{LaStein19}.
Space-time methods can be divided into 
two categories: the numerical solution and error estimation. In this work, we are primarily 
interested in the latter. 
After the previously mentioned early work, various problems have been considered in 
space-time formulations such as
incompressible flow \cite{tezduyar.etal1992},
first-order hyperbolic systems \cite{DoerflerFindeisenWieners+2016+409+428}, 
elastic wave equation \cite{hughes.hulbert1988,HuHu90,Koecher2015}, 
visco-acoustic/visco-elastic wave equations \cite{DoeWieZie21}, 
financial mathematics \cite{goll.etal2015}, 
the Biot equations in poroelasticity \cite{bause.etal2017},
and fluid-structure interaction 
\cite{hubner.etal2004,tezduyar.etal2006,tezduyar.etal2006a,TeSa07,TaTe11,failer.wick2018}.
Advancements for the numerical solution by means of space-time methods, 
such as for example multigrid methods, were undertaken in 
\cite{neumuller2013,gander.neumuller2016,Schaf22} and 
with space-time domain decomposition \cite{SINGH2018893}.

Classical norm-based a posteriori error estimation was done for parabolic problems
in \cite{ErJo91,ErJo95,Ve03,SteiYa18,LaSchaf20,LaSchaf21}.
Goal-oriented error estimation of space-time problems 
was performed in \cite{meidner2007,schmich.vexler2008,schmich2009,besier.rannacher2012}. Therein,
space-time formulations may serve three purposes: 
spatial error estimation, 
temporal error estimation 
\cite{MeiRi14,meidner.richter2015,failer.wick2018,goll.etal2015} 
or both simultaneously \cite{schmich.vexler2008,schmich2009,bangerth.etal2010,
besier.rannacher2012,DoerflerFindeisenWieners+2016+409+428,DoeWieZie21}. 
Decoupling space and time for rate-dependent problems in 
elasto-plasticity was considered in \cite{rannacher.suttmeier1999}. 
Moreover, we mention space-time developments in PDE-based optimization with 
and without a posteriori error control 
\cite{meidner2007,becker.etal2007,meidner.vexler2007,NeiVe12,failer.etal2016,Fai17,langer.etal2021,langer.etal2021a,LaSchaf21b,KhiSteiWi22_JCP}. 
A brief review of space-time concepts for goal-oriented a posteriori error estimation in fluid-structure interaction for deriving 
the adjoint in goal-oriented error estimation and optimization was conducted in \cite{Wi21_WCCM}.

Employing the dual-weighted residual method 
\cite{becker.rannacher1996,becker.rannacher2001}
for space-time goal-oriented error estimation comes with the challenge that
the adjoint problem is running backwards-in-time. For nonlinear problems, this means 
that the primal solution must be available at the respective time points. 
This can be done by simply storing all primal solutions in the RAM (random access memory) 
or hard disk, 
or by using checkpointing techniques \cite{schmich2009,meidner2007}. 

The main objective in this work is to combine space-time concepts 
from \cite{schmich.vexler2008} with an easy-to-implement partition-of-unity 
localization proposed in \cite{richter.wick2015}. 
The latter was established 
for stationary problems, which is extended in this study to 
space-time error estimation. Two error estimators are proposed: joint and split.
Because of Galerkin orthogonality, we need higher-order information in the adjoint problem 
for calculation of the primal residual estimator.  
There are different ways to achieve this. 
For stationary problems a mixed order approach is often used.
There, we discretize the primal problem with the low order 
$cG(s)dG(r)$ elements and the adjoint problem with high order $cG(s+1)dG(r+1)$ elements.
The notation was proposed in \cite{ErJoLo04,ErEstHaJoh09} and means that 
spatial discretization is based on $cG(s)$ or $cG(s+1)$ continuous Galerkin finite elements respectively,
where $s\in\mathbb{N}$ indicates the polynomial degree,
while temporal discretization is based on $dG(r)$, discontinuous Galerkin finite elements,
where $r\in\mathbb{N}_0$ indicates the temporal polynomial degree. 

This approach needs interpolation operators to calculate the low order adjoint solution. 
In the equal low order approach both problems are discretized with low order elements 
and the high order solutions are obtained by a suitable patch wise reconstruction operator. 
For the adjoint estimator higher-order information in the primal problem is needed as well. 
For the two earlier approaches this can again be calculated by patch wise reconstruction. 
Alternatively, in the equal high order approach both problems can be discretized with high order elements.
Then, only interpolation operators are needed but the whole solution becomes more expensive.
Additionally, these approaches can be mixed by using different approaches for temporal and spatial discretization.
From the resulting a posteriori error estimates, local error indicators 
are extracted to establish adaptive algorithms for both temporal and spatial
mesh refinement. For verification, error reductions and effectivity indices are observed.
Some preliminary results were published in the conference proceedings 
papers \cite{ThiWi21_PAMM} (heat equation) and \cite{ThiWi23_ICO} (low Mach number combustion). 
Moreover, the successfull application to incompressible flow 
is documented in \cite{RoThieKoeWi23} and a summary of all developments 
will appear in \cite{Thi23_phd}.
However, technical derivations and the theory have not been worked out therein. Moreover, 
the current work provides (for the first time) detailed computational comparisons 
of the joint and split error estimators in terms of effectivity indices as well as thorough 
investigations of the PU-DWR method in a space-time context.

The outline of this paper is as follows: In Section \ref{sec_space_time}, 
the primal problem statements are provided including their space-time discretizations.
Next, in Section \ref{Section: space-time DWR}, the dual-weighed residual 
method is recapitulated. Afterward in Section \ref{Section: Localization},
partition-of-unity DWR space-time goal-oriented a posteriori error estimators
are proposed. In Section \ref{Section: Numerical tests}, three numerical examples
are studied in order to substantiate our algorithmic developments.
We summarize our work in Section \ref{Section: Conclusions}.

\section{Space-time notation and problem formulations}
\label{sec_space_time}
In this section, we introduce notation and space-time spaces. 
Then, abstract forms 
on the continuous level, semi-discrete in time level, and fully discrete level are introduced.
These are subsequently realized with specific problem statements, namely the heat equation 
and a combustion problem, respectively. 
The road map of our developments is first based 
on abstract derivations as we believe that with this knowledge our results can be more easily applied to 
other problem statements such as incompressible flow, e.g., as already done in \cite{RoThieKoeWi23}, 
and further nonstationary, nonlinear, coupled problems.

\subsection{Notation and spaces}
The space time domain is denoted by $\Sigma\subset \mathbb{R}^{d+1}$
with $\Sigma = \{\Omega(t)\subset\mathbb{R}^d\colon t \in I \}$, with the temporal interval $I=(0,T)$ and
spatial domain $\Omega$.
In this paper, we consider time-independent domains $\Omega$, 
resulting in $\Sigma = \Omega\times I$.

Using $V=V(\Omega)= H^1_D(\Omega) = \{v\in H^1(\Omega) ; v|_{\Gamma_D}=0\} $ and $H = H(\Omega)= L^2(\Omega)$ we can define our space-time Hilbert space as
\begin{align}
X \coloneqq X(I,V)\coloneqq W((I,V(\Omega))\coloneqq \{v\colon v\in L^2(I,V(\Omega))\text{ and }\partial_t v\in L^2(I,V^*(\Omega)) \}
\end{align}
where $L^2(I,V(\Omega))$ is the Bochner space of $L^2$ functions over $I$ with values in $V(\Omega)$.\\
Here, $\Gamma_D\subseteq\partial\Omega$ denotes the Dirichlet boundary with condition $v(t)=0\text{ on }\Gamma_D$.
For inhomogeneous conditions i.e. $v(t)=g(t)\text{ on }\Gamma_D$ the ansatz space is $X(I,V)+g$, we will however limit 
derivations to the homogeneous case for the sake of brevity.
On $X(I,V)$ we can use the $L^2(I,L^2(\Omega))$ scalar product 
\begin{align}
(u,v)\coloneqq (u,v)_{L^2(I,H)} = \int\limits_0^T (u(t),v(t))_{H} \mathrm{d}t.
\end{align}

Since we want to use discontinuous Galerkin discretizations in time we have to define an infinite dimensional space 
$\widetilde{X}(\mathcal{T}_k,V(\Omega))$ that allows for jumps at the grid points of the temporal mesh $\mathcal{T}_k$.\\
To obtain $\mathcal{T}_k$ we decompose the temporal interval $I$ into $M$ open subintervals $I_m\coloneqq (t_{m-1},t_m)$ 
of length $k_m=t_m-t_{m-1}$, with the condition that
\begin{equation}
 \bar{I} = \bar{I}_1\cup\bar{I}_2\cup\dots\cup\bar{I}_M\quad\text{and}\quad I_i\cap I_j = \{ \}\; \forall i\neq j
\end{equation}
hold. Then, $\mathcal{T}_k$ is the collection of all intervals from $I_1$ to $I_M$.\\
Following the nomenclature of \cite{DiPietroErn2011} 
we can define the \textit{broken} Bochner space
\begin{equation}
 \widetilde{X}(\mathcal{T}_k,V(\Omega))\coloneqq\{v\in L^2(I,L^2(\Omega))\text{ and }v|_{I_m}\in W(I_m,V(\Omega))\;\forall I_m\in\mathcal{T}_k\}.
\end{equation}
For each local space the continuous embedding $W(I_m,V(\Omega))\subset C(\bar{I}_m,H(\Omega))$ holds such that the limits from above and below,
i.e. $v_m^\pm\coloneqq\lim\limits_{\varepsilon\to0} v(t_m\pm\varepsilon)$ are well defined.\\
Using these we can define the jump across temporal intervals as
\begin{equation}
 [v]_m\coloneqq v_m^+-v_m^-
\end{equation}
for all inner temporal grid points $t_m\in\mathcal{T}_k$.\\
Additionally, we introduce $[v]_0$ as a shorthand for the weakly imposed initial conditions, i.e.
\begin{equation}
 [v]_0\coloneqq v_0^+-v^0.
\end{equation}
Finally, let $A:X\times X\to\mathbb{R}$ be a semi-linear 
form representing the PDE (partial differential equation) in weak form, 
being nonlinear in the first argument and linear in the second argument.
Let $F:X\to\mathbb{R}$ be a linear form representing given right hand side data. Then,
the abstract problem reads: 
\begin{problem}[Abstract form on the continuous level]
Find $u\in X$ such that 
\begin{equation}
\label{eq_weak_form_abstract}
A(u)(\varphi) = F(\varphi) \quad\forall\varphi\in X.
\end{equation}
We provide strong forms and the respective weak forms in 
terms of \eqref{eq_weak_form_abstract} for the heat equation and a combustion 
problem in Section \ref{sec_heat} and Section \ref{sec_combustion}, respectively. 
\end{problem}

\subsection{Discretization}
In principle, the discretization steps are the same for all parabolic problems. 
Since we want to be able to have different trial functions in time and space, we split the 
discretization, starting with the temporal basis functions.

\subsubsection{Semi-discretization in time}
Given the temporal triangulation $\mathcal{T}_k$ as defined above we can obtain the 
semidiscrete space by discretization of the temporal functions into piecewise polynomials of degree 
$r\in\mathbb{N}_0$:
\begin{align}
\widetilde{X}_k^r(\mathcal{T}_k,V)&\coloneqq\left\{v_k\in L^2(I,H)\text{ and }v_k|_{I_m}\in\mathcal{P}_r(I_m,V)
	\right\}\subset \widetilde{X}(\mathcal{T}_k,V).
\end{align}

With $A(u_k)(\varphi)$ and $F(\varphi)$ depending on the actual problem,
the general time-discrete weak $dG(r)$ formulations reads:
\begin{problem}[Abstract form semi-discrete in time]
Find $u_k\in \widetilde{X}_k^r(\mathcal{T}_k,V)$, where $r\geq 0$, such that
\begin{align}
A(u_k)(\varphi_k) = F(\varphi_k)\;\forall\varphi_k\in\widetilde{X}_k^r(\mathcal{T}_k,V). 
\end{align}
\end{problem}

\subsubsection{Fully discrete abstract problem}
For the spatial discretization we use triangulations $\mathcal{T}_h^m$ of $\Omega$,
where $m=1,\ldots, M$ indicates each temporal subinterval $I_m$.
These are decomposed into 
quadrilateral/hexagonal elements $K$ and we use continuous 
test and trial functions of degree $s$ resulting in $V_h^s(\mathcal{T}_h^m)\subset V(\Omega)$ defined as
\begin{align}
 V_h^{s}(\mathcal{T}_h^m) \coloneqq \{v_h\in V\text{ and }v_h|_K\in Q_s(K) \forall K\in\mathcal{T}_h^m\}.
\end{align}
We notice that we can have different triangulations on different subintervals, 
resulting in time-dependent or dynamic meshes.
Using this, we can define the fully discrete function space :
\begin{align*}
\widetilde{X}_{k,h}^{r,s}(\mathcal{T}_k,\mathcal{T}_h^{1,\dots,M}) &:= 
\{ v_{kh} \in L^2(I,H)\text{ and }v_{kh}|_{I_m} \in P_r(I_m,V_h^{s}(\mathcal{T}_h^m))\;\forall I_m\in\mathcal{T}_k\}
\subset\widetilde{X}_k^r(\mathcal{T}_k,V).
\end{align*}
With these definitions, we obtain the 
fully discrete $cG(s)dG(r)$ formulation:
\begin{problem}[Abstract form fully discrete level]
Find $u_{kh}\in\widetilde{X}_{k,h}^{r,s}(\mathcal{T}_k,\mathcal{T}_h^{1,\dots,M})$, such that
\begin{align}
\label{eq_weak_form_abstract_fully_discrete}
A(u_{kh})(\varphi_{kh}) = F(\varphi_{kh})\quad\forall\varphi_{kh}\in
\widetilde{X}_{k,h}^{r,s}(\mathcal{T}_k,\mathcal{T}_h^{1,\dots,M}).
\end{align} 
Specific realizations are provided in the following two subsections, 
Section \ref{sec_heat} and Section \ref{sec_combustion}, respectively, by setting 
\begin{align*}
A(u_{kh})(\varphi_{kh}) &:= A_{\text{heat}}(u_{kh})(\varphi_{kh}), \qquad
F(\varphi_{kh}) := F_{\text{heat}}(\varphi_{kh}),\\
A(u_{kh})(\varphi_{kh}) &:=A_{\text{comb.}}(u_{kh})(\varphi_{kh}), \qquad
F(\varphi_{kh}) := F_{\text{comb.}}(\varphi_{kh}).
\end{align*}
\end{problem}

\subsection{Heat equation}
\label{sec_heat}
Having the abstract formulations on the continuous, semi-discrete (in time) and 
fully discrete levels at hand, we now proceed and provide two specific realizations. 
First, the heat equation is considered in this subsection and then a combustion problem is introduced in the next subsection.
Let $u:\Sigma\to\mathbb{R}$ be the solution of the heat equation
\begin{align}
\partial_t u - \Delta u &= f\text{ in }\Sigma,\\
u &= 0\text{ on }\partial\Omega\times I\nonumber,\\
u &= u^0\text{ on }\Omega\times\{t=0\},\nonumber
\end{align}
for a given initial condition $u^0\in H$ and right hand side function $f\in L^2(0,T;V^*)$. 
Using the discretization steps as described above, we obtain the following equations
\begin{align}
A_{\text{heat}}(u_{kh})(\varphi_{kh}) &\coloneqq\sum\limits_{m=1}^M\int\limits_{I_m}(\partial_tu_{kh},\varphi_{kh})_H\mathrm{d}t 
+ (\nabla u_{kh},\nabla\varphi_{kh})  \label{eq_heat_operator}\\
&\hphantom{\coloneqq}+ \sum\limits_{m=1}^{M-1}([u_{kh}]_m,\varphi_{kh,m}^+)_H
+(u_{kh,0}^+,\varphi_{kh,0}^+)_H\nonumber\\
F_{\text{heat}}(\varphi_{kh})&\coloneqq (f,\varphi_{kh})+(u^0,\varphi_{kh,0}^+)_H .\label{eqn:rhs_heat_dg_semidiscrete}
\end{align}

\subsection{Combustion}
\label{sec_combustion}
The following coupled nonlinear PDE describes the temperature 
dependent reaction and diffusion of a combustible substance 
without the influence of an additional fluid flow ($v\equiv0$); 
see \cite{Lang2001}.
Therefore, the low Mach number hypothesis holds 
and the fluid flow is not influenced by the reaction and can be ignored. 
The resulting equations for the dimensionless 
temperature $\theta:\Sigma\to\mathbb{R}$ 
and the species concentration $Y:\Sigma\to\mathbb{R}$ are
\begin{align}
 \partial_t \theta -\Delta\theta &= \omega(\theta,Y) \text{ in } \Sigma,\\
 \partial_t Y      -\frac{1}{Le}\Delta Y     &= -\omega(\theta,Y) \text{ in } \Sigma,
\end{align}
with the combustion reaction described by Arrhenius law
\begin{align}
\omega(u) := \omega(\theta,Y) \coloneqq \frac{\beta^2}{2Le}Y\exp\left(\frac{\beta(\theta-1)}{1+\alpha(\theta-1)}\right).
\end{align}
The Arrhenius law is parametrized by the Lewis number $Le>0$, 
the gas expansion $\alpha>0$ and the 
dimensionless activation energy $\beta>0$.
We want to be able to allow all three 
common types of boundary conditions i.e. those of Dirichlet, Neumann and Robin type.
For this we split the boundary $\partial\Omega$ 
into three non-overlapping parts $\Gamma_D$, $\Gamma_N$ and $\Gamma_R$. 
The Dirichlet boundary conditions 
are built into the function spaces, which is the usual approach.
The Neumann and Robin boundary conditions are given by
\begin{align}
 \partial_n \theta &= g_N^\theta\text{ on }\Gamma_N\times I\\
 \partial_n Y      &= g_N^Y     \text{ on }\Gamma_N\times I\\
 a_R^\theta \theta + b_R^\theta\partial_n \theta &= g_R^\theta\text{ on }\Gamma_R\times I\\
 a_R^Y      Y      + b_R^Y     \partial_n Y      &= g_R^Y     \text{ on }\Gamma_R\times I.
\end{align}
It remains to state the initial conditions:
\begin{align*}
\theta = \theta^0 \quad\text{on } \Omega\times \{t=0\},\\
Y = Y^0 \quad\text{on } \Omega\times \{t=0\}.
\end{align*}

By following the typical steps for the derivation of a weak formulation, 
integration by parts in space, and subsequent summation, we obtain
\begin{align}
 &(\partial_t \theta,\varphi^\theta ) + (\nabla\theta,\nabla\varphi^\theta) 
 - \int\limits_{I\times\partial\Omega} \partial_n \theta \varphi^\theta \mathrm{d}s\mathrm{d}t
 - (\omega(\theta,Y),\varphi^\theta) \nonumber \\
+&(\partial_t Y,\varphi^Y ) + (\nabla Y,\nabla\varphi^Y) - \int\limits_{I\times\partial\Omega} \partial_n Y \varphi^Y \mathrm{d}s\mathrm{d}t 
 + (\omega(\theta,Y),\varphi^Y) =0. \nonumber
\end{align}
Splitting the boundary integrals and considering 
homogeneous Dirichlet conditions on some parts, i.e., $\varphi^\theta|_{\Gamma_D}=\varphi^Y|_{\Gamma_D}=0$, we get
\begin{align}
 \int\limits_{I\times\partial\Omega} \partial_n \theta \mathrm{d}s\mathrm{d}t &= 
 \int\limits_{I\times\Gamma_N} g_N^\theta\varphi^\theta \mathrm{d}s\mathrm{d}t 
  + \int\limits_{I\times\Gamma_R} \frac{g_R^\theta}{b_R^\theta}\varphi^\theta
  -\frac{a_R^\theta}{b_R^\theta}\theta\varphi^\theta \mathrm{d}s\mathrm{d}t,\\
    \int\limits_{I\times\partial\Omega} \partial_n Y \mathrm{d}s\mathrm{d}t &= 
    \int\limits_{I\times\Gamma_N} g_N^Y\varphi^Y \mathrm{d}s \mathrm{d}t
    + \int\limits_{I\times\Gamma_R} \frac{g_R^Y}{b_R^Y}\varphi^Y-\frac{a_R^Y}{b_R^Y}Y\varphi^Y \mathrm{d}s\mathrm{d}t.
\end{align}

By introducing jump terms as described earlier, we obtain 
the following semi-linear and linear forms, respectively:
\begin{align}
 A_{\text{comb.}}(u_{kh})(\varphi_{kh})\coloneqq&
   \sum\limits_{m=1}^M\int\limits_{I_m} (\partial_t\theta_{kh},\varphi_{kh}^\theta)_H 
   +(\nabla\theta_{kh},\nabla\varphi_{kh}^\theta)_H \; \mathrm{d}t
   +\sum\limits_{m=1}^{M-1} ([\theta_{kh}]_{m},\varphi_{kh,m}^{\theta,+})_H \nonumber \\
   +&
   \sum\limits_{m=1}^M\int\limits_{I_m} (\partial_tY_{kh},\varphi_{kh}^Y)_H 
   +(\nabla Y_{kh},\nabla\varphi_{kh}^Y)_H \;\mathrm{d}t
   + \sum\limits_{m=1}^{M-1} ([Y_{kh}]_{m},\varphi_{kh,m}^{Y,+})_H \nonumber\\
    +&\int\limits_{I\times\Gamma_R} \frac{a_R^\theta}{b_R^\theta}\theta\varphi^\theta + \frac{a_R^Y}{b_R^Y}Y\varphi^Y \mathrm{d}s\mathrm{d}t
     + (\omega(u_{kh}),\varphi_{kh}^Y-\varphi_{kh}^\theta)\nonumber\\
     +&(\theta_{kh,0}^{+},\varphi_{kh,0}^{\theta,+})_H+(Y_{kh,0}^{+},\varphi_{kh,0}^{Y,+})_H, 
\label{eq_combustion_weak_form_lhs}\\
 F_{\text{comb.}}(\varphi_{kh})\coloneqq&\int\limits_{I\times\Gamma_N} g_N^\theta\varphi^\theta + g_N^Y\varphi^Y \mathrm{d}s\mathrm{d}t
 +\int\limits_{I\times\Gamma_R} \frac{g_R^\theta}{b_R^\theta}\varphi^\theta +\frac{g_R^Y}{b_R^Y}\varphi^Y \mathrm{d}s\mathrm{d}t \nonumber \\
 +&(\theta^{0},\varphi_{kh,0}^{\theta,+})_H+(Y^{0},\varphi_{kh,0}^{Y,+})_H
, \label{eq_combustion_weak_form_rhs}
\end{align}
with $u_{kh} = (\theta_{kh},Y_{kh})$ and $\varphi_{kh} = (\varphi_{kh}^\theta,\varphi_{kh}^Y)$.

\subsection{General formulation of parabolic problems}
A general parabolic weak formulation that includes the previous problem statements 
can be stated by
\begin{equation}
\begin{aligned}
 A_{\text{gen.}}(u_{kh})(\varphi_{kh})\coloneqq& \sum\limits_{m=1}^M 
 \int\limits_{I_m} (\partial_t u_{kh},\varphi_{kh})_H \mathrm{d}t+ a(u_{kh},\varphi_{kh})
 +\sum\limits_{m=1}^{M-1}([u_{kh}]_m,\varphi_{kh,m}^+)_H + (u_{kh,0}^+,\varphi_{kh,0}^+),\\
 F_{\text{gen.}}(\varphi_{kh})\coloneqq& (f,\varphi_{kh}) + (u^0,\varphi_{kh,0}^+)_H \label{eq_primal_problem}
\end{aligned}
\end{equation}
with an elliptic operator $a(u,\varphi)\coloneqq\int\limits_0^T \bar{a}(u(t),\varphi(t))\mathrm{d}t$.
Then, (non-)linearity of $A$ solely depends on the (non-)linearity of $\bar{a}$.
\begin{remark}
We notice that additional terms due to Neumann or Robin boundary conditions would appear inside $\bar{a}(u(t),\varphi(t))$ and/or $F(\cdot)$.
\end{remark}

\subsection{Numerical solution}
In the algorithmic realization, we notice that the choice of trial
and test spaces allow for a decoupling of the temporal discretization into slabs,
i.e. slices of the space-time cylinder. 
In the simplest case a slice just encompasses a single temporal interval, 
resulting in a sequential time-stepping scheme due to 
the dG(r) test functions,
and therefore effectively yielding classical time discretization schemes.
For dG(0) an implicit Euler-type scheme is recovered. 
At each time slab, the spatial problems are solved as described in the following.

For the basic implementation of the (linear) heat equation with a classical DWR error estimator, 
we refer to the \textit{dwr-diffusion} package \cite{kocher.etal2019} 
and the solvers implemented therein.
There, the sparse direct solver UMFPACK \cite{davis2004} 
is used for the linear equation systems.

For the nonlinear combustion equation we employ a classical Newton-type solver 
as briefly described in the following.
In the space-time setting we have to solve
\begin{equation*}
A (u_{kh})(\varphi_{kh} ) = 0
\quad \forall \varphi_{kh}\in \widetilde{X}^{r,s}_{kh}(\mathcal{T}_k,\mathcal{T}_h^{1,\dots,M}) ,
\end{equation*}
where $u_{kh}$ is the complete space-time solution over all intervals.
Given an initial guess $u_{kh}^{0}$, 
find the update $\delta u\in \widetilde{X}_{kh}^{r,s}(\mathcal{T}_k,\mathcal{T}_h^{1,\dots,M})$ of the
linear defect-correction problem for $j=0,1,2, \ldots$
\begin{align}
\label{dis_newton_defect}
A_u ' (u_{kh}^{j})(\delta u , \varphi_{kh} ) &=
-A (u_{kh}^{j})(\varphi_{kh} ),
\nonumber \\
u_{kh}^{j+1} &= u_{kh}^{j} + \alpha \delta u, \quad \alpha\in (0,1].
\end{align}
\begin{remark}
 As we have discontinuous test functions the nonlinear problem can be decoupled into one subproblem per time interval.
 Then, we obtain a time-stepping scheme with one Newton loop per interval.
\end{remark}
The defect-correction problems are solved using
the parallel sparse direct solver MUMPS \cite{MUMPS}.
The Jacobian $A'(\cdot)(\cdot,\cdot)$ is derived 
by using analytical expressions, e.g., \cite[Chapter 13]{Wi22_num_pde}, in order
to maintain superlinear (or quadratic) convergence. 
However, for $j>0$ reassembly of the matrix is omitted if 
the relative reduction of the residual is below a certain threshold.
This simplified Newton approach saves a lot of computational 
time as the factorization only needs to be done after a reassembly step.
This approach also benefits from 
preconditioned Krylow methods as the preconditioner is
only recomputed after reassembly.
The step size $\alpha$ is is determined by a damping line-search 
after solving the defect-correction problem.

\section{The dual weighted residual method in a space-time setting} 	
\label{Section: space-time DWR}
In this section, we review the general ideas of the DWR method.
We then derive joint and split error identities and corresponding error estimators.

\subsection{Error representation and estimation}
Let $J:X\to\mathbb{R}$ be some goal functional representing 
some quantity of interest (QoI). The general form reads
\begin{equation}
\label{eq_goal_functional}
J(u) = \int_{T_1}^{T_2} u(t)\, dt + u(T),
\end{equation}
where $T_1,T_2\in I$, e.g., $T_1 = 0$ and $T_2 = T$ and $T$ is the end time value.
We are interested in the discretization error $J(u) - J(u_{kh})$ and more specifically 
to minimize this error for a reasonable computational cost:
\[
\min J(u) - J(u_{kh})
\]
This becomes a constrained optimization problem since the solutions $u\in X$ 
and $u_{kh} \in \tilde X_{k,h}^{r,s}$ are obtained as PDE solutions from 
\eqref{eq_weak_form_abstract} and \eqref{eq_weak_form_abstract_fully_discrete},
respectively. These PDE problem statements are seen 
as constraint in terms of the optimization problem. Consequently,
for a given goal functional $J(u)$ we want to solve the following optimization problem 
\begin{align}
\label{eq_opt}
\min_{u\in X(I,V)}J(u) \quad
s.\,t.\; A(u)(\varphi)=F(\varphi)\quad\forall\varphi\in X(I,V).
\end{align}
This problem setting is the same in \cite{becker.rannacher2001}[Section 2.2, (2.13)]. Clearly, after the discretization, 
we can then measure the discretization error $J(u) - J(u_{kh})$. We notice 
that from an optimization viewpoint $J(u_{kh})$ is a constant in \eqref{eq_opt} and 
therefore implicitly contained therein modulo the constant shift.\\
We apply the method of Lagrange multipliers for this constrained optimization problem 
(see e.g., \cite{becker.rannacher2001}) and introduce the dual variable $z\in X(I,V)$.
To account for the discontinuities in the primal problem we define a discontinuous Lagrange functional
$\widetilde{\mathcal{L}}(u,z)\colon \widetilde{X}(\mathcal{T}_k,V)\times\widetilde{X}(\mathcal{T}_k,V)\to\mathbb{R}$, such that
\begin{equation}
 \widetilde{\mathcal{L}}(u,z)\coloneqq J(u)-A(u)(z)+F(z).
\end{equation}
Note that $A(u)(z)$ and $F(z)$ contain the jump terms and weakly imposed initial conditions.\\
For stationary points $(u,z)\in X(I,V)\times X(I,V)$ all jump terms vanish and the initial conditions are met exactly, 
such that $\widetilde{\mathcal{L}}(u,z)$ is consistent with the continuous functional $\mathcal{L}(u,z)$.
The first order optimality conditions yield the original primal problem $A(u)(\varphi)=F(\varphi)$ as well as the adjoint problem:\\
Find $z\in \widetilde{X}(\mathcal{T}_k,V)$ such that
\begin{align}
\label{eq_adjoint}
A_u'(u)(\psi,z) = J_u'(u)(\psi)\;\forall\psi\in \widetilde{X}(\mathcal{T}_k,V).
\end{align}

\begin{remark}
Note that the test and trial functions are switched in \eqref{eq_adjoint}.
Accordingly, the temporal derivative is now applied to the test function $\psi$. To rectify this, we apply integration by parts 
to the corresponding scalar product, obtaining:
\begin{equation}
 \begin{aligned}
&\sum\limits_{m=1}^M\int_{I_m}(\psi,-\partial_t z)+
\bar{a}'_u(u)(\psi,z)\;\mathrm{d}t+
\sum\limits_{m=1}^{M-1}(\psi_m^-,-[z]_m)_H\\ &+(\psi(T),z(T))_H 
= (\psi(T),z^M)_H +J_u'(u)(\psi)
\quad\forall\psi\in \widetilde{X}(\mathcal{T}_k,V).
\end{aligned}
\end{equation}
The negative sign of the temporal derivative means the adjoint problem has to be solved backwards in time, with an 
initial value $z^M$ depending on the goal functional.\\
For functionals evaluated on the whole temporal domain $z^M=0$ holds and for functionals defined only at $T$ 
$J_u'(u(T),\psi(T))$ can be reinterpreted as an initial value $z^M$ for details see e.g. \cite{meidner2007}.
\end{remark}
For linear PDEs and linear goal functionals we obtain the exact error representations (see \cite{becker.rannacher2001}):
\begin{align}
J(u)-J(u_{kh}) &= F(z-z_{kh})-A(u_{kh},z-z_{kh})&&\text{(primal error)}\label{eqn:primal_error_rep}\\
&= J(u-u_{kh})-A(u-u_{kh},z_{kh})&&\text{(adjoint error)}\label{eqn:adjoint_error_rep}.
\end{align}  
In the following, we focus on the primal error representation. As we can see we would need both the exact dual solution $z$
and the discrete dual solution $z_{kh}$. 
As this is infeasible for complicated problems we use a 
discrete solution of higher order for $z$.
In the equal low order approach both primal and adjoint problem 
are discretized using the same low order elements.
The higher order adjoint solution is then obtained by a patch wise reconstruction. 
This reconstruction is described in detail in \cite{schmich.vexler2008}.
In the mixed order approach the adjoint problem is discretized by higher order elements and 
the solution is used as the representation of the exact solution. 
The fully discrete solution is then obtained by interpolation into the lower order space.
It is also possible to use different approaches in time and space e.g. discretizing
the primal problem with $cG(1)dG(0)$ and the adjoint problem with $cG(2)dG(0)$.

Additionally, \eqref{eqn:primal_error_rep} can be split into a temporal and a spatial part by introducing 
the semidiscrete adjoint solution $z_k$ such that
\begin{align}
J(u)-J(u_{kh}) &= J(u)-J(u_k)+J(u_k)-J(u_{kh}),
\end{align}
where temporal and spatial errors are given by, respectively,
\begin{align}
J(u)-J(u_k) &= F(z-z_k)-A(u_{kh},z-z_k),\label{eqn:primal_temp_error_rep}\\
J(u_k)-J(u_{kh}) &= F(z_k-z_{kh})-A(u_{kh},z_k-z_{kh}).\label{eqn:primal_spat_error_rep}
\end{align}
That way \eqref{eqn:primal_temp_error_rep} can be used for temporal refinement and \eqref{eqn:primal_spat_error_rep} for spatial refinement.

\subsection{DWR for nonlinear time dependent problems }
\label{sec_dwr_nonlinear}

\subsubsection{Adjoint problem statements}

For $A_{\text{gen.}}$ the left hand side of the adjoint problem \eqref{eq_adjoint} in the space-time context explicitly reads as
\begin{equation}
 \begin{aligned}
A'_{u,\text{gen.}}(u)(\psi,z) \coloneqq &
\sum\limits_{m=1}^M \int_{I_m}(\psi,-\partial_t z_{kh})+ \bar{a}'_u(u(t))(\psi(t),z(t))\;\mathrm{d}t \\
&+\sum\limits_{m=1}^{M-1}(\psi_m^-,[z_{kh}]_m)_H+ (\psi(T),z_{kh}(T))_H.
\end{aligned} 
\end{equation}
As an example, for the combustion problem 
described in the Section \ref{sec_combustion} we obtain $F = F_{\text{comb.}}$ 
as well as the operator of the semi-linear form and its directional derivative, respectively,
\begin{align}
 \bar{a}(u(t),\varphi(t)) &= (\nabla\theta,\nabla\varphi^\theta)_H + (\nabla Y,\nabla\varphi^Y)_H + (\omega(u),\varphi^Y-\varphi^\theta)_H\\
 &\hphantom{=}+\int\limits_{\Gamma_R} \frac{a_R^\theta}{b_R^\theta}\theta\varphi^\theta + \frac{a_R^Y}{b_R^Y}Y\varphi^Y \mathrm{d}s\nonumber,\\
 \bar{a}_u'(u(t))(\psi,z) &= (\nabla\psi^\theta,\nabla z^\theta)_H + (\nabla\psi^Y,\nabla z^Y)_H + (\omega_\theta'(u)(\psi^\theta),z^Y-z^\theta)_H\\
 &\hphantom{=}+\int\limits_{\Gamma_R} \frac{a_R^\theta}{b_R^\theta}\psi^\theta z^\theta + \frac{a_R^Y}{b_R^Y}\psi^Y z^Y \mathrm{d}s+(\omega_Y'(u)(\psi^Y),z^Y-z^\theta)_H\nonumber.
\end{align}
Therein, $\omega'(u)(\psi)$ is the directional derivative of $\omega(u)$ into 
the direction $\psi$.

\subsubsection{Goal-oriented error representations}

As there are two ways to separate the full estimator 
into parts, we want to use a clear and precise terminology
to distinguish between those two. 

Firstly, we can separate by the problem residuals we compute. 
This gives the \textit{primal} and \textit{adjoint/dual} estimators. 
Oftentimes, the average of those two parts 
is also called \textit{mixed} estimator instead 
of \textit{full} estimator. 

Secondly, we can split the difference between 
the exact and the fully discrete solution by introducing the time-discrete solution. 
Not doing so we will call the resulting estimator 
the \textit{joint} estimator.
If we calculate both the temporal and spatial error estimator, we will call the sum of both the \textit{split} error estimator. 
As both seperations can be done simultaneously, combined expressions 
like split primal estimator or joint dual estimator are possible.

For nonlinear problems we obtain the following error representation \cite{becker.rannacher2001}:
\begin{theorem}[Joint error identity]
\label{theorem_joint_error}
Let the primal problem and adjoint problem be given. 
Let $(u,z)\in \widetilde{X}(\mathcal{T}_k,V)\times \widetilde{X}(\mathcal{T}_k,V)$, 
$(u_k,z_k)\in \widetilde{X}_k^r(\mathcal{T}_k,V) \times \widetilde{X}_k^r(\mathcal{T}_k,V)$ and 
$(u_{kh},z_{kh})\in \widetilde X_{k,h}^{r,s}(\mathcal{T}_k,\mathcal{T}_h^{1,\dots,M}) 
\times \widetilde X_{k,h}^{r,s}(\mathcal{T}_k,\mathcal{T}_h^{1,\dots,M})$.
Then, we have the space-time joint error identity
\begin{align}
\label{eq_DWR_theorem}
J(u)-J(u_{kh}) &= \frac12 \rho(u_{kh})(z-z_{kh}) + 
\frac12 \rho^*(u_{kh},z_{kh})(u-u_{kh})+ \mathcal{R}_{kh}, 
\end{align}
with the primal error estimator $\rho$ and the adjoint error estimator $\rho^*$ 
\begin{align*}
\rho(u_{kh})(z-z_{kh}) &:= F(z-z_{kh}) - A(u_{kh},z-z_{kh}),\\
\rho^*(u_{kh},z_{kh})(u-u_{kh}) &:= J'(u_{kh})(u-u_{kh}) - A'_u(u_{kh})(u-u_{kh},z_{kh}),
\end{align*}
as well as a remainder term $\mathcal{R}_{kh}$ of higher order.
\end{theorem}
\begin{proof}
With $\widetilde X_{k,h}^{r,s}(\mathcal{T}_k,\mathcal{T}_h^{1,\dots,M})
\subset\widetilde{X}(\mathcal{T}_k,V)$ and $J(u_{kh})=\widetilde{\mathcal{L}}(u_{kh},z_{kh})$ 
the assumptions of \cite{schmich.vexler2008}[Proposition 3.1] hold, proving the representation.
\end{proof}

\begin{theorem}[Split error identity]
\label{theorem_split_error}
With the previous assumptions, we have the split error identity
\begin{align*}
J(u)-J(u_{kh}) = (J(u) - J(u_k)) + (J(u_k) - J(u_{kh})),
\end{align*}
with 
\begin{align*}
J(u) - J(u_k) &= \frac12 \rho(u_{k})(z-z_{k}) 
+ \frac12 \rho^*(u_{k},z_{k})(u-u_{k}) + R_k,\\
J(u_k) - J(u_{hk}) &= \frac12 \rho(u_{kh})(z_k-z_{kh}) 
+ \frac12 \rho^*(u_{kh},z_{kh})(u_k-u_{kh}) + R_h.
\end{align*}
\end{theorem}
\begin{proof}
The proof follows the same ideas as before but with 
$\widetilde X_{k,h}^{r,s}(\mathcal{T}_k,\mathcal{T}_h^{1,\dots,M})
\subset\widetilde{X}_k^r(\mathcal{T}_k,V)
\subset\widetilde{X}(\mathcal{T}_k,V)$ 
as well as  $J(u_{kh})=\widetilde{\mathcal{L}}_k(u_{kh},z_{kh})$ 
and $J(u_k) = \widetilde{\mathcal{L}}(u_k,z_k)$.
\end{proof}

\subsubsection{Error estimators}
From the previous error identities, we obtain 
error estimators in four variants.
First, we have the full error estimator
\begin{equation}
\label{eq_DWR_theorem_est}
\eta := \frac12 \rho(u_{kh})(z-z_{kh}) + \frac12 \rho^*(u_{kh},z_{kh})(u-u_{kh})+ \mathcal{R}_{kh}.
\end{equation}
However, the unknown solutions $u$ and $z$ still enter. 
This is already an error estimator, because for cases where 
$u$ and $z$ are known, we can already estimate (discretization) errors in goal 
functionals. Of course, for most problems in practice, this first version does 
not play a role.\\
To this end,
higher-order approximations $\widetilde u\in \widetilde{X}$ and $\widetilde z\in \widetilde{X}$ 
are introduced \cite{becker.rannacher2001}. Examples 
of such approximations are $\widetilde u := u_{kh}^{r+1,s+1}$ and $\widetilde z := z_{kh}^{r+1,s+1}$
such that we obtain the computable error estimator
\begin{align}
\label{eq_comp_est}
\eta^{(r+1,s+1)} &:= \frac12 \rho(u_{kh})(z_{kh}^{r+1,s+1}-z_{kh}) 
+ \frac12 \rho^*(u_{kh},z_{kh})(u_{kh}^{r+1,s+1}-u_{kh}) + \mathcal{R}_{kh}.
\end{align}
If the remainder term is omitted (which is indeed usually done in practice) 
we obtain the practical error estimator
\begin{align}
\label{eq_prac_est}
\eta_h^{(r+1,s+1)} &:= \frac12 \rho(u_{kh})(z_{kh}^{r+1,s+1}-z_{kh}) 
+ \frac12 \rho^*(u_{kh},z_{kh})(u_{kh}^{r+1,s+1}-u_{kh}).
\end{align}
Finally, we also introduce the primal-based error estimator:
\begin{align}
\label{eq_primal_error_estimator}
\eta_{prim}^{(r+1,s+1)} := \rho(u_{kh})(z_{kh}^{r+1,s+1} - z_{kh}).
\end{align} 
As we also need higher order information of the primal problem for \eqref{eq_comp_est} and \eqref{eq_prac_est}
to calculate the adjoint estimator we now have three possible approaches.
In addition to the two previous discretization approaches the
equal high order approach uses a higher order element discretization for both the primal and adjoint problem.
Then, interpolation into a lower order element space yields both $u_{kh}$ and $z_{kh}$.
However, inserting the interpolated $u_{kh}$ into the goal functional yields worse results
compared to a native low order solution, so ideally the primal problem should also be solved in low order to calculate the
functional values. Various algorithmic realizations 
with corresponding theoretical results, and performance 
analyses for stationary problems were recently established in \cite{EndtLaWi21_smart}.

\section{Error localization}
\label{Section: Localization}
In this section, we address our key development, namely the construction 
of a space-time partition-of-unity (PU) localization of goal-oriented 
a posteriori error estimators. In the published literature
as mentioned in the introduction, so far only stationary cases have been addressed with 
the PU localization. Here, we first extend the idea to time-dependent problems.
Since we want to use the DWR error estimator for grid refinement, 
we need to split the estimator
into element- or DoF-wise error contributions.
Three known approaches are the classical integration 
by parts \cite{becker.rannacher2001,bangerth.rannacher2003}, 
a variational filtering operator over patches of elements \cite{braack.ern2003}
and a variational partition-of-unity localization \cite{richter.wick2015}.
For stationary problems, the effectivity of these 
localizations was established and numerically substantiated 
in \cite{richter.wick2015}.

First, in Section \ref{sec_PU}, we exemplarily derive the variational 
partition of unity approach
for our space-time estimator of the heat equation.
Second, Section \ref{sec:eval_pu_dwr} focuses 
on details of the actual evaluation 
including the needed interpolation operations.  
Next, we list in Section \ref{sec:error_indi} the resulting 
error indicators, finally followed by the adaptive algorithms 
designed in Section \ref{sec_adaptive_algo}.

\subsection{The partition-of-unity approach for the heat equation}
\label{sec_PU}
In this key section, we extend the ideas from \cite{richter.wick2015}
and apply a partition-of-unity (PU) localization to a space-time 
error estimator.
To this end, we first design the PU space. 
The simplest choice is $V_{PU} = \widetilde{X}_{k,h}^{0,1}$, i.e.
a $cG(1)dG(0)$ discretization.
Effectively, this yields one spatial partition of unity 
$(\chi_{i,m})_{i=1}^{\#DoFs(\mathcal{T}_h^m)}\in Q_1(\mathcal{T}_h^m)$ 
per time interval $I_m$ for $m=1,\ldots, M$.
As this is a Lagrangian finite element, we have immediately
\begin{proposition}
\label{prop_PU}
For a function $\chi\in V_{PU}$, it holds
\begin{align}
\label{eq_PU}
 \sum_{m=1}^M \sum_{i=1}^{\#DoFs(\mathcal{T}_h^m)} \chi_{i,m} \equiv1.
\end{align}
\end{proposition}
\begin{proof}
Follows immediately from the properties of the finite element functions.
\end{proof}
\begin{remark}
Another common choice for the PU space is $V_{PU} = X_{kh}^{1,1}$, i.e. a $cG(1)cG(1)$ discretization, 
for example in native $d+1$-dimensional discretizations \cite{endtmayer2023goaloriented}.
In general this ensures a coupling between neighboring temporal elements to address the 
problem shown in \cite{CaVe99}. However, for discontinuous Galerkin discretizations the dominating edge residuals, i.e. jump terms, are explicitly 
included in the estimator.
\end{remark}
\begin{remark}
Clearly, using a $dG$ discretization in time yields a natural decoupling for which the space-time PU reduces effectively to a PU in space. Nonetheless, we formulated  our concepts using a space-time PU as the methodology applies (see again 
\cite{endtmayer2023goaloriented}) to a larger class of problems. As our work is one of the first into this direction in the literature, our aim is to provide 
the full methodology in order to have a starting point for further future work.
\end{remark}
In the following, $\tilde{u}\in \widetilde{X}$ and $\tilde{z}\in \widetilde{X}$ denote approximations of the exact solutions.
In principle the joint and split estimators only differ in the interpolation 
difference i.e $\tilde{z}-z_{kh}$ or
$\tilde{z}-z_k$ and $z_k-z_{kh}$ respectively.

For the joint estimator the local contributions are summed over 
all DoFs of a fixed interval to obtain the corresponding temporal estimator. 
Subsequent summation over all time intervals yields the global estimator.
For the split estimators the spatial estimator is summed over all space-time 
DoFs and the temporal estimator is summed over all time intervals.
The total error estimator is then the sum of these two error parts. 

\begin{proposition}[Primal joint error estimator for the heat equation]
\label{prop_primal_joint_estimator_heat}
For the space-time formulation of the heat equation, we have the following 
a posteriori joint error estimator with partition-of-unity localization:
\begin{align}
|J(u) - J(u_{kh}) | \leq& |\eta_{\text{joint}}| \coloneqq 
\left| \sum\limits_m \eta_{kh}^m \right|, 
\qquad \text{with }\;
\eta_{kh}^m\coloneqq\sum\limits_{i\in\mathcal{T}_h^m} \eta_{kh}^{i,m},
\end{align}
with the error indicators
\begin{equation}
\begin{aligned}
\eta_{kh}^{i,m} &:= \int\limits_{I_m} (f,(\tilde{z}-z_{kh})\chi_{i,m})_H \;\mathrm{d}t
       -\int\limits_{I_m}(\nabla u_{kh}, \nabla ((\tilde{z}-z_{kh})\chi_{i,m}))_H\;\mathrm{d}t\\
       &\quad 
       -\int\limits_{I_m}(\partial_t u_{kh}, (\tilde{z}-z_{kh})\chi_{i,m})_H \;\mathrm{d}t
       -([u_{kh}]_{m-1},(\tilde{z}^+(t_{m-1})-z_{kh}^+(t_{m-1}))\chi_{i,m})_H.
\end{aligned}
\end{equation}
\end{proposition}
\begin{proof}
We start from Theorem \ref{theorem_joint_error} with 
\[
J(u)-J(u_{kh}) = \frac12 \rho(u_{kh})(z-z_{kh}) + 
\frac12 \rho^*(u_{kh},z_{kh})(u-u_{kh})+ \mathcal{R}_{kh}, 
\]
which yields 
\[
|J(u)-J(u_{kh})| \leq |\frac12 \rho(u_{kh})(z-z_{kh}) + 
\frac12 \rho^*(u_{kh},z_{kh})(u-u_{kh}) + \mathcal{R}_{kh}|.
\]
Considering the primal part only (see \eqref{eq_primal_error_estimator}), gives us
\[
|J(u)-J(u_{kh})| \leq |\eta_{joint}| := |\rho(u_{kh})(z-z_{kh})|.
\]
Inserting the PU \eqref{eq_PU} yields
\[
|J(u)-J(u_{kh})| \leq |\eta_{joint}| := \left|\sum_{m=1}^M \sum_{i=1}^{\#DoFs(\mathcal{T}_h^m)} 
\rho(u_{kh})\big((z-z_{kh})\chi_{i,m}\big) \right|.
\]
Then, employing the definition of the primal residual leads to
\[
|J(u)-J(u_{kh})| \leq |\eta_{joint}| := \left|\sum_{m=1}^M \sum_{i=1}^{\#DoFs(\mathcal{T}_h^m)} 
F\big( (z-z_{kh})\chi_{i,m} \big) - A(u_{kh})\big( (z-z_{kh})\chi_{i,m} \big) \right|.
\]
Here, we employ the left hand side and right hand side of the heat equation, 
namely \eqref{eq_heat_operator} and \eqref{eqn:rhs_heat_dg_semidiscrete}, respectively,
by replacing the test function $\varphi_{kh}$ by the PU-weighted adjoint sensitivity 
measure~$(z-z_{kh})\chi_{i,m}$. 
Finally, the (unknown) solution $z$ is approximated by some higher order 
representation $\tilde z$ from which we obtain the assertion.
\end{proof}
\begin{definition}[Effectivity and indicator indices]
We notice that the effecitivity index is defined as 
\[
I_{eff} := \frac{|J(u)-J(u_{kh})|}{|\eta_{joint}|}.
\]
Applying the triangle inequality on $\eta_{joint}$ 
yields a more strict criterion, i.e., here
\[
|\eta_{joint}| = \left| \sum\limits_m \eta_{kh}^m \right| \leq \sum\limits_m \sum\limits_{i\in\mathcal{T}_h^m} |\eta_{kh}^{i,m}|,
\]
from which the so-called indicator index 
\[
I_{ind} := \frac{|J(u)-J(u_{kh})|}{\sum\limits_m \sum\limits_{i\in\mathcal{T}_h^m} |\eta_{kh}^{i,m}|}
\] 
can be defined.
\end{definition}

\begin{remark}
Note that we only use the temporal part for marking 
time steps and calculating the global estimator.
Since the indicators for each time step are obtained by summing over all elements in said time step
the spatial PU $\chi_{i,m}$ effectively cancels due to the PU property. 
As a consequence, the spatial PU can be omitted directly in the computation of the temporal indicators.
\end{remark}

\begin{proposition}[Primal split error estimator for the heat equation]
\label{prop_primal_split_heat}
For the space-time formulation of the heat equation, we have the following 
a posteriori split error estimator with partition-of-unity localization:
\begin{align}
|J(u) - J(u_{kh})| 
\leq  |\eta_{\text{split}}| \coloneqq \left|\sum\limits_m \left( \eta_k^m 
+ \sum_{i\in\mathcal{T}_h^m} \eta_{h}^{i,m} \right) \right|,
\end{align}
with the temporal error indicators
\begin{equation}
\begin{aligned}
\eta_k^m &:=\int\limits_{I_m} \Big( (f,\tilde{z}-z_{k})_H 
       -(\partial_t u_{kh}, \tilde{z}-z_{k})_H 
       -(\nabla u_{kh}, \nabla (\tilde{z}-z_{k}))_H \Big) \;\mathrm{d}t\\
       &\quad -([u_{kh}]_{m-1},\tilde{z}^+(t_{m-1})-z_{k}^+(t_{m-1}))_H,
\end{aligned}
\end{equation}
and the spatial error indicators
\begin{equation}
\begin{aligned}
\eta_{h}^{i,m} &:=  \int\limits_{I_m}(f,(z_k-z_{kh})\chi_{i,m})_H \;\mathrm{d}t
       -\int\limits_{I_m}(\nabla u_{kh}, \nabla ((z_k-z_{kh})\chi_{i,m}))_H\;\mathrm{d}t \\
       &\quad-\int\limits_{I_m}(\partial_t u_{kh}, (z_k-z_{kh})\chi_{i,m})_H \;\mathrm{d}t
      -([u_{kh}]_{m-1},(z_k^+(t_{m-1})-z_{kh}^+(t_{m-1}))\chi_{i,m})_H.
\end{aligned}
\end{equation}
\end{proposition}
\begin{proof}
For the primal split error estimator we start with Theorem \eqref{theorem_split_error}.
Then, we proceed for both parts with the primal estimator:
\begin{align*}
|J(u) - J(u_k)| &\leq |\rho(u_{k})(z-z_{k})|,\\
|J(u_k) - J(u_{hk})| &\leq |\rho(u_{kh})(z_k-z_{kh})|.
\end{align*}
Combining the left hand sides to the full error estimator and utilizing the PU in a similar 
way as the proof of Proposition \eqref{prop_primal_joint_estimator_heat} yields
\begin{align*}
|J(u) - J(u_{kh})| &\leq |\rho(u_{k})(z-z_{k}) + \rho(u_{kh})(z_k-z_{kh})|\\
&\leq \left| \sum\limits_m \left( \eta_k^m 
+ \sum_{i\in\mathcal{T}_h^m} \eta_{h}^{i,m} \right) \right|\\ 
&=: |\eta_{\text{split}}|.
\end{align*}
Here, the error indicators are obtained as 
in the proof of Proposition \eqref{prop_primal_joint_estimator_heat}, 
which yields the assertion. 
\end{proof}

\begin{remark}
Note that the basis functions are globally defined, 
so the DoF-wise errors contain an implicit sum over all elements in practice.
However, calculating element based estimators by constraining the 
spatial integrals to each element $K$ yields 
the unlocalized estimator, as the sum over all $\chi_{i,m}$ on a 
single element is $1$ and effectively cancels out.
To use well-known element based marking strategies 
the DoF-estimators have to be calculated globally.
Afterwards element wise estimators can be calculated by summing 
all estimators belonging to the DoFs of the corresponding element:
\begin{align}
 \eta_K^m = \sum\limits_{i\in K} \eta_{\bullet}^{i,m},
\end{align}  
where $\bullet$ stands for $h$ or $kh$.
\end{remark}

\begin{proposition}[Adjoint joint error estimator for the heat equation]
For the space-time formulation of the heat equation, we have the following 
a posteriori joint error estimator with partition-of-unity localization:
\begin{align}
|J(u) - J(u_{kh}) | \leq& |\eta_{\text{joint}}^*| \coloneqq |\sum\limits_m \eta_{kh}^{m,*}|, 
\quad \text{with }\;
\eta_{kh}^{m,*}\coloneqq\sum\limits_{i\in\mathcal{T}_h^m} \eta_{kh}^{i,m,*},
\end{align}
with the error indicators
\begin{equation}
\begin{aligned}
\eta_{kh}^{i,m,*} &:= \int\limits_{I_m}J_u'(u_{kh})((\tilde{u}-u_{kh})\chi_{i,m}) \;\mathrm{d}t
       -\int\limits_{I_m}(\nabla ((\tilde{u}-u_{kh})\chi_{i,m}),\nabla z_{kh})_H\;\mathrm{d}t\\
       &\quad +\int\limits_{I_m}((\tilde{u}-u_{kh})\chi_{i,m},\partial_t z_{kh})_H \;\mathrm{d}t
       +((\tilde{u}^-(t_m)-u_{kh}^-(t_m))\chi_{i,m},[z_{kh}]_m)_H.
\end{aligned}
\end{equation}
\end{proposition}
\begin{proof}
We start from Theorem \ref{theorem_joint_error} with 
\[
J(u)-J(u_{kh}) = \frac12 \rho(u_{kh})(z-z_{kh}) + 
\frac12 \rho^*(u_{kh},z_{kh})(u-u_{kh})+ \mathcal{R}_{kh}, 
\]
which yields 
\[
|J(u)-J(u_{kh})| \leq |\frac12 \rho(u_{kh})(z-z_{kh}) + 
\frac12 \rho^*(u_{kh},z_{kh})(u-u_{kh})+ \mathcal{R}_{kh}|.
\]
Considering the adjoint part only, gives us
\[
|J(u)-J(u_{kh})| \leq |\eta_{joint}^*| := |\rho^*(u_{kh},z_{kh})(u-u_{kh})|
= |J'(u_{kh})(u-u_{kh}) - A'_u(u_{kh})(u-u_{kh},z_{kh})|.
\]
Utilizing the PU as in the proof of Proposition \eqref{prop_primal_joint_estimator_heat}
and the respective definition of $J$ in \eqref{eq_goal_functional} and the adjoint 
$A'_u$ of the heat equation by again approximating $u$ by some higher-order 
approximation $\tilde u$ yields the assertion.
\end{proof}

\begin{proposition}[Adjoint split error estimator for the heat equation]
\label{prop_adjoint_spit_estimator_heat}
For the space-time formulation of the heat equation, we have the following 
a posteriori split error estimator with partition-of-unity localization:
\begin{align}
|J(u) - J(u_{kh}) | \leq |\eta_{\text{split}}^*|
\coloneqq \left| \sum\limits_m \left( \eta_k^{m,*}
+ \sum_{i\in\mathcal{T}_h^m} \eta_{h}^{i,m,*} \right) \right|,
\end{align}
with the temporal error indicators
\begin{equation}
\begin{aligned}
\eta_k^{m,*} &:=\int\limits_{I_m} \Big( J_u'(u_{kh})(\tilde{u}-u_{k})  
       +(\tilde{u}-u_k,\partial_t z_{kh})_H 
       -(\nabla (\tilde{u}-u_{k}),\nabla u_{kh})_H \Big) \;\mathrm{d}t\\
       &\quad -(\tilde{u}^-(t_m)-u_{k}^-(t_m),[z_{kh}]_m,)_H,
\end{aligned}
\end{equation}
and the spatial error indicators
\begin{equation}
\begin{aligned}
\eta_{h}^{i,m,*} &:=  \int\limits_{I_m}(J_u'(u_{kh})((u_k-u_{kh})\chi_{i,m}) \;\mathrm{d}t
       -\int\limits_{I_m}(\nabla ((u_k-u_{kh})\chi_{i,m}),\nabla z_{kh})_H\;\mathrm{d}t \\
       &\quad 
       +\int\limits_{I_m}((u_k-u_{kh})\chi_{i,m},\partial_t z_{kh})_H \;\mathrm{d}t
       -((u_k^-(t_m)-u_{kh}^-(t_m))\chi_{i,m},[z_{kh}]_m)_H.
\end{aligned}
\end{equation}
\end{proposition}
\begin{proof}
The proof starts as in Proposition \ref{prop_primal_split_heat}, but now taking the adjoint 
residual. The rest is then a combination of the previous three propositions and follows 
conceptionally the same lines.
\end{proof}
\begin{remark}
We emphasize that in this work, we estimate discretization 
errors only. The space-time extension of stationary 
PU-DWR versions such as \cite{MeiRaVih109,RaVih13b,EndtLaWi18} to linear or nonlinear iteration 
errors in our current space-time setting is part of future work.
\end{remark}

\subsection{The partition-of-unity approach for the combustion problem}
\label{sec_PU_combustion}
 In this section, we state the error estimator of the combustion problem:
 \begin{proposition}[Primal split error estimator for combustion]
\label{prop_primal_split_combustion}
 Let us assume homogeneous boundary conditions on $\Gamma_N$ as well 
 as for the species concentration on $\Gamma_R$. 
 For the temperature we have the cooling conditions
 $\kappa\theta+\partial_n\theta=0$ on the Robin boundary, 
 such that $g_N^\theta = g_N^Y  = g_R^\theta = g_R^Y \equiv 0$, 
 $a_R^\theta = \kappa$, $b_R^\theta=1$ and $a_R^Y=b_R^Y=0$ hold.
 Then, we have the following 
 a posteriori primal split error estimator with partition-of-unity localization
 for the space-time formulation of the time-dependent combustion problem:
 \begin{align}
 |J(\{\theta,Y\}) - &J(\{\theta,Y\}_{kh}) | \leq 
 |\eta| := \left|\sum\limits_m \left( \eta_{k}^m 
 +\sum\limits_{i\in\mathcal{T}_h^m} \eta_{i,h}^m\right) \right|,
 \end{align}
 with the temporal error indicators
 \begin{equation}
 \begin{aligned}
 \eta_{k}^m = &-\int\limits_{I_m} \Big( (\partial_t \theta_{kh},z^\theta-z^\theta_k)_H
 + (\nabla\theta_{kh},\nabla(z^\theta-z^\theta_k))_H \Big) \mathrm{d}t\\
 &+\int\limits_{I_m} \int\limits_{\Gamma_R}\kappa\theta(z^\theta-z^\theta_k)\mathrm{d}s \mathrm{d}t\\
 &+\int\limits_{I_m} (\partial_t Y_{kh},z^Y-z^Y_k)_H + (\nabla Y_{kh},\nabla(z^Y-z^Y_k))_H \mathrm{d}t\\
 &-\int\limits_{I_m} (\omega(\theta_{kh},Y_{kh}), z^\theta-z^\theta_k)_H
 +(\omega(\theta_{kh},Y_{kh}), z^Y-z^Y_k)_H\mathrm{d}t\\
 &-([\theta_{kh}]_{m-1},z^{\theta,+}(t_{m-1})-z^{\theta,+}_{k}(t_{m-1}))_H\\
 &-([Y_{kh}]_{m-1},z^{Y,+}(t_{m-1})-z^{Y,+}_{k}(t_{m-1}))_H,
 \end{aligned}
 \end{equation}
 and the spatial error indicators
 \begin{equation}
 \begin{aligned}
 \eta_{i,h}^m = &-\int\limits_{I_m}(\partial_t \theta_{kh},(z^\theta_k-z^\theta_{kh})\chi_{i,m})_H 
 + (\nabla\theta_{kh},\nabla((z^\theta_k-z^\theta_{kh})\chi_{i,m}))_H\mathrm{d}s\mathrm{d}t\\
 &+\int\limits_{I_m}\int\limits_{\Gamma_R}\kappa\theta(z^\theta_k-z^\theta_{kh})
 \chi_{i,m}\mathrm{d}s\mathrm{d}t\\
 &+ \int\limits_{I_m} (\partial_t Y_{kh},(z^Y_k-z^Y_{kh})\chi_{i,m})_H 
 + (\nabla Y_{kh},\nabla((z^Y_k-z^Y_{kh})\chi_{i,m}))_H \mathrm{d}s \mathrm{d}t \\
 &- \int\limits_{I_m} (\omega(\theta_{kh},Y_{kh}), (z^\theta_k-z^\theta_{kh})\chi_{i,m})_H
 +(\omega(\theta_{kh},Y_{kh}), (z^Y_k-z^Y_{kh})\chi_{i,m})_H\mathrm{d}t \\
 &-([\theta_{kh}]_{m-1},(z^{\theta,+}_{k}(t_{m-1})-z^{\theta,+}_{kh}(t_{m-1}))\chi_{i,m})_H \\
 &-([Y_{kh}]_{k,m-1},(z^{Y,+}_k(t_{m-1})-z^{Y,+}_{kh}(t_{m-1}))\chi_{i,m})_H.
 \end{aligned}
 \end{equation}
 \end{proposition}
\begin{proof}
We start as in Proposition \ref{prop_primal_split_heat} and employ 
for the primal residual and the right hand side the weak form of 
the combustion problem, i.e., \eqref{eq_combustion_weak_form_lhs}
and \eqref{eq_combustion_weak_form_rhs}, respectively.
\end{proof}

\begin{proposition}[Adjoint split error estimator for combustion]
\label{prop_adjoint_split_estimator_combustion}
 Using the same boundary conditions, we have the following 
 a posteriori adjoint split error estimator with partition-of-unity localization
 for the space-time formulation of the time-dependent combustion problem:
 \begin{align}
 |J(\{\theta,Y\}) - &J(\{\theta,Y\}_{kh}) | \leq |\eta| := 
\left| \sum\limits_m\left(\eta_{k}^{m,*} +\sum\limits_{i\in\mathcal{T}_h^m} \eta_{i,h}^{m,*}\right) \right|,
 \end{align}
 with the temporal error indicators
 \begin{align*}
 \eta_{k}^{m,*} = &
 J'_\theta (\{\theta_{kh},Y_{kh}\})(\theta-\theta_k)|_{I_m}+J'_Y (\{\theta_{kh},Y_{kh}\})(Y-Y_k)|_{I_m}\\
 &-\int\limits_{I_m} (\theta-\theta_k,-\partial_t z_{kh}^\theta)_H + (\nabla(\theta-\theta_k),\nabla z_{kh}^\theta)_H
 +\int\limits_{\Gamma_R} \kappa(\theta-\theta_k)z_{kh}^\theta\mathrm{d}s\mathrm{d}t\\
  &-\int\limits_{I_m} (Y-Y_k,\partial_t z_{kh}^Y)_H+ (\nabla(Y-Y_k),\nabla z_{kh}^Y))_H\mathrm{d}t\\
  &-\int\limits_{I_m} (\omega_\theta'(\theta_{kh},Y_{kh})(\theta-\theta_k)+\omega_Y'(\theta_{kh},Y_{kh})(Y-Y_k),z_{kh}^Y-z_{kh}^\theta)_H
  \mathrm{d}t\\
  &-((\theta^-(t_m)-\theta_{k}^-(t_m)),[z_{kh}^\theta]_m)-((Y^-(t_m)-Y_{k}^-(t_m)),[z_{kh}^Y]_m),
 \end{align*}
 and the spatial indicators
 \begin{align*}
 \eta_{i,h}^{m,*}=&
 J'_\theta (\{\theta_{kh},Y_{kh}\})((\theta_k-\theta_{kh})\chi_{i,m})|_{I_m}+J'_Y (\{\theta_{kh},Y_{kh}\})((Y_k-Y_{kh})\chi_{i,m})|_{I_m}\\
 &-\int\limits_{I_m} ((\theta_k-\theta_{kh})\chi_{i,m},-\partial_t z_{kh}^\theta)_H + (\nabla((\theta_k-\theta_{kh})\chi_{i,m}),\nabla z_{kh}^\theta)_H
 +\int\limits_{\Gamma_R} \kappa(\theta_k-\theta_{kh})\chi_{i,m}z_{kh}^\theta\mathrm{d}s\mathrm{d}t\\
  &-\int\limits_{I_m} ((Y_k-Y_{kh})\chi_{i,m},\partial_t z_{kh}^Y)_H+ (\nabla((Y_k-Y_{kh})\chi_{i,m}),\nabla z_{kh}^Y))_H\mathrm{d}t\\
  &-\int\limits_{I_m} (\omega_\theta'((\theta_k-\theta_{kh})\chi_{i,m})(\theta_k-\theta_{kh})
  +\omega_Y'(\theta_{kh},Y_{kh})((Y_k-Y_{kh})\chi_{i,m}),z_{kh}^Y-z_{kh}^\theta)_H\mathrm{d}t\\
  &-((\theta_{k}^-(t_m)-\theta_{kh}^-(t_m))\chi_{i,m},[z_{kh}^\theta]_m)-((Y_{k}^-(t_m)-Y_{kh}^-(t_m))\chi_{i,m},[z_{kh}^Y]_m).
 \end{align*}
 \end{proposition}
\begin{proof}
We start as in Proposition \ref{prop_adjoint_spit_estimator_heat}
and the respective definition of $J$ in \eqref{eq_goal_functional}. In contrast, 
the adjoint $A'_u$ is now derived from the combustion system 
by approximating $u=(\theta,Y)$ 
by some higher-order 
approximation $\tilde u = (\tilde\theta,\tilde Y)$ yields the assertion.
\end{proof}

\subsection{Evaluation of the space-time PU-DWR}
\label{sec:eval_pu_dwr}
For the practical evaluation we need to 
properly define the interpolation differences.
Depending on the approach, we need interpolations 
from a high order space into a low order space 
and reconstructions the other way around. 
In space we denote them as
\[
i_h^{(s+1)}\colon \widetilde{X}_{k,h}^{r,s+1}\mapsto \widetilde{X}_{k,h}^{r,s}
\quad\text{and}\quad
i_{2h}^{(s+1)}\colon \widetilde{X}_{k,h}^{r,s}\mapsto \widetilde{X}_{k,h}^{r,s+1}.
\]
In time we use 
\[
i_k^{(r+1)}\colon \widetilde{X}_{k,h}^{r+1,s}\mapsto \widetilde{X}_{k,h}^{r,s}
\quad\text{and}\quad
i_{2k}^{(r+1)}\colon \widetilde{X}_{k,h}^{r,s}\mapsto \widetilde{X}_{k,h}^{r+1,s}.
\]
In the following, we take a closer look at the interpolation 
difference for a higher order solution as used in 
the mixed and equal high order approach, i.e. the interpolations $i_h^{(s+1)}$ and $i_k^{(r+1)}$.
After that we write down the resulting localized error estimators 
for each PU-DoF.
For good visual representations of the high order reconstructions based on a low order solution
see \cite{richter.wick2015} and \cite{schmich.vexler2008} for the spatial part $i_{2h}^{(s+1)}$ and temporal part $i_{2k}^{(r+1)}$ respectively.

\begin{figure}[h!]
\centering
\subfloat[][$\tilde{z} \in X_{k,h}^{1,2}$]{
\includegraphics[trim = 170 300 250 450, clip, width=.3\textwidth ]{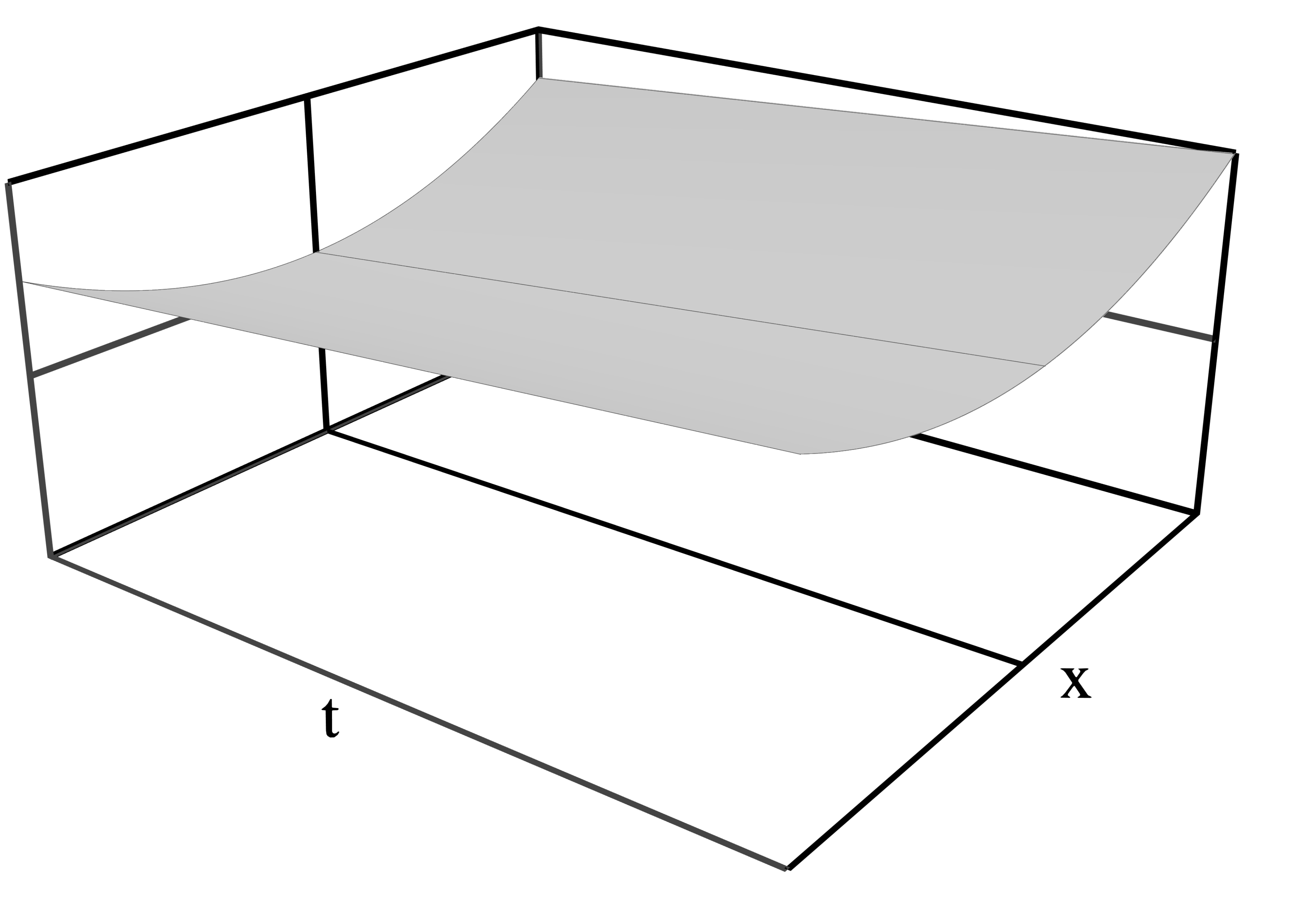}
}
\subfloat[][$i_kz\in \widetilde{X}_{k,h}^{0,2}$]{
\includegraphics[trim = 170 300 250 450, clip, width=.3\textwidth ]{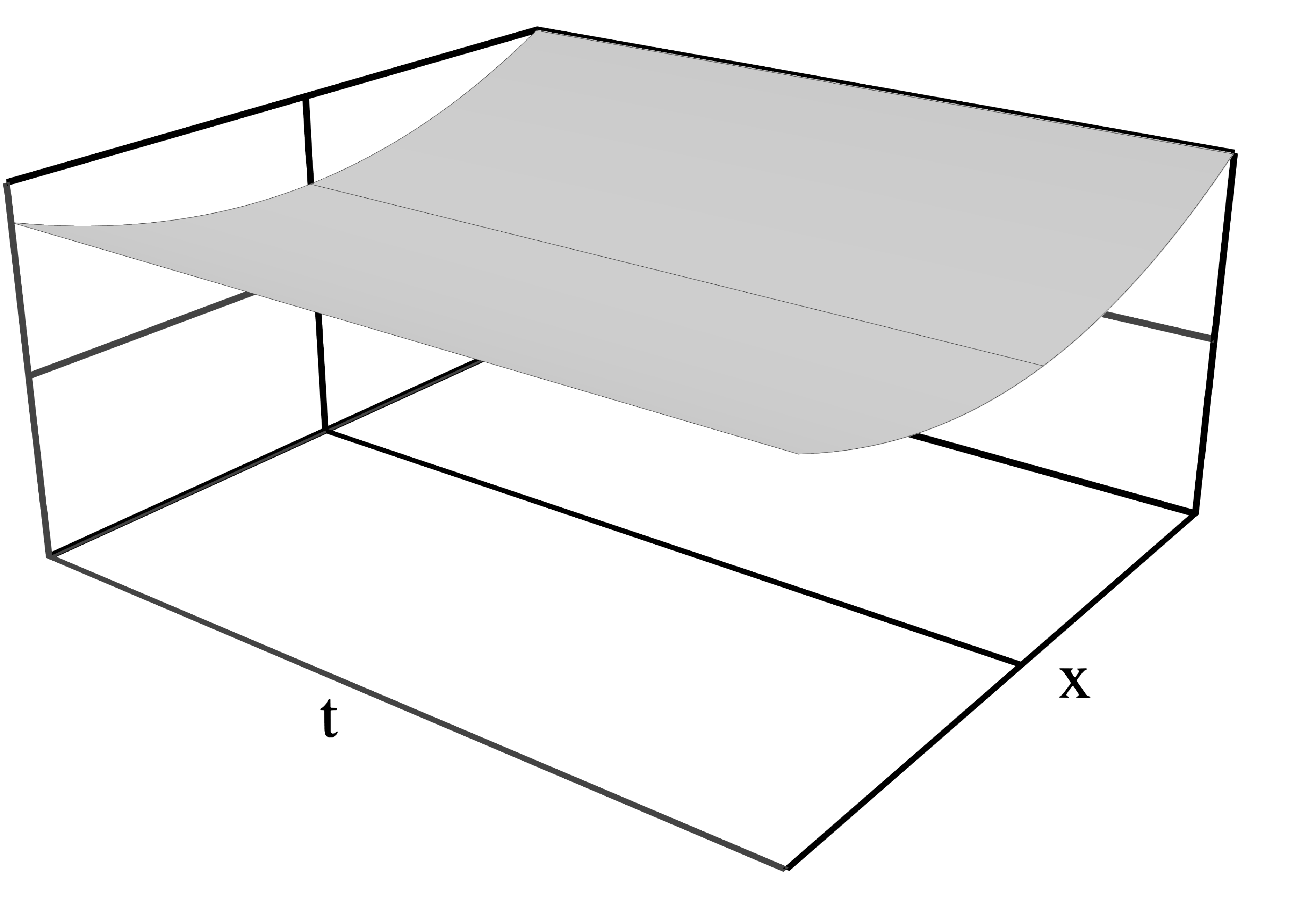}
}
\subfloat[][$i_{kh}z\in \widetilde{X}_{k,h}^{0,1}$]{
\includegraphics[trim = 170 300 250 450, clip, width=.3\textwidth ]{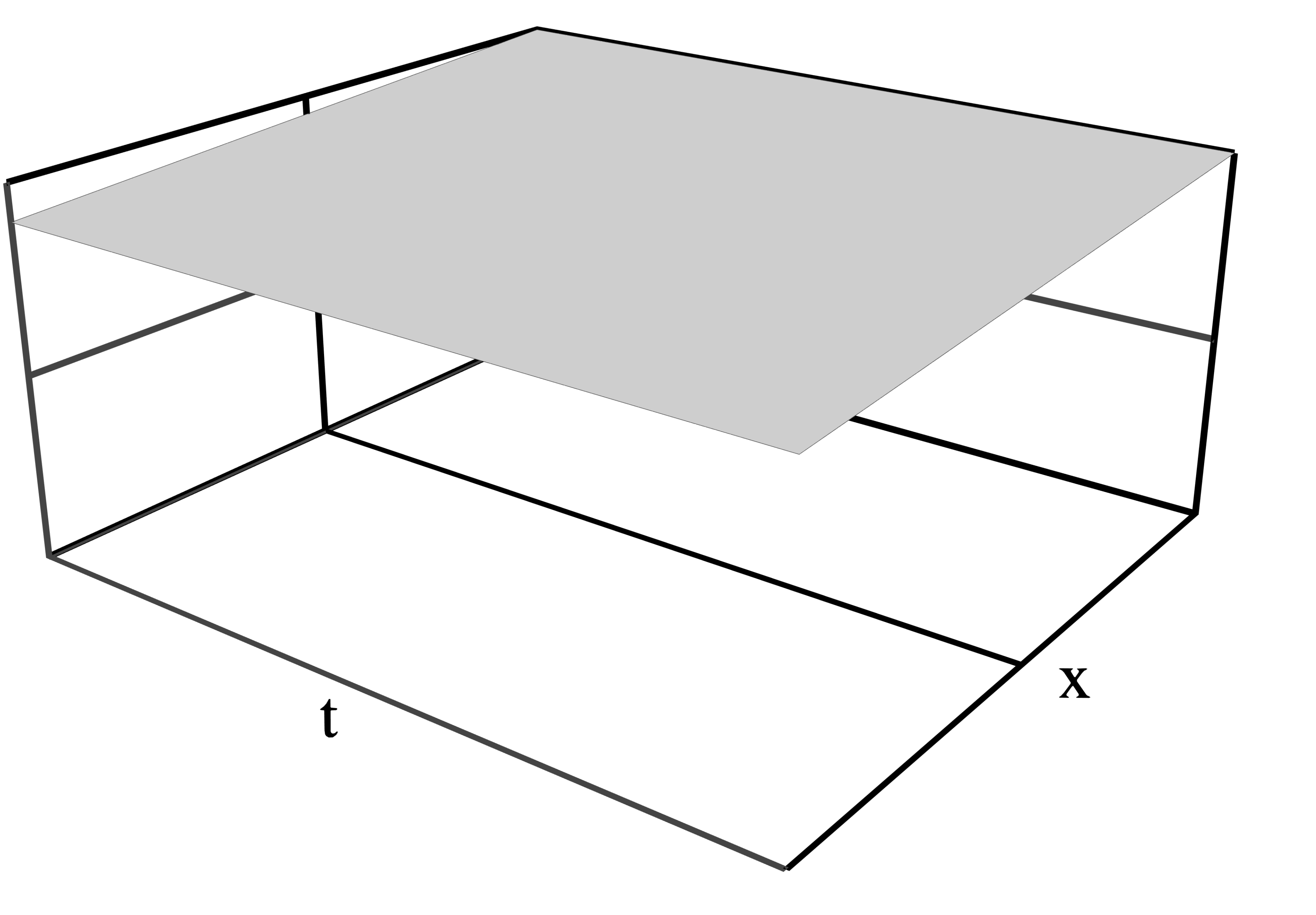}
}
\caption{Different interpolation levels on a single 1+1D space-time element
}
\label{fig:interpolations}
\end{figure}

For our visualization the high order space $\widetilde{X}_{k,h}^{1,2}$ 
and the low order space $\widetilde{X}_{k,h}^{0,1}$ are used.
Then, $u_k$ and $z_k$ are elements of $\widetilde{X}_{k,h}^{0,2}$.
Furthermore, we notice that in this specific case of piecewise constant discrete solutions the identities
\(u_{k}^{-}(t_{m})=u_{k}^{+}(t_{m-1})\) and \(u_{kh}^{-}(t_{m})=u_{kh}^{+}(t_{m-1})\) hold.
With this, we obtain the following operators:
\begin{definition}[Temporal interpolation operators for $r=0$]\ \\
The interpolation operator from piecewise linear elements to 
piecewise constant elements reads as
\begin{align}
i_k^1 \tilde{z}(t) = \begin{cases}
\tilde{z}^{-}(t_m) &\text{ for } t\in I_m,\\
\tilde{z}^{-}(t_1) &\text{ for } t = 0.
\end{cases}  
\end{align}
The reconstruction of the high order solution, i.\,e.\ the other way around, is obtained by linear interpolation
\[
\tilde{z}|_{I_m}(t)=\frac{t_m-t}{k_m}z_k^{-}(t_{m-1})+\frac{t-t_{m-1}}{k_m} z_k^{-}(t_m).
\] 
\end{definition}
For the spatial interpolation operator $i_h^2$ we use a linear finite element ansatz with the
vertex DoFs of the spatial triangulation. Using $i_h^2$ 
on the temporally interpolated solution 
$i_k \tilde{z}$ yields $i_{kh}\tilde{z}$ 
and vice versa.
To illustrate this, Figure \ref{fig:interpolations} 
shows the different interpolation levels for a single 1$+$1D finite element. 

\begin{definition}[Application of operators depending on the choice of solution spaces]
Let $\hat{u}$ and $\hat{z}$ denote the approximated solutions of the primal and adjoint problems,
depending on the choice of finite elements for each problem. 
Then, we obtain our terms by the following evaluations:
\begin{equation}
\label{interpolations_kh}
  u_{kh}^{-}(t_m) =\begin{cases}
	\hat{u}^{-}(t_m) &\text{for } \hat{u}\in \widetilde{X}_{k,h}^{0,1},\\
	i_k^1\hat{u}(t) &\text{for } \hat{u}\in X_{k,h}^{1,1},\\
	i_h^2\hat{u}^-(t_m) &\text{for } \hat{u}\in \widetilde{X}_{k,h}^{0,2},\\
	i_h^2i_k^1\hat{u}(t)  &\text{for } \hat{u}\in X_{k,h}^{1,2},
                \end{cases}                
\end{equation}
\begin{equation}
\label{interpolations_k}
z_{k}^{-}(t_m) = \begin{cases}
	i_{2h}^2\hat{z}^{-}(t_m) &\text{for } \hat{z}\in \widetilde{X}_{k,h}^{0,1},\\
	i_k^1i_{2h}^2\hat{z} &\text{for } \hat{z}\in X_{k,h}^{1,1},\\
	\hat{z}^{-}(t_m) &\text{for } \hat{z}\in \widetilde{X}_{k,h}^{0,2},\\
	i_k^1\hat{z} &\text{for } \hat{z}\in X_{k,h}^{1,2}.
               \end{cases} 
\end{equation}
For $z_{kh}^{-}(t_m)$ we have the same interpolations on $\hat{z}$ as for $u_{kh}^{-}(t_m)$ on $\hat{u}$. 
\end{definition}
For these finite element spaces we can use the midpoint 
rule with $t_0 = (t_{m+1}-t_m)/2$ for all temporal integrals when the 
resulting terms are linear in time. For temporal nonlinearities
and
higher-order right hand side functions $f$, higher-order quadrature
rules, usually Gauss quadratures, have to be used.
\begin{remark}[Reconstructions for the adjoint estimator]
 For the adjoint estimator we additionally need $u_k$ and $u$.
 The semi-discrete $u_k$ can be obtained the same way as $z_k$, but for $u$ we need to change the interpolation direction, which results in
 \[
   \tilde{u}|_{I_m}(t)=\frac{t_m-t}{k_m}u_k^{-}(t_m)+\frac{t-t_{m-1}}{k_m}u_k^{-}(t_{m+1}).
 \] 
\end{remark}
\begin{remark}
Finally, we comment on the treatment of pointwise evaluations such as 
in 
Proposition \ref{prop_adjoint_spit_estimator_heat} and
Proposition \ref{prop_adjoint_split_estimator_combustion}
in the terms 
$u_k^-(t_m)-u_{kh}^-(t_m)$
and
$Y_{k}^-(t_m)-Y_{kh}^-(t_m)$, respectively.
As an example, we explain more details for the heat equation. Due to the 
reverse construction into a higher order space using $i_{2k}^{(r+1)}$, it holds 
$u^-(t_m) = u_k^{-}(t_{m+1})$ for which we deal with a jump 
at $t_m$, which in general is not identically equal to zero.
Moreover, we note that these jump terms include a spatial integration
which is done by Gauss-Legendre quadrature, i.e. with quadrature points 
that do not lie on the boundaries of the spatial elements. 
Therefore, $((u_k^{-}(t_m)-u_{kh}^{-}(t_m))\chi_{i,m},[z_{kh}]_m)$ and corresponding terms
will in general be nonzero as well, irrespective of the spatial interpolation. 
Conversely, if the difference is zero then the low order solution is an exact representation 
of the high order solution such that no refinement is needed.
\end{remark}

\subsection{Error indicators in space and time}
\label{sec:error_indi}
With the previous evaluations, we can now define the respective indicators 
in space and time for the heat equation and the combustion problem.
Note that the temporal derivative $\partial_t u_{kh}$ 
vanishes for $u_{kh}\in \widetilde{X}_{k,h}^{0,s}$,
i.e. piecewise constant elements in time.

\subsubsection{Natural PU $cG(1)dG(0)$}
\label{sec_natural_PU_error_indicators}
Employing the previously introduced PU, namely $cG(1)dG(0)$ yields 
to the following results for the error indicators.
\begin{proposition}[Joint primal error indicator for the heat equation]
We have the following joint error indicator
for the heat equation
\begin{equation}
\begin{aligned}
\eta_{kh,\text{heat}}^{i,m} =& \int\limits_{t_{m-1}}^{t_m}(f(t),(\tilde{z}(t)-z_{kh}^{-}(t_m))\chi_{i,m})_H\mathrm{d}t\ 
-(u_{kh}^{-}(t_m)-u_{kh}^{-}(t_{m-1}),(z_k^{-}(t_{m-1})-z_{kh}^{-}(t_m))\chi_{i,m})_H\\
-& \frac{k_m}{2} \cdot(\nabla u_{kh}^{-}(t_m),(\nabla z_k^{-}(t_{m-1})+\nabla z_k^{-}(t_m)-2\nabla z_{kh}^{-}(t_m))\chi_{i,m})_{H}\\
-& \frac{k_m}{2} \cdot(\nabla u_{kh}^{-}(t_m),(z_k^{-}(t_{m-1})+z_k^{-}(t_m)-2z_{kh}^{-}(t_m))\nabla\chi_{i,m})_H
\end{aligned}
\end{equation}
with the time step size $k_m=t_m-t_{m-1}$.
\end{proposition}
\begin{proof}
Starting from Proposition \ref{prop_primal_joint_estimator_heat} and the 
representations of $\eta_{kh}^{i,m}$, we evaluate the temporal integrals and employ 
the interpolation operators derived in Section \ref{sec:eval_pu_dwr} and 
obtain the joint error indicator $\eta_{kh,\text{heat}}^{i,m}$.
\end{proof}
Accordingly, we have
\begin{proposition}[Split primal error indicators for the heat equation]
\label{prop:split_primal_heat_indicator}
The split indicators $\eta_{k,\text{heat}}^{m}$ and $\eta_{h,\text{heat}}^{i,m}$ for the heat equation are given by
\begin{align}
\eta_{k,\text{heat}}^{m} &= \int\limits_{t_{m-1}}^{t_m}(f(t),(\tilde{z}(t)-z_k^{-}(t_m)))_H \mathrm{d}t
- k_m/2\cdot(\nabla u_{kh}^{-}(t_m),\nabla (z_k^{-}(t_{m-1})-z_k^{-}(t_m)))_H\nonumber\\
&\hphantom{=} -(u_{kh}^{-}(t_m)-u_{kh}^{-}(t_{m-1}),(z_k^{-}(t_{m-1})-z_k^{-}(t_m)))_H,
\end{align}
and
\begin{align}
\eta_{h,\text{heat}}^{i,m}&= \int\limits_{t_{m-1}}^{t_m}(f(t),(z_k^{-}(t_m)-z_{kh}^{-}(t_m))\chi_{i,m})_H\mathrm{d}t\nonumber\\ 
&\hphantom{=}- k_m\cdot(\nabla u_{kh}^{-}(t_m),(\nabla z_k^{-}(t_m)-\nabla z_{kh}^{-}(t_m))\chi_{i,m}
  +(z_k^{-}(t_m)-z_{kh}^{-}(t_m))\nabla\chi_{i,m})_H\nonumber\\
&\hphantom{=}  -(u_{kh}^{-}(t_m)-u_{kh}^{-}(t_{m-1}),(z_k^{-}(t_m)-z_{kh}^{-}(t_m))\chi_{i,m})_H.
\end{align}
\end{proposition}
\begin{proof}
Starting from Proposition \ref{prop_primal_split_heat} and the 
representations of $\eta_{k}^{m}$ and $\eta_{h}^{i,m}$, 
we evaluate the temporal integrals and employ 
the interpolation operators derived in Section \ref{sec:eval_pu_dwr} and 
obtain the error indicators $\eta_{k,\text{heat}}^{m}$ and $\eta_{h,\text{heat}}^{i,m}$.
\end{proof}

Using the same interpolations we obtain
\begin{proposition}[Split primal error indicators for combustion]
\label{prop:split_primal_combustion_indicator}
For the combustion problem we have the following primal error indicators
\begin{equation}
\begin{aligned}
 \eta_{k,\text{combustion}}^{m} = -k_m/2\Big[ &(\nabla\theta_{kh}^{-}(t_m),\nabla(z_k^{\theta,-}(t_{m-1})-z_k^{\theta,-}(t_m)))_H\\
 &(\nabla Y_{kh}^{-}(t_m),\nabla(z_k^{Y,-}(t_{m-1})-z_k^{Y,-}(t_m)))_H \\
 &+ \int\limits_{\Gamma_R}\kappa\theta_{kh}^{-}(t_m)(z_k^{\theta,-}(t_{m-1})-z_k^{\theta,-}(t_m))\mathrm{d}s \\
 &-(\omega(\theta_{kh}^{-}(t_m),Y_{kh}^{-}(t_m)),z_k^{\theta,-}(t_{m-1})-z_k^{\theta,-}(t_m))_H\\
 &+(\omega(\theta_{kh}^{-}(t_m),Y_{kh}^{-}(t_m)),z_k^{Y,-}(t_{m-1})-z_k^{Y,-}(t_m))_H\Big]\\
 &-(\theta_{kh}^{-}(t_m)-\theta_{kh}^{-}(t_{m-1}),z_k^{\theta,-}(t_{m-1})-z_k^{\theta,-}(t_m))_H\\
 &-(Y_{kh}^{-}(t_m)-Y_{kh}^{-}(t_{m-1}),z_k^{Y,-}(t_{m-1})-z_k^{Y,-}(t_m))_H
\end{aligned}
\end{equation}
and
\begin{equation}
\begin{aligned}
 \eta_{h,\text{combustion}}^{i,m} = -k_m\Big[&(\nabla\theta_{kh}^{-}(t_m),\nabla((z_k^{\theta,-}(t_{m})-z_{kh}^{\theta,-}(t_m))\chi_{i,m}))_H\\
 &(\nabla Y_{kh}^{-}(t_m),\nabla((z_k^{Y,-}(t_{m})-z_{kh}^{Y,-}(t_m))\chi_{i,m}))_H \\
 &+ \int\limits_{\Gamma_R}\kappa\theta_{kh}^{-}(t_m)(z_k^{\theta,-}(t_{m})-z_{kh}^{\theta,-}(t_m))\chi_{i,m}\mathrm{d}s \\
 &-(\omega(\theta_{kh}^{-}(t_m),Y_{kh}^{-}(t_m)),(z_k^{\theta,-}(t_{m})-z_{kh}^{\theta,-}(t_m))\chi_{i,m})_H\\
 &+(\omega(\theta_{kh}^{-}(t_m),Y_{kh}^{-}(t_m)),(z_k^{Y,-}(t_{m})-z_{kh}^{Y,-}(t_m))\chi_{i,m})_H\Big]\\
 &-(\theta_{kh}^{-}(t_m)-\theta_{kh}^{-}(t_{m-1}),(z_k^{\theta,-}(t_{m})-z_{kh}^{\theta,-}(t_m))\chi_{i,m})_H\\
 &-(Y_{kh}^{-}(t_m)-Y_{kh}^{-}(t_{m-1}),(z_k^{Y,-}(t_{m})-z_{kh}^{Y,-}(t_m))\chi_{i,m})_H.
\end{aligned}
\end{equation}
\end{proposition}
\begin{proof}
Starting from Proposition \ref{prop_primal_split_combustion} and the 
representations of $\eta_{k}^{m}$ and $\eta_{h}^{i,m}$, 
we evaluate the temporal integrals and employ 
the interpolation operators derived in Section \ref{sec:eval_pu_dwr} and 
obtain the error indicators $\eta_{k,\text{combustion}}^{m}$ and $\eta_{h,\text{combustion}}^{i,m}$.
\end{proof}

\begin{remark}[Identity of the indicator variants]
Since $\sum\limits_{i\in\mathcal{T}_h^m} \chi_{i,m} \equiv 1$ in Proposition \ref{prop_PU} holds,
we 
have
\begin{equation}
 \eta_k^m + \sum\limits_{i\in\mathcal{T}_h^m} \eta_{h}^{i,m}
 = \sum\limits_{i\in\mathcal{T}_h^m} \eta_{kh}^{i,m},
\end{equation}
such that the choice between the two indicator variants is only important for adaptive refinement.
If one is only interested in estimating the error then both variants are identical. 
\end{remark}

\subsubsection{Alternative PU $cG(1)cG(1)$}
To investigate the impact of the choice of the PU space we derive 
the split primal indicators for the heat equation based on 
a $cG(1)cG(1)$ PU. Now, the right hand side integration might 
need even higher order quadrature rules. Apart from that, the
highest temporal order in the temporal 
indicators is quadratic (constant primal solution, linear adjoint 
solution and PU) so we apply Simpsons rule instead.
Additionally, we obtain two sets of spatial and temporal 
indicators per interval $I_m$, i.e. 
\[
\eta_{k,m,\text{heat}}^{m-1}, \quad\eta_{k,m,\text{heat}}^{m}, \quad
\eta_{h,m,\text{heat}}^{i,m-1}, \quad\eta_{h,m,\text{heat}}^{i,m}.
\]
The temporal indicators for refinement of $I_m$ are then obtained by
\begin{equation}
 \eta_{k,\text{heat},cG(1)}^m = 
 \sum\limits_{i\in V_h^1(\mathcal{T}_h^{m-1})} \eta_{k,{m-1},\text{heat}}^{m-1}
 +\sum\limits_{i\in V_h^1(\mathcal{T}_h^{m})} \eta_{k,{m},\text{heat}}^{m-1} +\eta_{k,{m},\text{heat}}^{m}
 +\sum\limits_{i\in V_h^1(\mathcal{T}_h^{m+1})} \eta_{k,{m+1},\text{heat}}^{m}.
\end{equation}
For the spatial element indicators we first have to interpolate the indicator vectors 
$(\eta_{h,m-1,\text{heat}}^{i,m-1})_{i\in V_h^1(\mathcal{T}_h^{m-1})}$ and
$(\eta_{h,m+1,\text{heat}}^{i,m})_{i\in V_h^1(\mathcal{T}_h^{m+1})}$ to $\mathcal{T}_h^m$.
Then, the element indicator for $K\in\mathcal{T}_h^m$ is calculated as
\begin{equation}
 \eta_{h,\text{heat},cG(1)}^{K,m} = \sum\limits_{i\in K} 
  \eta_{h,m-1,\text{heat}}^{i,m-1}  + \eta_{h,m,\text{heat}}^{i,m-1}+ \eta_{h,m,\text{heat}}^{i,m} + \eta_{h,m+1,\text{heat}}^{i,m}.
\end{equation}
Employing these derivations for the new choice of the PU, Simpson's rule for 
quadrature in time, 
and then proceeding 
as in Section \ref{sec_natural_PU_error_indicators}, we obtain the following 
results.
\begin{proposition}[Split primal error indicators for the heat equation with $cG(1)cG(1)$ PU]
\label{prop:split_primal_heat_indicator_cG1_PU}
The split indicators for the heat equation are given by
\begin{equation}
\begin{aligned}
\eta_{k,m,\text{heat}}^{m-1} &= \int\limits_{t_{m-1}}^{t_m}(f(t),(\tilde{z}(t)-z_k^{-}(t_m))\frac{t_m-t}{t_m-t_{m-1}})_H \mathrm{d}t
- 2k_m/6\cdot(\nabla u_{kh}^{-}(t_m),\nabla (z_k^{-}(t_{m-1})-z_k^{-}(t_m)))_H\\
&\hphantom{=} -(u_{kh}^{-}(t_m)-u_{kh}^{-}(t_{m-1}),(z_k^{-}(t_{m-1})-z_k^{-}(t_m)))_H,
\end{aligned}
\end{equation}
and
\begin{equation}
\begin{aligned}
\eta_{k,m,\text{heat}}^{m} =& \int\limits_{t_{m-1}}^{t_m}(f(t),(\tilde{z}(t)-z_k^{-}(t_m))\frac{t-t_{m-1}}{t_m-t_{m-1}})_H \mathrm{d}t\\
&- k_m/6\cdot(\nabla u_{kh}^{-}(t_m),\nabla (z_k^{-}(t_{m-1})-z_k^{-}(t_m)))_H,
\end{aligned}
\end{equation}
as well as
\begin{equation}
\begin{aligned}
\eta_{h,m,\text{heat}}^{i,m-1}&= \int\limits_{t_{m-1}}^{t_m}(f(t),(z_k^{-}(t_m)-z_{kh}^{-}(t_m))\frac{t_m-t}{t_m-t_{m-1}}\chi_{i,m})_H\mathrm{d}t\\ 
&\hphantom{=}- k_m/2\cdot(\nabla u_{kh}^{-}(t_m),(\nabla z_k^{-}(t_m)-\nabla z_{kh}^{-}(t_m))\chi_{i,m}+(z_k^{-}(t_m)-z_{kh}^{-}(t_m))\nabla\chi_{i,m})_H\\
&\hphantom{=}  -(u_{kh}^{-}(t_m)-u_{kh}^{-}(t_{m-1}),(z_k^{-}(t_m)-z_{kh}^{-}(t_m))\chi_{i,m})_H,
\end{aligned}
\end{equation}
and
\begin{equation}
\begin{aligned}
\eta_{h,m,\text{heat}}^{i,m}&= \int\limits_{t_{m-1}}^{t_m}(f(t),(z_k^{-}(t_m)-z_{kh}^{-}(t_m))\frac{t-t_{m-1}}{t_m-t_{m-1}}\chi_{i,m})_H\mathrm{d}t\\ 
&\hphantom{=}- k_m/2\cdot(\nabla u_{kh}^{-}(t_m),(\nabla z_k^{-}(t_m)-\nabla z_{kh}^{-}(t_m))\chi_{i,m}+(z_k^{-}(t_m)-z_{kh}^{-}(t_m))\nabla\chi_{i,m})_H\\
\end{aligned}
\end{equation}
\end{proposition}

For the general adjoint estimator we also need $z_{kh}$ and $u_k$ from $I_{m+1}$
which can be obtained by the interpolations described in \eqref{interpolations_kh} and \eqref{interpolations_k} respectively.
Additionally we only look at goal functionals of the types
\begin{align*}
 J_1(u)(\varphi) &= \int\limits_0^T (\bar{J}_1(u),\varphi)_H \mathrm{d}t, \qquad
 J_2(u)(\varphi) = \int\limits_0^T \int\limits_{\partial\Omega} \bar{J}_2(u)\varphi \mathrm{d}s\mathrm{d}t,
\end{align*}
which are essentially interchangeable in the following formulas. 
Therefore, we only write down the indicators for $J_1$.

\begin{proposition}[Joint adjoint error indicator for the heat equation]
We have the following joint error indicator
for the heat equation
\begin{equation}
\begin{aligned}
 \eta_{kh,\text{heat}}^{i,m,*} &= \int\limits_{t_{m-1}}^{t_m} J_u'(u_{kh})((\tilde{u}(t)-u_{kh}^{-}(t_m))\chi_{i,m})\mathrm{d}t
+((u_k^{-}(t_{m+1})-u_{kh}^{-}(t_m))\chi_{i,m},
z_{kh}^{-}(t_{m+1})-z_{kh}^{-}(t_{m}))_H\\
&\hphantom{=}- \frac{k_m}{2} \cdot((\nabla u_k^{-}(t_{m+1})+\nabla u_k^{-}(t_m)-2\nabla u_{kh}^{-}(t_m))\chi_{i,m},\nabla z_{kh}^{-}(t_m))_H.
\end{aligned}
\end{equation}
\end{proposition}

\begin{proposition}[Split adjoint error indicators for the heat equation]
\label{prop:split_adjoint_heat_indicator}
The split indicators $\eta_k^{m,*}$ and $\eta_{h,i}^{m,*}$ for the heat equation are given by
\begin{align}
 \eta_{k,\text{heat}}^{m,*} &= \int\limits_{t_{m-1}}^{t_m}J_u'(u_{kh})(\tilde{u}(t)-u_{k}^{-}(t_m)) \mathrm{d}t
- k_m/2\cdot(\nabla (u_k^{-}(t_{m+1})-u_k^{-}(t_m)),\nabla z_{kh}^{-}(t_m))_H\nonumber\\
&\hphantom{=} +(u_{k}^{-}(t_{m+1})-u_{k}^{-}(t_m),z_{kh}^{-}(t_{m+1})-z_{kh}^{-}(t_m))_H,
\end{align}
and
\begin{align}
 \eta_{h,\text{heat}}^{i,m,*}&= \int\limits_{t_{m-1}}^{t_m}J_u'(u_{kh})((u_k^{-}(t_m)-u_{kh}^{-}(t_m))\chi_{i,m})\mathrm{d}t\nonumber\\ 
&\hphantom{=}- k_m\cdot((\nabla u_k^{-}(t_m)-\nabla u_{kh}^{-}(t_m))\chi_{i,m}
+(u_k^{-}(t_m)-u_{kh}^{-}(t_m)\nabla\chi_{i,m},\nabla z_{kh}^{-}(t_m))_H\nonumber\\
&\hphantom{=}  +((u_k^{-}(t_m)-u_{kh}^{-}(t_m))\chi_{i,m},z_{kh}^{-}(t_{m+1})-z_{kh}^{-}(t_m))_H.
\end{align}
\end{proposition}

All functionals we want to examine for the combustion equation are only dependent on $u_{kh}$ and the primal weight,
which is at most linear in time. Therefore, we can simplify the estimator by also applying the midpoint rule to the functional.

\begin{proposition}[Split adjoint error indicators for combustion]
\label{prop:split_adjoint_combustion_indicator}
For the combustion problem we have the following adjoint error indicators
\begin{equation}
\begin{aligned}
\eta_{k,\text{combustion}}^{m,*} =
 k_m/2\Big[&(J_{1,\theta}'(\theta_{kh}^{-}(t_m),Y_{kh}^{-}(t_m)),\theta_k^{-}(t_{m+1})-\theta_k^{-}(t_m))_H\\
    &+ (J_{1,Y}'(\theta_{kh}^{-}(t_m),Y_{kh}^{-}(t_m)),Y_k^{-}(t_{m+1})-Y_k^{-}(t_m))_H\\
    &-(\nabla(\theta_k^{-}(t_{m+1})-\theta_k^{-}(t_m)),\nabla z_{kh}^{\theta,-}(t_m))_H
    -\int\limits_{\Gamma_R}\kappa (\theta_k^{-}(t_{m+1})-\theta_k^{-}(t_m)) z_{kh}^{\theta,-}(t_m) \mathrm{d}s\\
    &
    -(\nabla(Y_k^{-}(t_{m+1})-Y_k^{-}(t_m)),\nabla z_{kh}^{Y,-}(t_m))_H\\
    &-(\omega_\theta'(\theta_{kh}^{-}(t_m),Y_{kh}^{-}(t_m))(\theta_k^{-}(t_{m+1})-\theta_k^{-}(t_m)),z_{kh}^{Y,-}(t_m)-z_{kh}^{\theta,-}(t_m))_H\\
    &+(\omega_Y'(\theta_{kh}^{-}(t_m),Y_{kh}^{-}(t_m))(Y_k^{-}(t_{m+1})-Y_k^{-}(t_m)),z_{kh}^{Y,-}(t_m)-z_{kh}^{\theta,-}(t_m))_H
 \Big]\\
 &+(\theta_{k}^{-}(t_{m+1})-\theta_{k}^{-}(t_m),z_{kh}^{\theta,-}(t_{m+1})-z_{kh}^{\theta,-}(t_m))_H\\
 &+(Y_{k}^{-}(t_{m+1})-Y_{k}^{-}(t_m),z_{kh}^{Y,-}(t_{m+1})-z_{kh}^{Y,-}(t_m))_H,
\end{aligned}
\end{equation}
and
\begin{equation}
\begin{aligned}
\eta_{h,\text{combustion}}^{i,m,*} =
 k_m\Big[&(J_{1,\theta}'(\theta_{kh}^{-}(t_m),Y_{kh}^{-}(t_m)),(\theta_k^{-}(t_m)-\theta_{kh}^{-}(t_m))\chi_{i,m})_H\\
    &+ (J_{1,Y}'(\theta_{kh}^{-}(t_m),Y_{kh}^{-}(t_m)),(Y_k^{-}(t_m)-Y_{kh}^{-}(t_m))\chi_{i,m})_H\\
    &-(\nabla((\theta_k^{-}(t_m)-\theta_{kh}^{-}(t_m))\chi_{i,m}),\nabla z_{kh}^{\theta,-}(t_m))_H \\
    &-(\nabla((Y_k^{-}(t_m)-Y_{kh}^{-}(t_m))\chi_{i,m}),\nabla z_{kh}^{Y,-}(t_m))_H\\
    &-\int\limits_{\Gamma_R}\kappa (\theta_k^{-}(t_m)-\theta_{kh}^{-}(t_m))\chi_{i,m}z_{kh}^{\theta,-}(t_m) \mathrm{d}s\\
    &-(\omega_\theta'(\theta_{kh}^{-}(t_m),Y_{kh}^{-}(t_m))(\theta_k^{-}(t_m)-\theta_{kh}^{-}(t_m))\chi_{i,m},z_{kh}^{Y,-}(t_m)-z_{kh}^{\theta,-}(t_m))_H\\
    &+(\omega_Y'(\theta_{kh}^{-}(t_m),Y_{kh}^{-}(t_m))(Y_k^{-}(t_m)-Y_{kh}^{-}(t_m))\chi_{i,m},z_{kh}^{Y,-}(t_m)-z_{kh}^{\theta,-}(t_m))_H
 \Big]\\
 &+((\theta_k^{-}(t_m)-\theta_{kh}^{-}(t_m))\chi_{i,m},z_{kh}^{\theta,-}(t_{m+1})-z_{kh}^{\theta,-}(t_m))_H\\
 &+((Y_k^{-}(t_m)-Y_{kh}^{-}(t_m))\chi_{i,m},z_{kh}^{Y,-}(t_{m+1})-z_{kh}^{Y,-}(t_m))_H.
\end{aligned}
\end{equation}
\end{proposition}

\subsection{Adaptive algorithm}
\label{sec_adaptive_algo}
There are multiple options for solving the primal and adjoint problems, 
i.e. time-stepping, time-slabbing and fully simultaneous space-time,
which all need adjustments to marking and refinement. However, 
all simulations shown here are time-stepping based so we limit ourselves
to this approach.
Again, we note that in a space-time context the 
adjoint problem runs backward in time. Then, all information is 
collected to evaluate the error estimators.

Since one of the error components might dominate, we employ an 
equilibration for time-stepping as proposed in
\cite{schmich.vexler2008}. This will sometimes restrict 
refinement to space or time.
The overall procedure follows the typical loop: solve, estimate, mark, and refine.

For error estimation in time-stepping we obtain 
Algorithm~\ref{algo_estimate_timestep}.
There, the main choice is in whether to use the split or joint estimators. 
Note that in the case of the joint estimator, the indicators 
$\eta_{kh}^m$ and $\eta_{kh}^{m,*}$ are 
calculated by summation of the spatial indicators.

\begin{algorithm}
\caption{ESTIMATE on a single interval $I_m$, i.e. time-stepping based}
\begin{algorithmic}
\Require $\hat{u}_{kh}$ on $I_{m-1}$ and $I_m$ and $\hat{z}_{kh}$ on $I_{m-1}$, $I_m$ and $I_{m+1}$
\State interpolate/reconstruct $\tilde{u}$, $u_k$, $u_{kh}$ as well as 
$\tilde{z}$, $z_k$, $z_{kh}$ at quadrature points.
\State Calculate $\eta_{\star}^m$ and $\eta_\star^{m,*}$, 
where $\star$ denotes $k$ or $kh$
\State Calculate $\eta_{\bullet}^{i,m}$ and $\eta_{\bullet}^{i,m,*}$ for each PU-DoF $i$ on $\mathcal{T}_h^m$, 
where, again, $\bullet$ denotes $h$ or $kh$
\end{algorithmic}
\label{algo_estimate_timestep}
\end{algorithm}
Having calculated the estimators we mark and refine elements by following Algorithm \ref{algo_mark_refine_timestep}.
There, the equilibration is done by first calculating the global estimators. 
Note that these coincide in the joint case, such that no equilibration is 
performed.

\begin{algorithm}
\caption{MARK and REFINE for the time-stepping approach}
\begin{algorithmic}
\Require indicators on each interval $I_m$ and equilibration factor $c>0$
\State Calculate global temporal estimator 
$\eta_k= \frac12 \sum\limits_{m=1}^M (\eta_\star^m+ \eta_\star^{m,*})$
\State Calculate global spatial estimator 
$\eta_h= \frac12 \sum\limits_{m=1}^M \sum\limits_{i\in\mathcal{T}_h^m} 
(\eta_\bullet^{i,m} + \eta_\bullet^{i,m,*})$
\If{$|\eta_k|*c \geq |\eta_h|$} 
\State mark $I_m$ for temporal refinement based on chosen strategy
\EndIf
\If{$|\eta_h|*c \geq |\eta_k|$}
\For{$m=1,\dots,M$}
\State mark and refine elements in $\mathcal{T}_h^m$ based on chosen strategy
\EndFor
\EndIf
\For{$m=1,\dots,M$}
\If{$I_m$ is marked}
\State Split/Refine $I_m$ into two intervals with (possibly new) mesh $\mathcal{T}_h^m$
\EndIf
\EndFor
\end{algorithmic}
\label{algo_mark_refine_timestep}
\end{algorithm}

\newpage
\section{Numerical tests}
\label{Section: Numerical tests}
In this final section, we substantiate our 
space-time error estimators and algorithms
with the help of three numerical experiments.
In the first configuration, a $2+1D$ heat equation with manufactured solution 
is considered. This allows us to investigate in detail effectivity 
indices. Next, in Configuration 2, again a $2+1D$ heat equation
is utilized, but with a dynamic manufactured solution inspired 
by Hartmann \cite{hartmann1998}.
In our final configuration, we consider 
a nonlinear coupled problem, namely nonlinear combustion. Therein,
very detailed comparisons of different polynomial degrees, 
primal, adjoint and full estimators are undertaken.
The computations are based on extensions of the 
DTM package \textit{dwr-diffusion} \cite{kocher.etal2019},
which itself is based on deal.II \cite{arndt.etal2020}.
The programming codes are open-source can be found on 
\url{https://github.com/jpthiele/pu-dwr-diffusion} and 
\url{https://github.com/jpthiele/pu-dwr-combustion} respectively,
and follow good practices of sustainable research software developments \cite{Anztetal20}.

\subsection{Configuration 1: 2+1D Heat equation with a simple manufactured solution}
\label{sec:config_simple_heat}
\subsubsection{Problem statement}
To test the 2+1D implementation and the derived estimators we prescribe the following solution for the heat equation 
introduced in Section \ref{sec_heat}
\begin{equation}
 u(t,x,y) = -\frac{(x^2-x)(y^2-y)}4t.
\end{equation}
Inserting the solution into the PDE yields the right hand side function
\begin{equation}
 f(t,x,y) = -\frac{(x^2-x)(y^2-y)}4 + \frac{(x^2-x)}2t  + \frac{(y^2-y)}2t.
\end{equation}
\subsubsection{Configuration}
The PDE is solved on the unit square and the temporal interval $(0,1)$, i.e. $T=1$. 
Inserting $x=0$, $y=0$ or $t=0$ yields $u = 0$, resulting in homogeneous Dirichlet boundary conditions and an initial condition 
of $u^0 \equiv 0$.

\subsubsection{Goal functionals}
To test whether the error identity holds, we need a linear goal functional.
A simple choice is the averaged solution
\begin{equation}
 J(u) = \frac{1}{|\Omega|T} \int\limits_0^T \int\limits_{\Omega}  u(t,x)\mathrm{d}t\mathrm{d}x.
\end{equation}
Inserting the analytical solution for an arbitrary $T$ we obtain
\begin{equation}
 J(u) = -\frac{T}{288},
\end{equation}
as reference value.
\subsubsection{Discussion of findings}
We expect \eqref{eqn:primal_error_rep} and \eqref{eqn:adjoint_error_rep} to hold, which is identical to $I_{\text{eff}}= 1$.
We observe the effectivity indices for all three approaches in space defined as
\[
 I_{\text{eff}}^{s/\widetilde{s}} = \frac{\eta_{kh}^{s/\widetilde{s}}}{J(u)-J(u_{kh})}
\]
with $\eta_{kh}^{s/\widetilde{s}}$ as the general estimator for $\hat{u}_{kh}\in \widetilde{X}_{k,h}^{0,s}$ and $\hat{z}_{kh}\in \widetilde{X}_{k,h}^{0,\widetilde{s}}$
as defined in Propositions \ref{prop:split_primal_heat_indicator} and \ref{prop:split_adjoint_heat_indicator} for the 
primal and adjoint estimators respectively. All estimators are computed using Algorithm \ref{algo_estimate_timestep}.

Table \ref{table:simple_I_effs} shows that all approaches yield effectivity indices very close to $1$.
We also notice, that the biggest difference between primal and adjoint estimator is obtained for the mixed order approach.\\
However, this is not surprising as the approach is tailored to the primal estimator and calculating the 
adjoint estimator could be seen as questionable for the two following reasons. 
When interpolating the higher order solution for the adjoint problem for $z_{kh}$
we do not obtain the optimal $z_{kh}$ compared to solving with bilinear elements directly.
Additionally, reconstructing the higher order primal solution yields a worse approximation compared to solving directly with biquadratic
finite elements.\\
Furthermore, we notice that we use a lot more elements in time than in space. 
In space-time the estimator is dependent on a good balance between space and time discretization. 
For the sake of brevity we investigate this further in the following configuration as the 
solution is much more interesting.

\begin{table}[H]
\centering
\caption{Section \ref{sec:config_simple_heat}: Performance of the primal (left) and adjoint (right) error estimators under global refinement for temporal $dG(0)$ discretization of the adjoint equation. 
Averaged solution functional.}
\label{table:simple_I_effs}
\begin{tabular}{l|l|l|l|l}
\hline
 $M$ & $N$ & $I_{\text{eff}}^{1/1}$ & $I_{\text{eff}}^{1/2}$ & $I_{\text{eff}}^{2/2}$  \\
\hline
 $1000$ &    $64$ & $1.008726$ & $1.000479$ & $1.012815$ \\
 $2000$ &   $256$ & $1.002999$ & $1.001123$ & $1.004264$ \\
 $4000$ &  $1024$ & $1.002757$ & $1.002329$ & $1.003109$ \\
 $8000$ &  $4096$ & $1.004806$ & $1.004705$ & $1.004896$ \\
\hline
\end{tabular}\qquad
\begin{tabular}{l|l|l|l|l}
\hline
 $M$ & $N$ & $I_{\text{eff}}^{1/1}$ & $I_{\text{eff}}^{1/2}$ & $I_{\text{eff}}^{2/2}$  \\
\hline
 $1000$ &    $64$ & $1.018668$ & $1.031224$ & $1.012812$ \\
 $2000$ &   $256$ & $1.005793$ & $1.008949$ & $1.004263$ \\
 $4000$ &  $1024$ & $1.003508$ & $1.004289$ & $1.003109$ \\
 $8000$ &  $4096$ & $1.005001$ & $1.005192$ & $1.004896$ \\
\hline
\end{tabular}
\end{table}

\subsection{Configuration 2: 2+1D Heat equation with a dynamic manufactured solution}
\label{sec:config_hartmann_heat}
\subsubsection{Problem statement}
This test case was designed in \cite{hartmann1998}. We again solve the heat equation.
The manufactured solution is a rotating hill on a unit-square spatial domain 
$\Omega = (0,1)^2$ in the time interval $(0,T),\;T=1$. \\
The manufactured solution is given as
\begin{align}
 u(x,y,t) = \frac{1}{1+50((x-x_0(t))^2+(y-y_0(t))^2)},\\
 x_0(t) = \frac12+\frac14\cos(2\pi t),\\
 y_0(t) = \frac12+\frac14\sin(2\pi t).
\end{align}
The right hand side of the problem is obtained as in Section \ref{sec:config_simple_heat} by inserting this solution into the 
heat equation. Additionally, all definitions for $\eta_{kh}^{s/\widetilde{s}}$ and $I_{\text{eff}}^{s/\widetilde{s}}$ are the same 
using again Algorithm \ref{algo_estimate_timestep} to compute the indicators.
\subsubsection{Goal functional}
Since we are interested in capturing the local behaviour of the solution,
we choose the $L_2$-error as functional of interest, i.e.
\begin{equation}
 J(u_{kh}) = (u-u_{kh},u-u_{kh})^{1/2}.
\end{equation}

\subsubsection{Comparison of the different spatial approaches}
We start by looking at the exact error and the resulting estimators for different choices of $M_{\text{initial}}$.
We see that the estimators in Tables \ref{table:hartmann_comp_Minit_exact1}, \ref{table:hartmann_comp_Minit_eta_low}, 
\ref{table:hartmann_comp_Minit_eta_mixed}
and \ref{table:hartmann_comp_Minit_eta_high} are converging to relatively
stable values with rising $M_{\text{initial}}$.
We can also see that since the temporal elements are only doubled while the spatial elements are quadrupled, the initial number of 
temporal elements has to be large enough for uniform refinement.
In adaptive simulations we can control this better as we can choose different fractions of spatial and temporal elements to be 
marked for refinement.

The equal low order estimators (Table \ref{table:hartmann_comp_Minit_eta_low}) are underestimating the error with an effectivity of roughly $0.34-0.40$ for the coarsest mesh.
This gets much better after the two refinements where the estimator is close to the exact $L_2$-error. 
However, as the reconstruction effectively works on an even coarser mesh this is not surprising.

The mixed order estimators (Table \ref{table:hartmann_comp_Minit_eta_mixed}) are less dependent on the spatial mesh size 
but overestimate the error with an effectivity of  $1.20-1.48$. We can also see that the overestimation gets smaller with rising $M_{\text{initial}}$ but this 
is of course an additional cost factor.

The equal high order estimators (Table \ref{table:hartmann_comp_Minit_eta_high}) are also less dependent on the spatial mesh size,
but they also use double the amount of degrees of freedom for the primal problem. 
They also overestimate the error but for large enough $M_{\text{initial}}$ the effectivity is less than $1.10$.
However, Table \ref{table:hartmann_comp_Minit_exact21} shows that solving the primal problem with biquadratic elements
and using a bilinear interpolation between the vertices does not recover the best approximation $u_{kh}$.
In practice this means that the primal problem should be solved natively 
with bilinear elements to obtain the 
solution for which the error is actually estimated, which leads to additional costs.
This gets especially expensive for nonlinear problems. Note that in this particular case the resulting adjoint solution and consequently 
the estimators would also be different as the $L_2$ error itself factors into $J'_u$, so the accuracy of the estimator could be better.

In conclusion, both reconstructing a higher order solution and natively solving the primal or adjoint problem with higher order elements in space
work well for linear problems, but the reconstruction is cheaper and leads to a better estimator for fine enough meshes.\\
How the low and mixed order approach perform for higher order finite elements and corresponding errors could be subject of further studies.

\begin{table}[H]
\centering
\caption{Section \ref{sec:config_hartmann_heat}: The exact $L_2$ error under global refinement for different initial (uniform) temporal grids and bilinear finite elements in space.}
\label{table:hartmann_comp_Minit_exact1}
\begin{tabular}{l|l|l|l|l}
\hline
 $N$ & $M_{\text{init}}=100$ & $M_{\text{init}}=400$ & $M_{\text{init}}=800$ & $M_{\text{init}}=1600$ \\
\hline
   $64$ & $1.68717e-02$ & $1.65372e-02$ & $1.64827e-02$ & $1.64555e-02$ \\ 
  $256$ & $4.72445e-03$ & $4.51407e-03$ & $4.48052e-03$ & $4.46392e-03$ \\
 $1024$ & $1.28165e-03$ & $1.16287e-03$ & $1.14457e-03$ & $1.13561e-03$ \\
 $4096$ & $3.67307e-04$ & $3.01477e-04$ & $2.91814e-04$ & $2.87159e-04$ \\
\hline
\end{tabular}
\end{table}

\begin{table}[H]
\centering
\caption{Section \ref{sec:config_hartmann_heat}: The primal equal low order error estimator $\eta_{kh}^{1/1}$ under global refinement for different initial (uniform) temporal grids.}
\label{table:hartmann_comp_Minit_eta_low}
\begin{tabular}{l|l|l|l|l}
\hline
 $N$ & $M_{\text{init}}=100$ & $M_{\text{init}}=400$ & $M_{\text{init}}=800$ & $M_{\text{init}}=1600$ \\
\hline
   $64$ & $6.74336e-03$ & $5.77335e-03$ & $5.65418e-03$ & $5.60003e-03$ \\
  $256$ & $3.66291e-03$ & $3.30266e-03$ & $3.22538e-03$ & $3.18913e-03$ \\ 
 $1024$ & $1.41226e-03$ & $1.10956e-03$ & $1.06963e-03$ & $1.05100e-03$ \\
 $4096$ & $4.84373e-04$ & $3.13938e-04$ & $2.92451e-04$ & $2.82709e-04$ \\ 
\hline
\end{tabular}
\end{table}

\begin{table}[H]
\centering
\caption{Section \ref{sec:config_hartmann_heat}: The primal mixed order error estimator $\eta_{kh}^{1/2}$ under global refinement for different initial (uniform) temporal grids.}
\label{table:hartmann_comp_Minit_eta_mixed}
\begin{tabular}{l|l|l|l|l}
\hline
 $N$ & $M_{\text{init}}=100$ & $M_{\text{init}}=400$ & $M_{\text{init}}=800$ & $M_{\text{init}}=1600$ \\
\hline
   $64$ & $2.05887e-02$ & $1.99047e-02$ & $1.98049e-02$ & $1.97567e-02$ \\
  $256$ & $6.13544e-03$ & $5.64178e-03$ & $5.57640e-03$ & $5.54573e-03$ \\
 $1024$ & $1.74073e-03$ & $1.46754e-03$ & $1.43250e-03$ & $1.41630e-03$ \\
 $4096$ & $5.44978e-04$ & $3.87165e-04$ & $3.67990e-04$ & $3.59415e-04$ \\
\hline
\end{tabular}
\end{table}

\begin{table}[H]
\centering
\caption{Section \ref{sec:config_hartmann_heat}: The primal equal high order error estimator $\eta_{kh}^{2/2}$ under global refinement for different initial (uniform) temporal grids.}
\label{table:hartmann_comp_Minit_eta_high}
\begin{tabular}{l|l|l|l|l}
\hline
 $N$ & $M_{\text{init}}=100$ & $M_{\text{init}}=400$ & $M_{\text{init}}=800$ & $M_{\text{init}}=1600$ \\
\hline
   $64$ & $1.82979e-02$ & $1.77462e-02$  & $1.76630e-02$  & $1.76224e-02$ \\
  $256$ & $5.38291e-03$ & $4.95363e-03$  & $4.89527e-03$  & $4.86777e-03$ \\
 $1024$ & $1.52650e-03$ & $1.28219e-03$  & $1.25015e-03$  & $1.23518e-03$ \\
 $4096$ & $4.80895e-04$ & $3.38434e-04$  & $3.20913e-04$  & $3.12975e-04$ \\
\hline
\end{tabular}
\end{table}

\begin{table}[H]
\centering
\caption{Section \ref{sec:config_hartmann_heat}: The approximated $L_2$ error under global refinement for different initial (uniform) temporal grids 
and biquadratic finite elements in space, where the solution is interpolated down to bilinear elements.}
\label{table:hartmann_comp_Minit_exact21}
\begin{tabular}{l|l|l|l|l}
\hline
 $N$ & $M_{\text{init}}=100$ & $M_{\text{init}}=400$ & $M_{\text{init}}=800$ & $M_{\text{init}}=1600$ \\
\hline
   $64$ & $2.16806e-02$ & $2.13230e-02$ & $2.12638e-02$ & $2.12343e-02$ \\ 
  $256$ & $6.38393e-03$ & $6.17212e-03$ & $6.13780e-03$ & $6.12076e-03$ \\
 $1024$ & $1.72805e-03$ & $1.61279e-03$ & $1.59468e-03$ & $1.58575e-03$ \\
 $4096$ & $4.78600e-04$ & $4.16052e-04$ & $4.06646e-04$ & $4.02072e-04$ \\
\hline
\end{tabular}
\end{table}

\subsubsection{Comparison of adaptive refinement to the original computations}
For this configuration we want to compare our results 
to those described by Hartmann \cite{hartmann1998},
where the manufactured solution was 
first formulated. To our knowledge this configuration was not reproduced and published so far, so the original thesis is the only point of comparison.
There, the classical estimator is used, which is obtained by partial integration to obtain a strong form with jump terms in space.
There, $Q_1$ elements in space and $dG(0)$ elements in time are used as well, but it is unknown which quadrature formula was used in time 
for the nonlinear $f$. We used the right box rule as this is what corresponds to the implicit Euler scheme and got our error ($1.92e-02$) closest to 
the error of the original results ($1.75e-02$).
Table \ref{table:hartmann_adaptive_thiwi} 
shows our results with the split and joint estimators. These 
perform very well in comparison to the original Hartmann results from
\cite[Table 3.4]{hartmann1998}.

\begin{figure}[H]
 \centering
  \includegraphics[width=0.8\textwidth,page = 1]{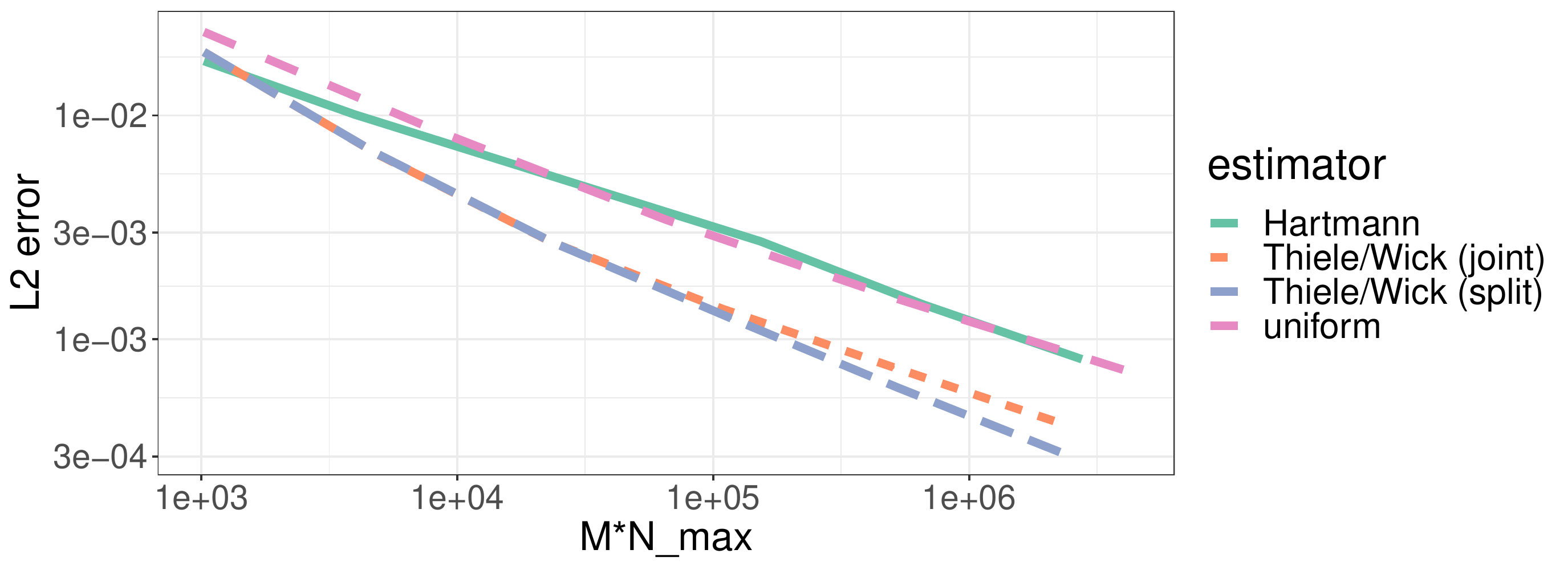}
 \caption{Section \ref{sec:config_hartmann_heat}: Error convergence of the Hartmann testcase}
\label{figure:hartmann_convergence_comparison}
\end{figure}

\begin{figure}[H]
 \centering
 
 \subfloat{
  \includegraphics[clip, trim = 362 50 390 53, width=0.22\textwidth]{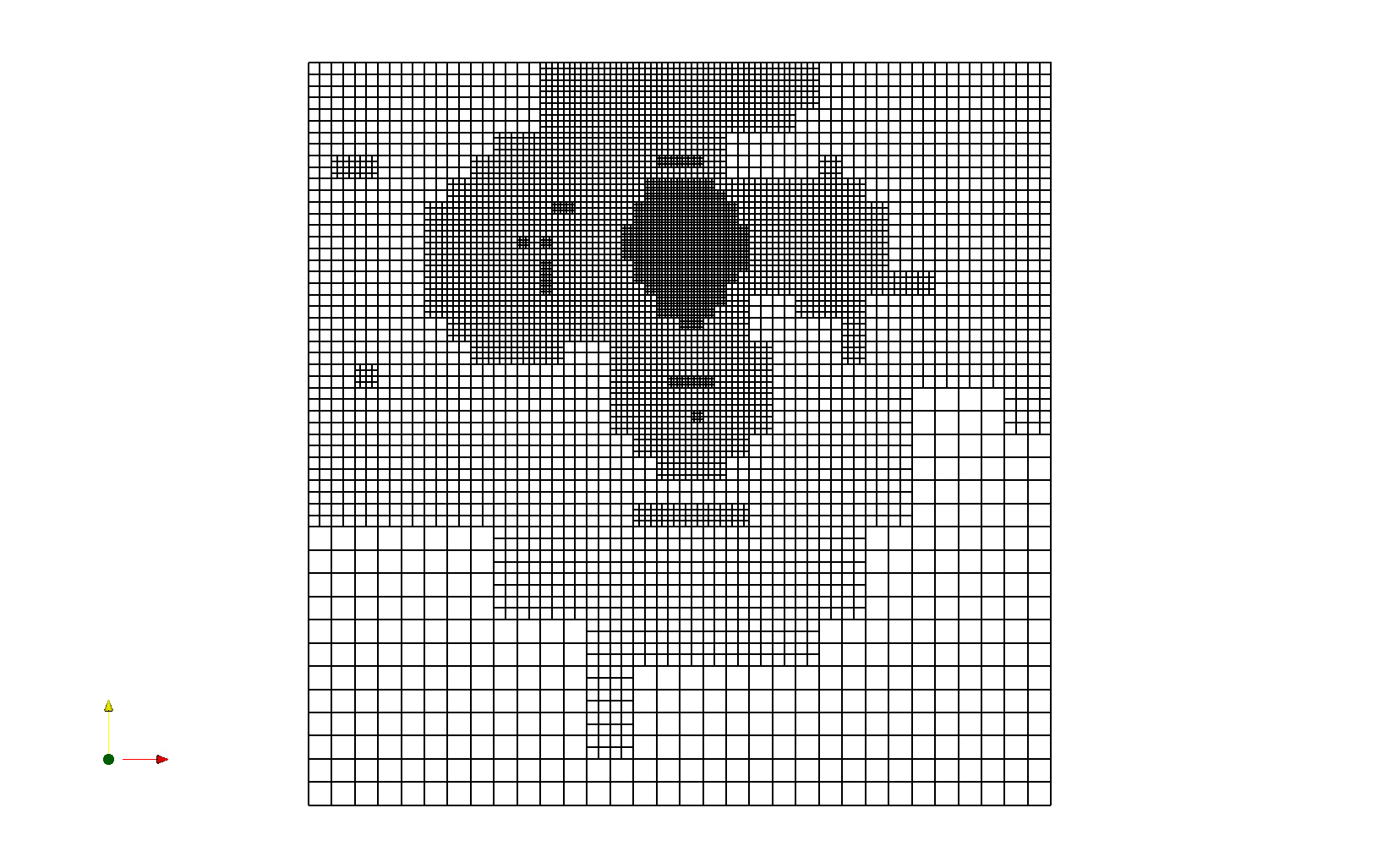}
  }
 \subfloat{
   \includegraphics[clip, trim = 362 50 390 53, width=0.22\textwidth]{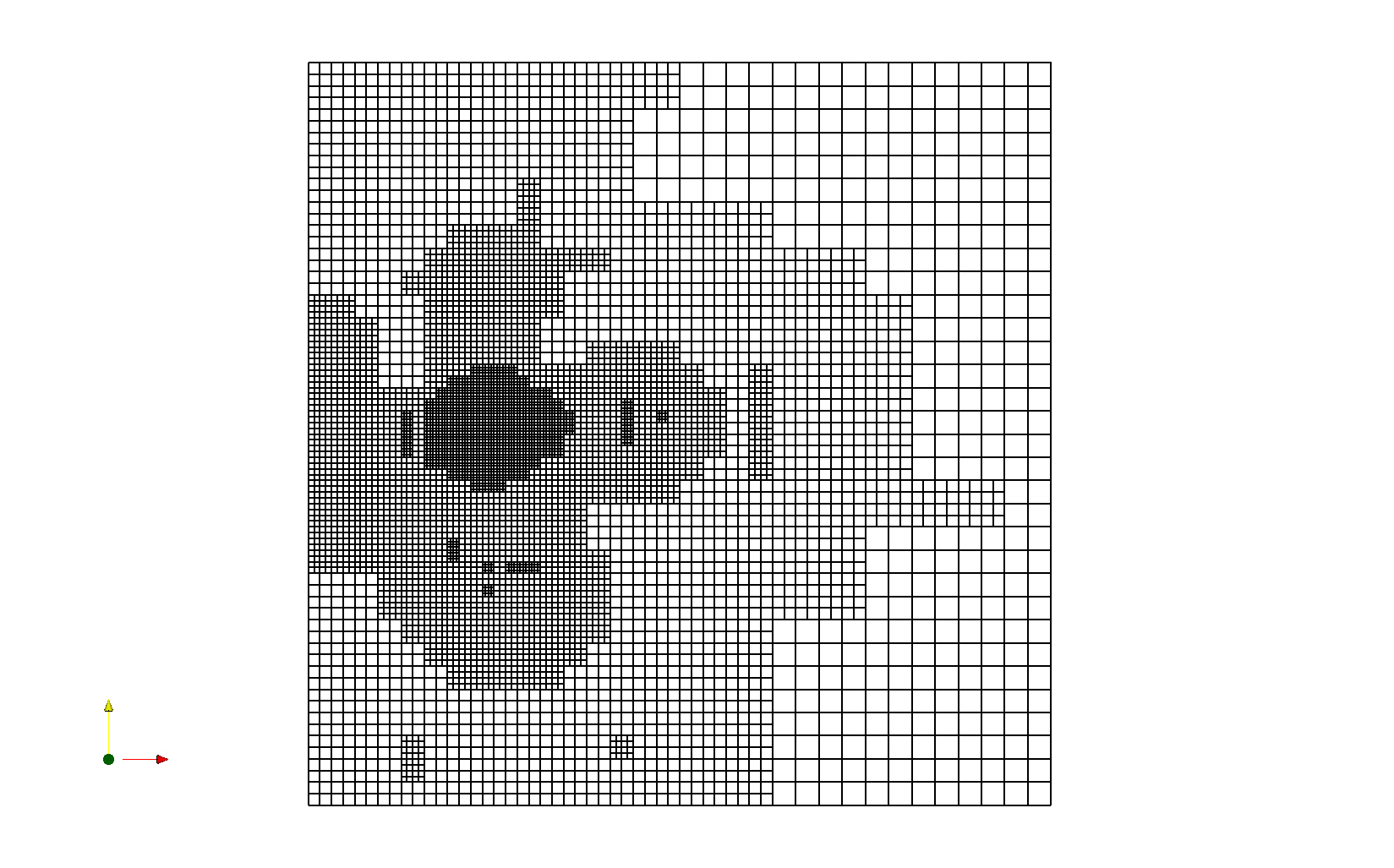}
 }
 \subfloat{
   \includegraphics[clip, trim = 362 50 390 53, width=0.22\textwidth]{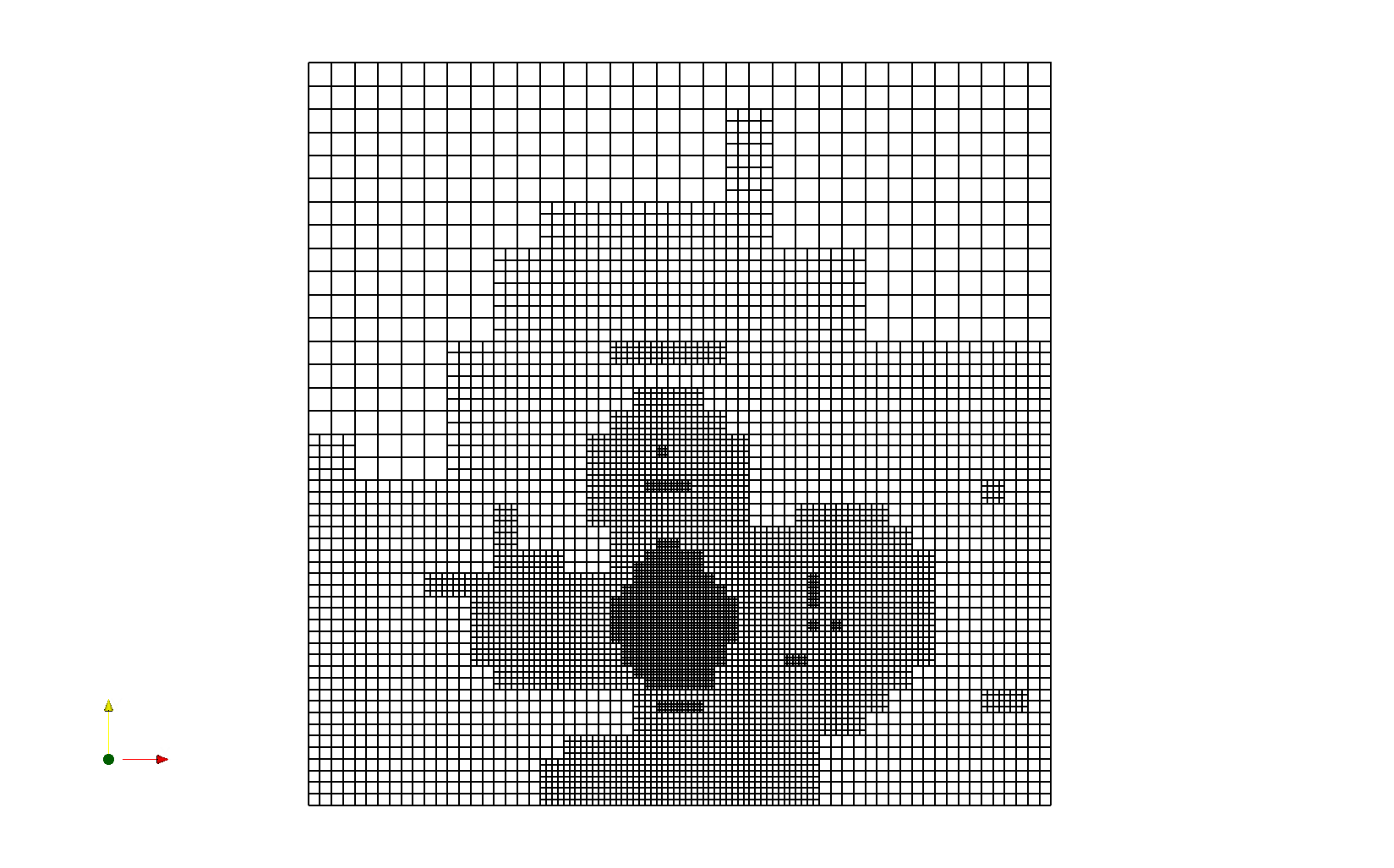}
 }
 \subfloat{
   \includegraphics[clip, trim = 362 50 390 53, width=0.22\textwidth]{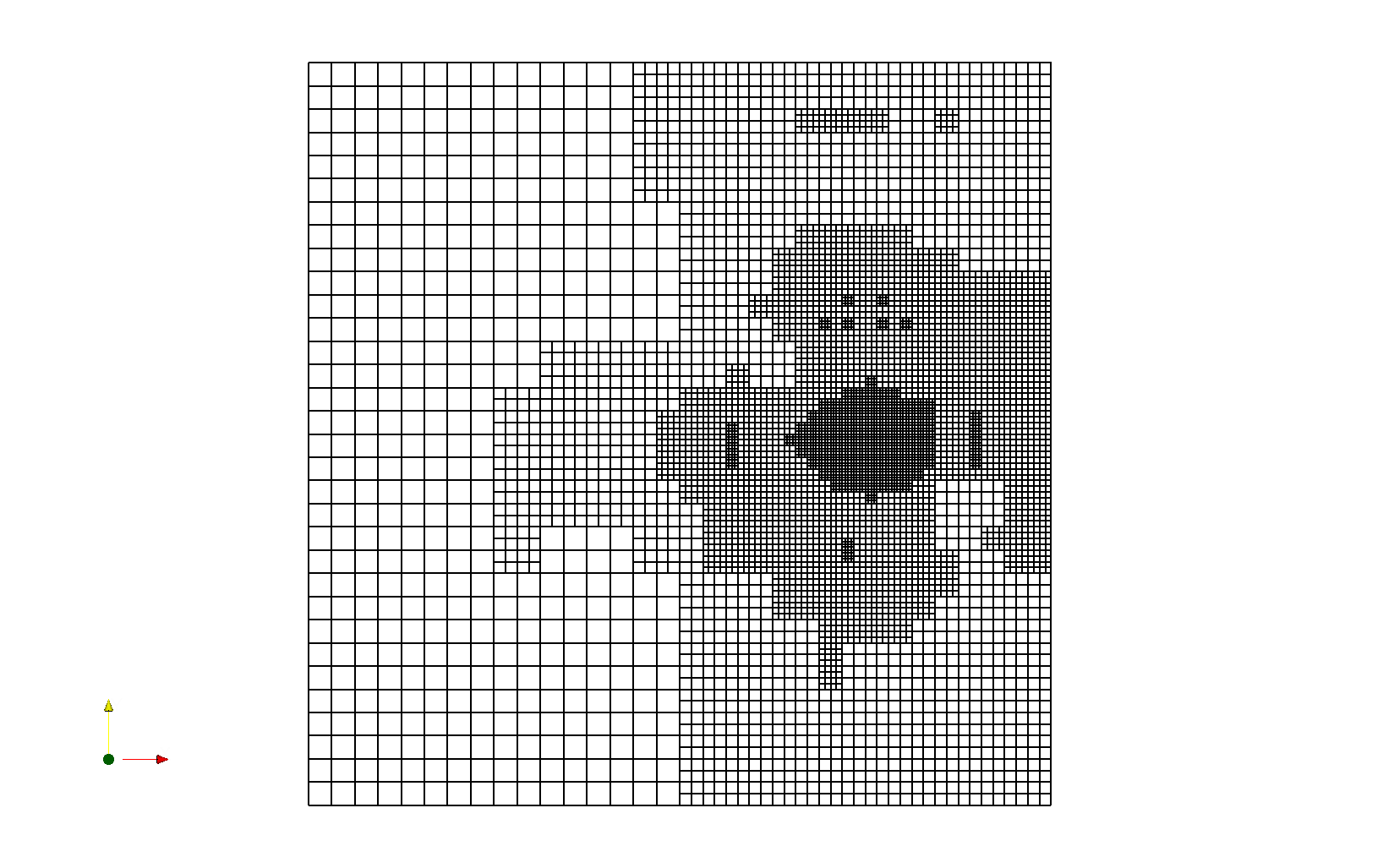}
 }
 \caption{Section \ref{sec:config_hartmann_heat}: Grid after 4 refinement loops with the split PU-DWR estimator at $t=i/4$, $i\in\{1,2,3,4\}$}
\label{figure:Hartmann_meshes}
\end{figure}

\begin{figure}[H]
 \centering
 \subfloat{
  \includegraphics[clip, trim = 275 0 460 0, width=0.22\textwidth]{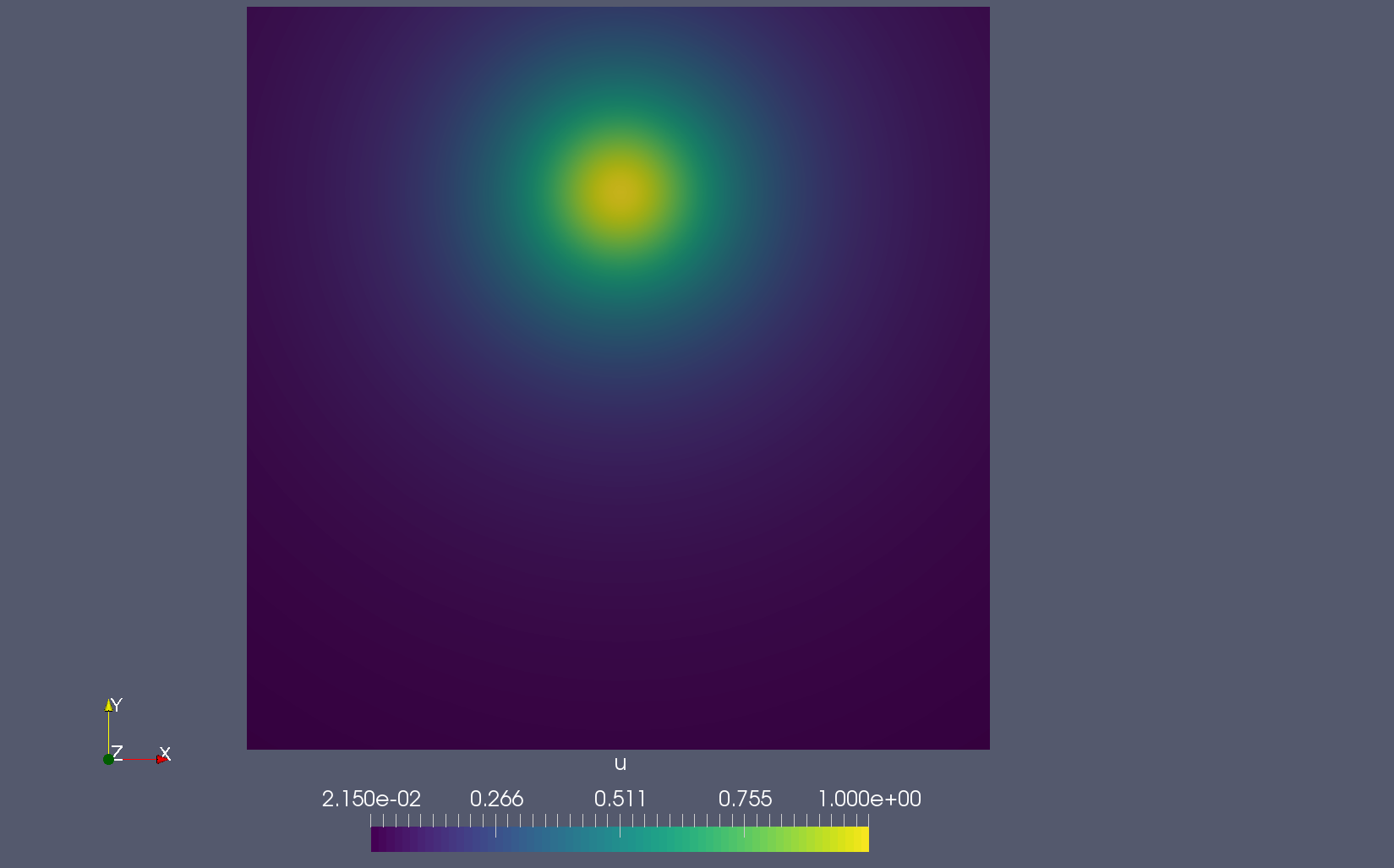}
  }
 \subfloat{
   \includegraphics[clip, trim = 275 0 460 0, width=0.22\textwidth]{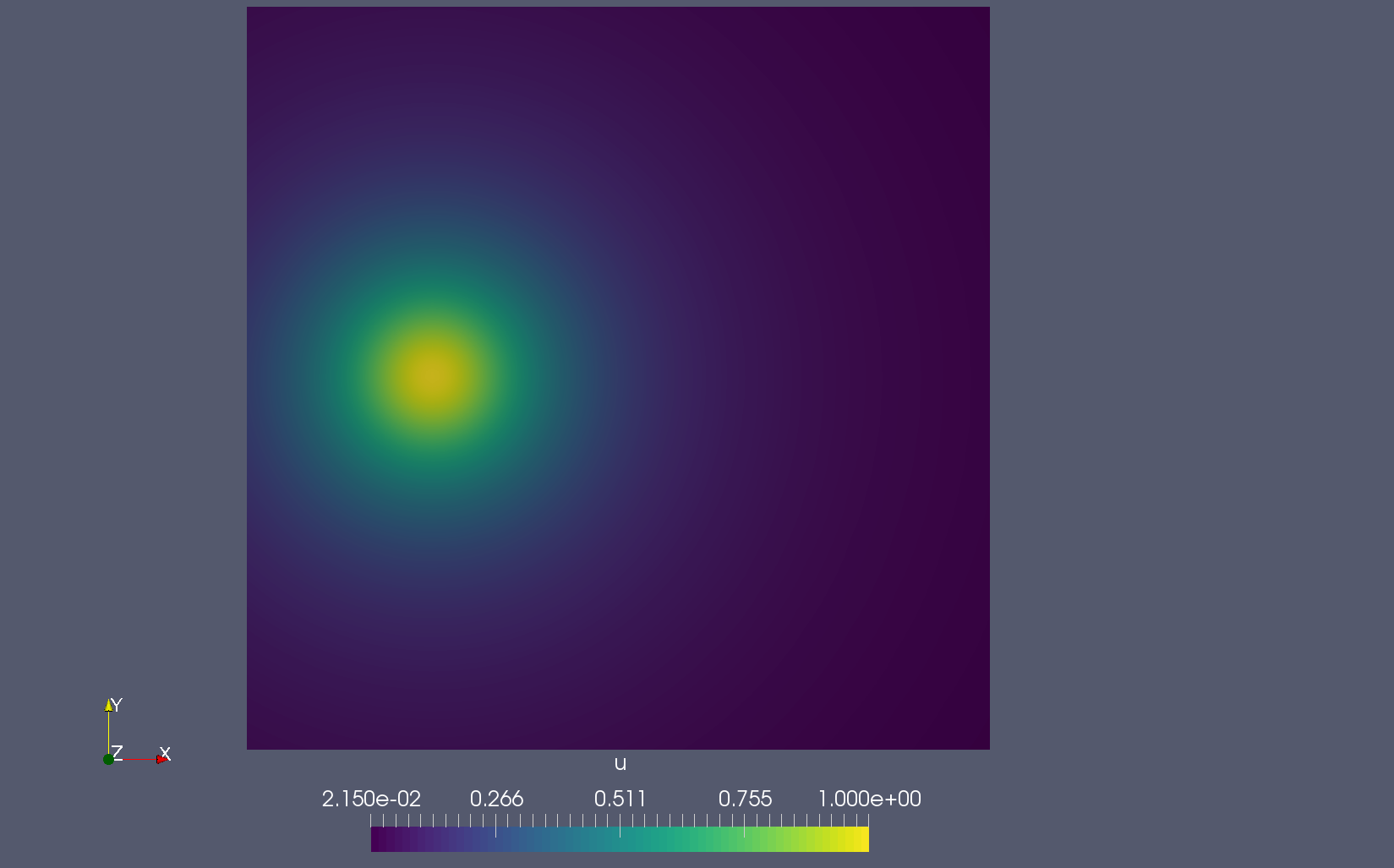}
 }
 \subfloat{
   \includegraphics[clip, trim = 275 0 460 0, width=0.22\textwidth]{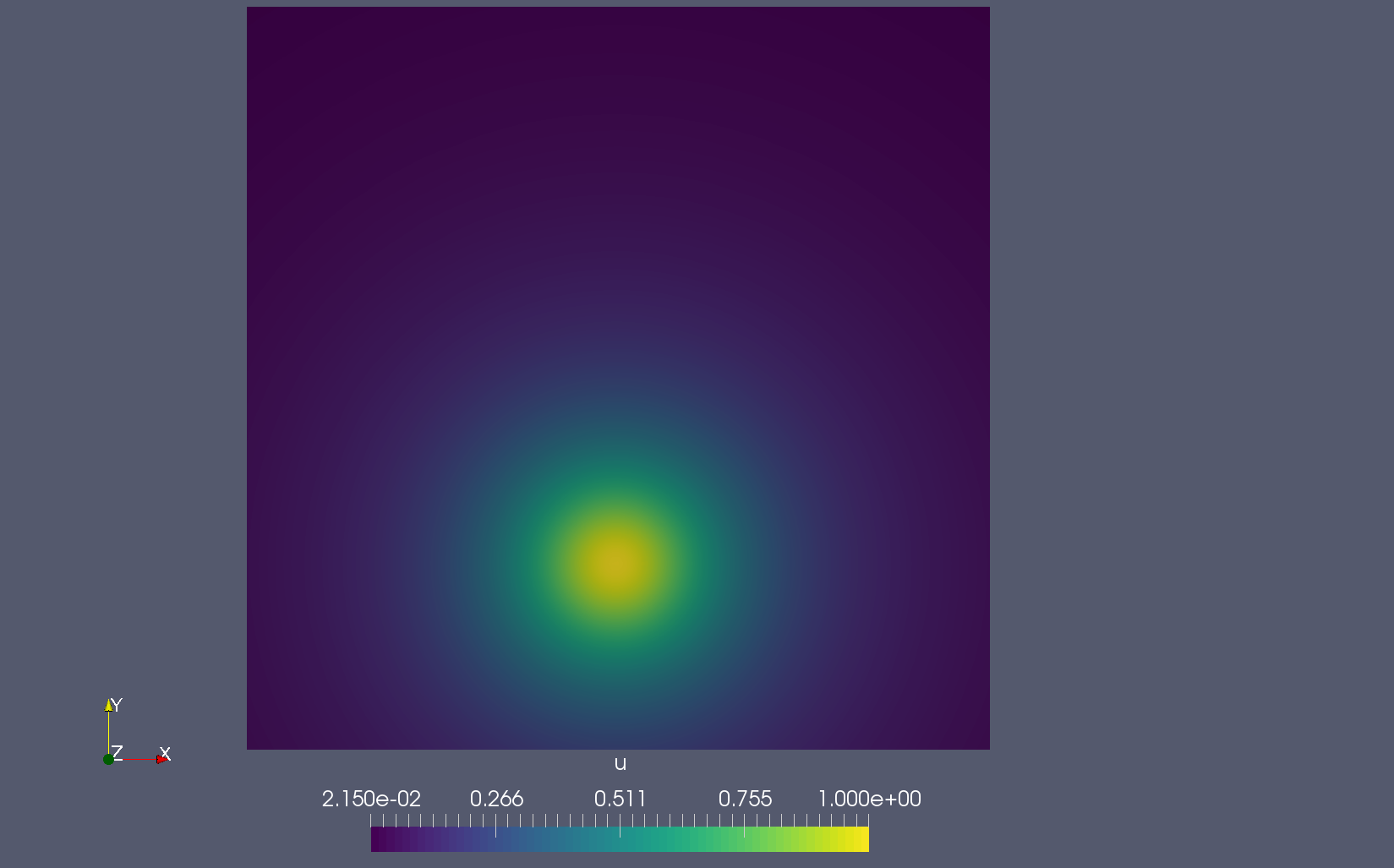}
 }
 \subfloat{
   \includegraphics[clip, trim = 275 0 460 0, width=0.22\textwidth]{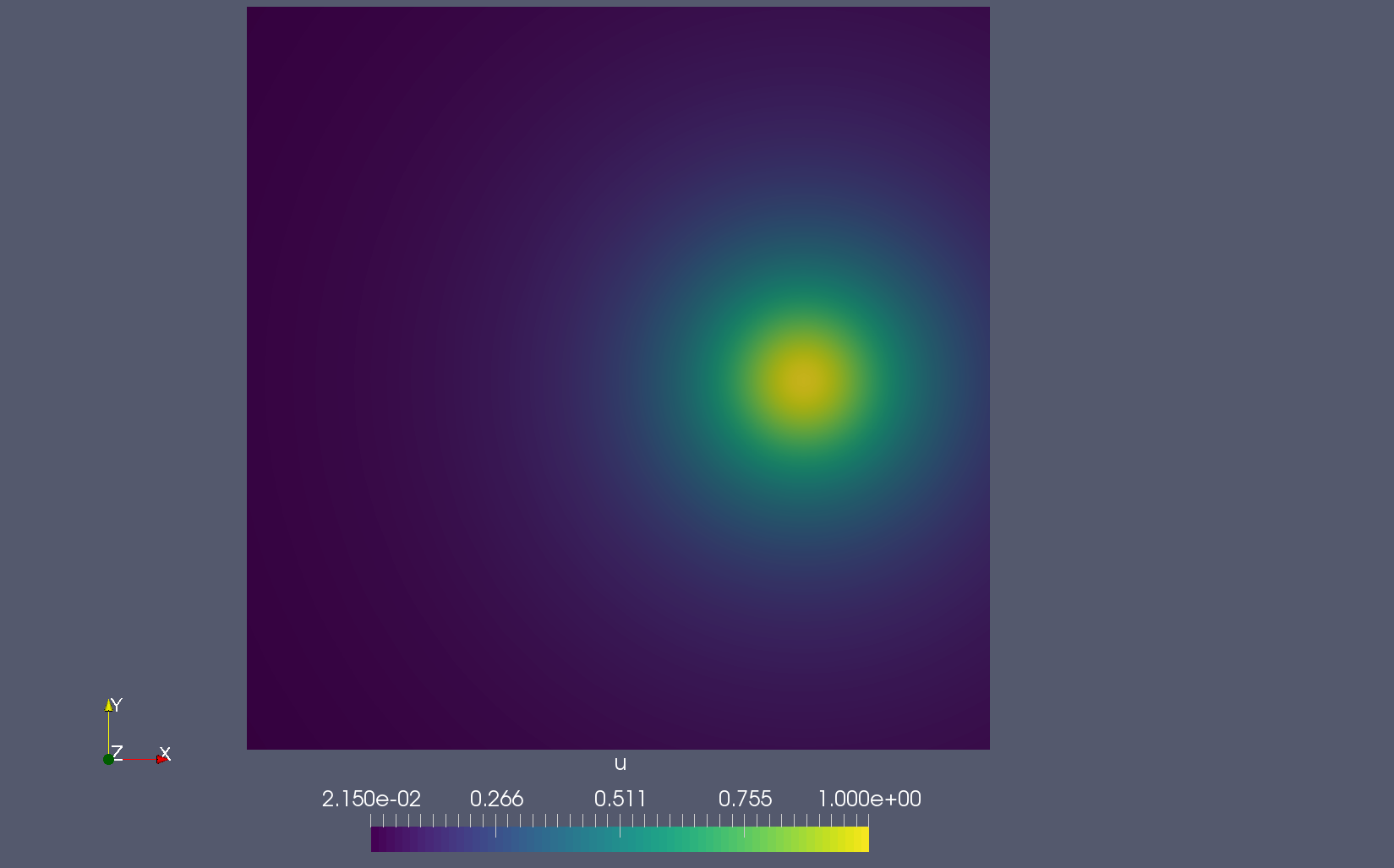}
 }
 \caption{Section \ref{sec:config_hartmann_heat}: Solution after 4 refinement loops with the split PU-DWR estimator at $t=i/4$, $i\in\{1,2,3,4\}$}
\label{figure:Hartmann_solution}
\end{figure}

Even though the marking strategy used by Hartmann is unknown, we got close to the number of temporal elements and the maximum
number of spatial elements with fixed rate marking in Algorithm \ref{algo_mark_refine_timestep}.
In comparison, our estimators better localize the error as we get comparable errors with one less loop that additionally
has a smaller maximum number of spatial elements. This can be seen in Figure \ref{figure:hartmann_convergence_comparison},
where the original results are only performing roughly as well as our computation with uniform refinement,
while both PU-DWR estimators yield better convergence. 
We plotted the $L_2$ error against $M*N_{\max}$ which is an upper bound for the actual number of space-time elements
as no further information was available from the original computations. However, Table \ref{table:hartmann_adaptive_elements}
shows that at least for our simulations the actual number is not too far from the upper bound. 
Additionally, we can see that, as expected, the split estimator outperforms the joint estimator.\\
Finally, Figure \ref{figure:Hartmann_meshes} shows that the local refinement nicely matches the corresponding solution 
of Figure \ref{figure:Hartmann_solution} and that the meshes are indeed changing over time.

\begin{table}[H]
\centering
\caption{Section \ref{sec:config_hartmann_heat}: Our results with the equal low order primal split (left) and joint (right) 
PU-DWR estimator with fixed rate marking of $95\%$ in time and $40\%$ in space.}
\label{table:hartmann_adaptive_thiwi}
\begin{tabular}{l|l|l|l|l}
 \hline
 $M$ & $N_{\max}$ & $||e||_{L_2(\Sigma)}$ & $\eta^{1/1}_{split}$ & $I_{\text{eff}}$ \\
 \hline 
  $16$ &   $64$ & $1.92e-02$ & $2.57e-02$ & $1.34$ \\
  $32$ &  $151$ & $6.77e-03$ & $7.64e-03$ & $1.13$ \\
  $64$ &  $364$ & $2.73e-03$ & $3.32e-03$ & $1.22$ \\
 $127$ &  $898$ & $1.26e-03$ & $1.57e-03$ & $1.24$ \\
 $250$ & $2110$ & $6.08e-04$ & $9.24e-04$ & $1.52$ \\ 
 $490$ & $4840$ & $3.07e-04$ & $6.15e-04$ & $2.00$ \\
 \hline
\end{tabular}
\qquad
\begin{tabular}{l|l|l|l|l}
 \hline
 $M$ & $N_{\max}$ & $||e||_{L_2(\Sigma)}$ & $\eta^{1/1}_{joint}$ & $I_{\text{eff}}$ \\
 \hline 
  $16$ &   $64$ & $1.92e-02$ & $4.15e-03$ & $0.22$ \\
  $32$ &  $154$ & $6.67e-03$ & $3.92e-03$ & $0.59$ \\
  $64$ &  $367$ & $2.71e-03$ & $3.32e-03$ & $1.22$ \\
 $126$ &  $859$ & $1.37e-03$ & $1.60e-03$ & $1.17$ \\
 $248$ & $1987$ & $7.58e-04$ & $9.15e-04$ & $1.21$ \\ 
 $486$ & $4543$ & $4.27e-04$ & $5.09e-04$ & $1.19$ \\
 \hline
\end{tabular}
\end{table}

\begin{table}[H]
\centering
\caption{Section \ref{sec:config_hartmann_heat}: Comparison of actual number of space time elements and estimation by $M*N_{\max}$}
\label{table:hartmann_adaptive_elements}
\begin{tabular}{l|l|l|l|l}
 \hline
 $loop$ & $M * N_{\max} (split)$ & $\#elements (split)$ & $M * N_{\max} (joint) $ & $\#elements (joint)$ \\
 \hline 
 $0$ &    $1024$ &    $1024$ &    $1024$ &    $1024$ \\
 $1$ &    $4832$ &    $4646$ &    $4928$ &    $4760$ \\
 $2$ &   $23296$ &   $21910$ &   $23488$ &   $22144$ \\
 $3$ &  $114046$ &  $104962$ &  $108234$ &  $100512$ \\
 $4$ &  $527500$ &  $486745$ &  $492776$ &  $451433$ \\
 $5$ & $2371600$ & $2186416$ & $2207898$ & $1993332$ \\
 \hline
\end{tabular}
\end{table}

\subsubsection{Comparison of PU spaces}
\begin{figure}[H]
 \centering
  \includegraphics[width=0.8\textwidth,page = 1]{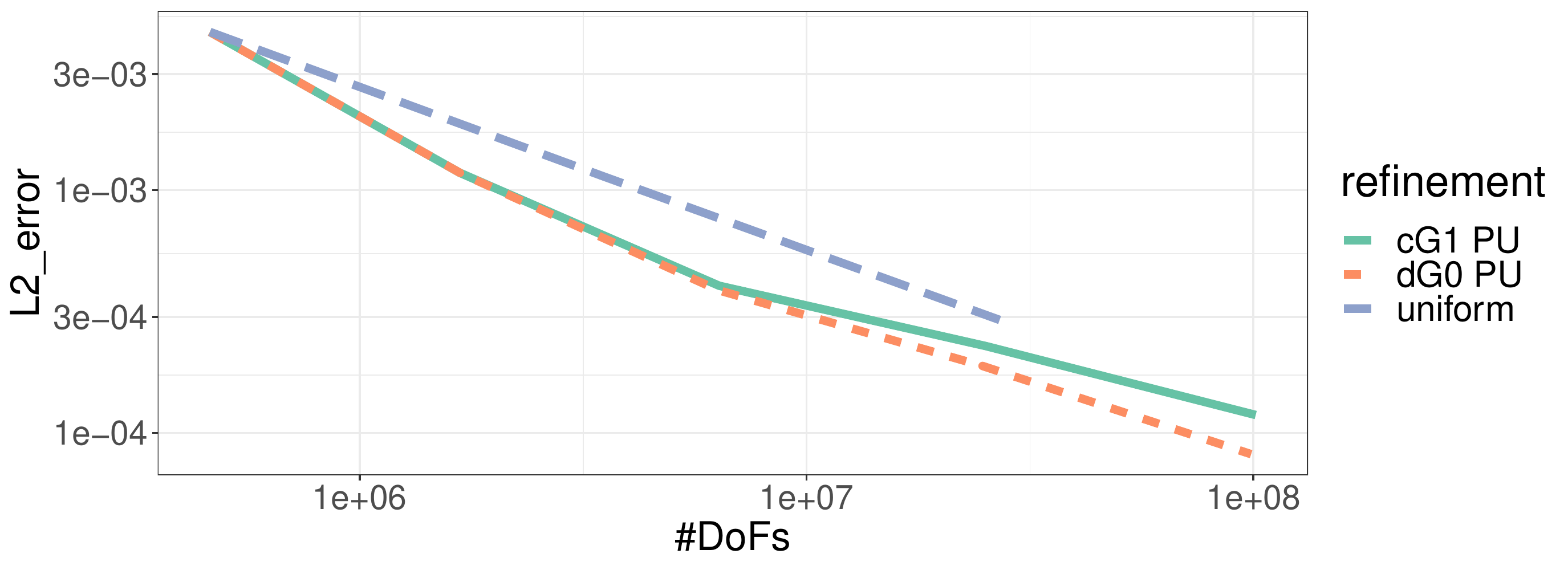}
 \caption{Section \ref{sec:config_hartmann_heat}: Error convergence of the Hartmann testcase for $\eta^{1/1}$ and $M_\text{init} = 1600$}
\label{figure:hartmann_pu_comparison_lo}
\end{figure}

\begin{figure}[H]
 \centering
  \includegraphics[width=0.8\textwidth,page = 2]{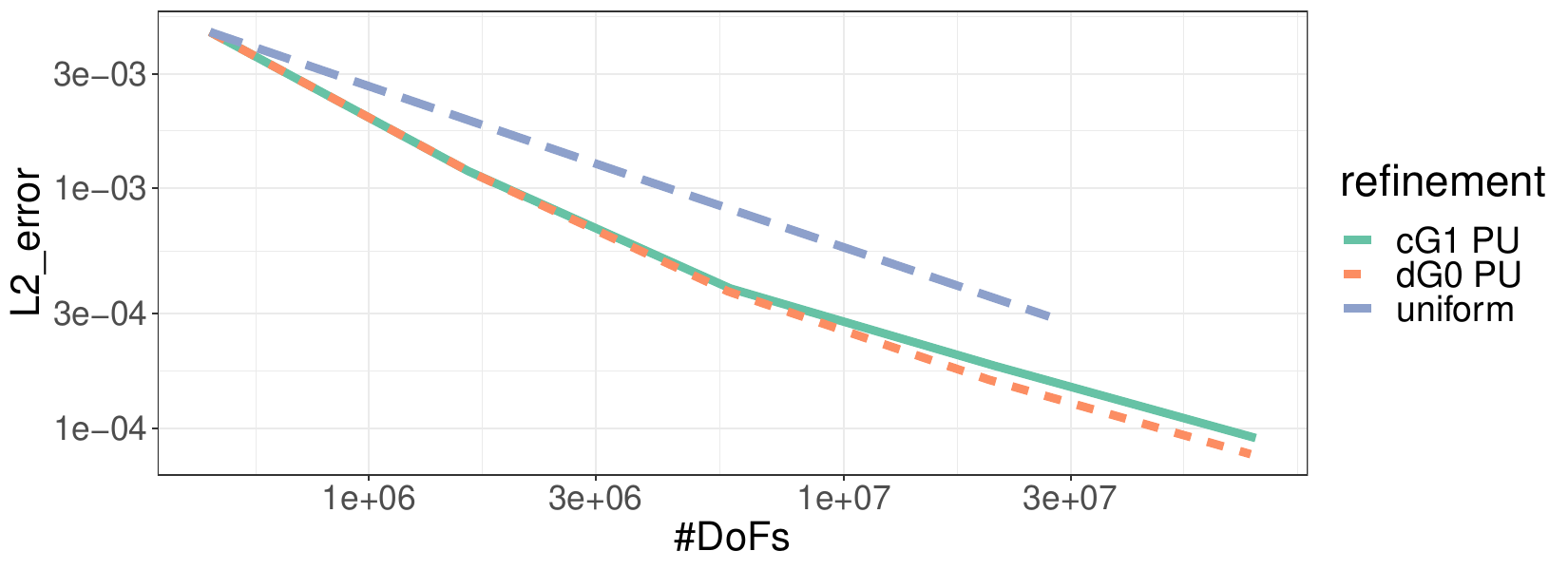}
 \caption{Section \ref{sec:config_hartmann_heat}: Error convergence of the Hartmann testcase for $\eta^{1/2}$ and $M_\text{init} = 1600$}
\label{figure:hartmann_pu_comparison_mo}
\end{figure}

\begin{figure}[H]
 \centering
  \includegraphics[width=0.8\textwidth,page = 3]{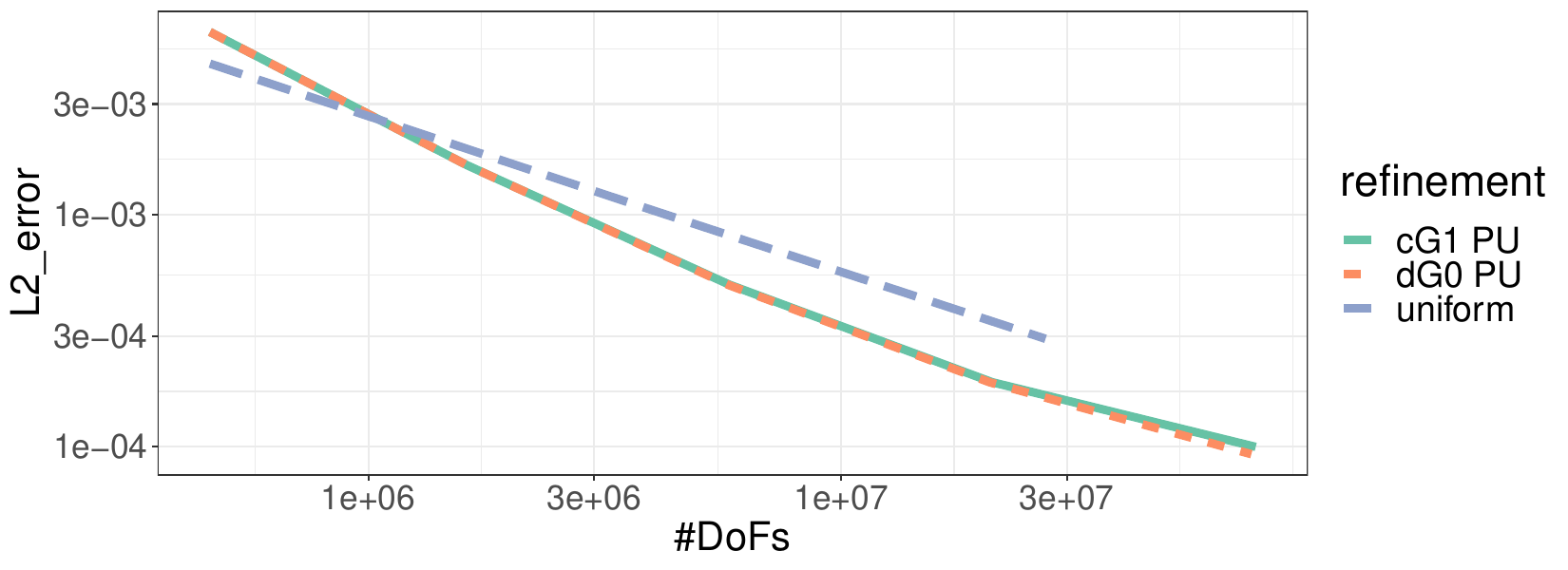}
 \caption{Section \ref{sec:config_hartmann_heat}: Error convergence of the Hartmann testcase for $\eta^{2/2}$ and $M_\text{init} = 1600$}
\label{figure:hartmann_pu_comparison_ho}
\end{figure}

\begin{figure}[H]
 \centering
  \includegraphics[width=0.8\textwidth,page = 4]{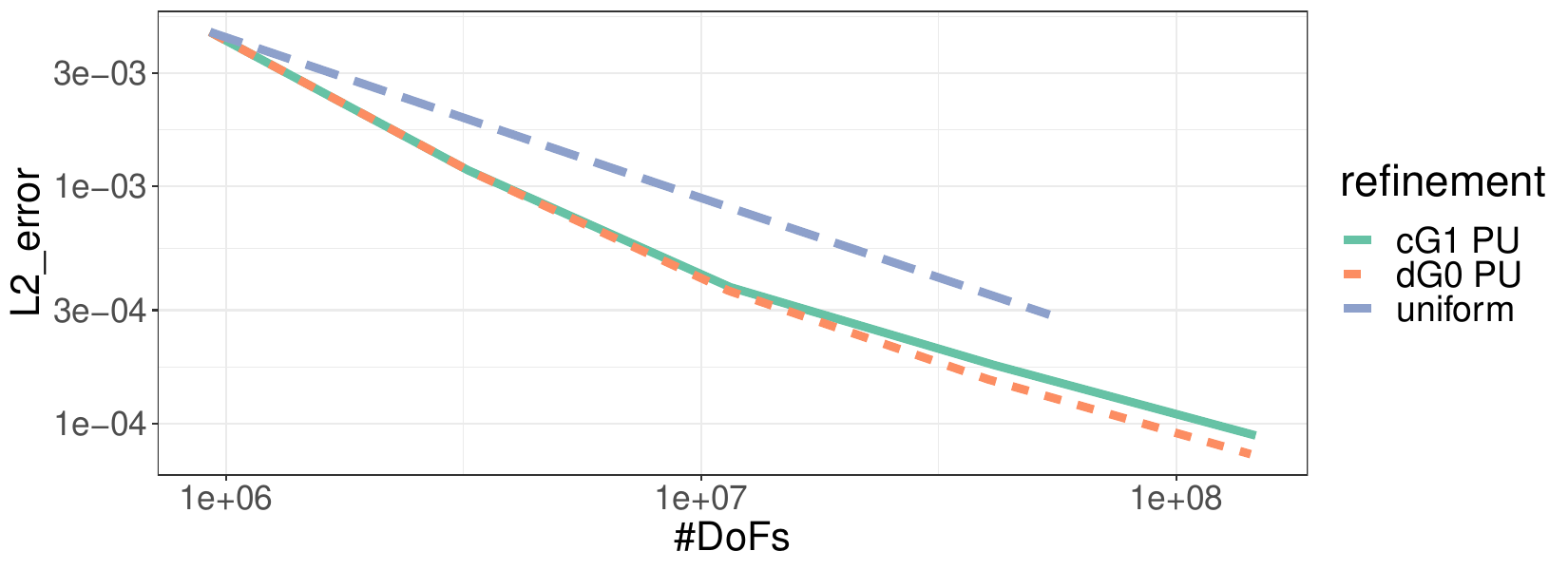}
 \caption{Section \ref{sec:config_hartmann_heat}: Error convergence of the Hartmann testcase for $\eta^{1/2}$ and $M_\text{init} = 3200$}
\label{figure:hartmann_pu_comparison_mo_2M}
\end{figure}

Here, we want to examine whether the additional coupling with a $cG(1)$ partition-of-unity in time is beneficial for adaptive refinement.
For this, we performed multiple adaptive simulations with the corresponding low, mixed and high order estimators.
Figures \ref{figure:hartmann_pu_comparison_lo}, \ref{figure:hartmann_pu_comparison_mo} and \ref{figure:hartmann_pu_comparison_ho} 
show the exemplary results with an initial temporal mesh of $1600$ elements with fixed rate marking of $60\%$ in time and $40\%$ in space.
In all three cases the $dG(0)$ PU performs better than the $cG(1)$ PU, with the lowest difference for the equal high order estimator.
For the equal high order case the $L^2$-error is calculated on the interpolated primal solution, which does not recover the actual
best approximation solution of directly solving the primal solution in the low order space, such that the initial error is higher than 
for uniform refinement.
Finally, Figure \ref{figure:hartmann_pu_comparison_mo_2M} shows the results for the mixed order estimator with a finer initial 
temporal mesh, which are qualitatively the same as for Figure \ref{figure:hartmann_pu_comparison_mo}.\\
Overall we conclude that for discontinuous Galerkin discretizations in time the $dG(0)$ partition-of-unity is completely sufficient in 
addition to being cheaper to compute and easier to implement.

\subsection{Configuration 3: nonlinear combustion}
\label{sec:config_combustion}
\subsubsection{Problem statement}
The final test case is as described in \cite{schmich.vexler2008}
(originally based on \cite{Lang2001})
and some preliminary results were published in our prior 
work \cite{ThiWi23_ICO}. Here, we solve the nonlinear combustion equations described in Section \ref{sec_combustion}.

\subsubsection{Configuration}
The reaction is simulated in a rectangular channel of length $L=60$ and height $H=16$ in which
two cooled rods of length $L/4$ and height $H/4$ are inserted into both channel walls at $L/4$.
The reaction is solved for a total of $T=60$ with $256$ time and $896$ space DoFs initially.\\
The cooling of $\Gamma_R$ is described by the Robin boundary condition $\partial_n\theta = -0.1\theta$,
with homogeneous Neumann conditions for the species concentration.
The left wall $\Gamma_D$ is kept at a constant temperature of $\theta_D = 1$ without any combustible species $Y_D = 0$. 
All other walls $\Gamma_N$ are described by homogeneous Neumann boundary conditions.
An initial flame front is described by
\begin{align}
 \theta^0 = \begin{cases} 
              1,& x\leq 9\\
              \exp(9-x),& x > 9
            \end{cases}\\
 Y^0 = \begin{cases}
         0, & x\leq 9\\
         1-\exp(Le(9-x)), & x > 9 .
       \end{cases}
\end{align}
\subsubsection{Parameters}
The reaction parameters are a Lewis number of $Le = 1$, a gas expension of $\alpha = 0.8$ and a dimensionless energy of $\beta = 10$.
\subsubsection{Goal Functionals}
The first functional we investigate is the average reaction rate i.e. 
\begin{equation}
 J_1(\theta,Y) = \frac{1}{T|\Omega|} \int\limits_0^T\int\limits_\Omega \omega(\theta,Y) \mathrm{d}x\mathrm{d}t.
\end{equation}
This nonlinear functional is defined on the whole space-time domain.
For the second functional we calculate the average species concentration on the cooled rods $\Gamma_R$ i.e.
\begin{equation}
 J_2(\theta,Y) = \frac{1}{T|\Gamma_R|} \int\limits_0^T\int\limits_{\Gamma_R} Y \mathrm{d}s\mathrm{d}t.
\end{equation}
This is a linear functional, but it is only defined on part of the boundary.

\subsubsection{Discussion of findings for $J_1$}
For both functionals the indicators for 
$\eta_{kh}^{{s/\widetilde{s}}}$ are computed using 
Propositions \ref{prop:split_primal_combustion_indicator}
and \ref{prop:split_adjoint_combustion_indicator} and Algorithm \ref{algo_estimate_timestep} with 
$\widetilde u$, $u_k$, $u_{kh}$ and 
$\widetilde z$, $z_k$, $z_{kh}$ following
from \eqref{interpolations_kh} and 
\eqref{interpolations_k} with
$\hat{u}_{kh}\in \widetilde{X}^{0,s}_{k,h}$ 
and $\hat{z}_{kh}\in \widetilde{X}^{0,\widetilde s}_{k,h}$. 
Then, the full estimators are computed by taking the averages of the primal and adjoint indicators at each 
space-time Dof and summing over all DoFs.

Tables \ref{table:comb_glob_prim_reaction}, \ref{table:comb_glob_adj_reaction} 
and \ref{table:comb_glob_full_reaction} show the behaviour of the primal, adjoint and full estimators for $J_1$ respectively.
We can see that all estimators with the mixed order and equal high order approach behave similarly. 
For both, the adjoint estimator and the resulting full estimator are overestimating the error by about two orders of magnitude.
However, the primal estimators are not too far off.
The equal low order approach yields the best results on the finest level, and the primal and adjoint estimators are comparable.
It can be inferred that there is no benefit from calculating the full estimator in this case. 
Therefore, we use Algorithm \ref{algo_mark_refine_timestep} with fixed rate marking refining $50\%$ of all
temporal elements and on each interval $30\%$ of all spatial elements based only on the primal indicators.

\begin{figure}[H]
 \centering
 \subfloat{
  \includegraphics[width=0.48\textwidth]{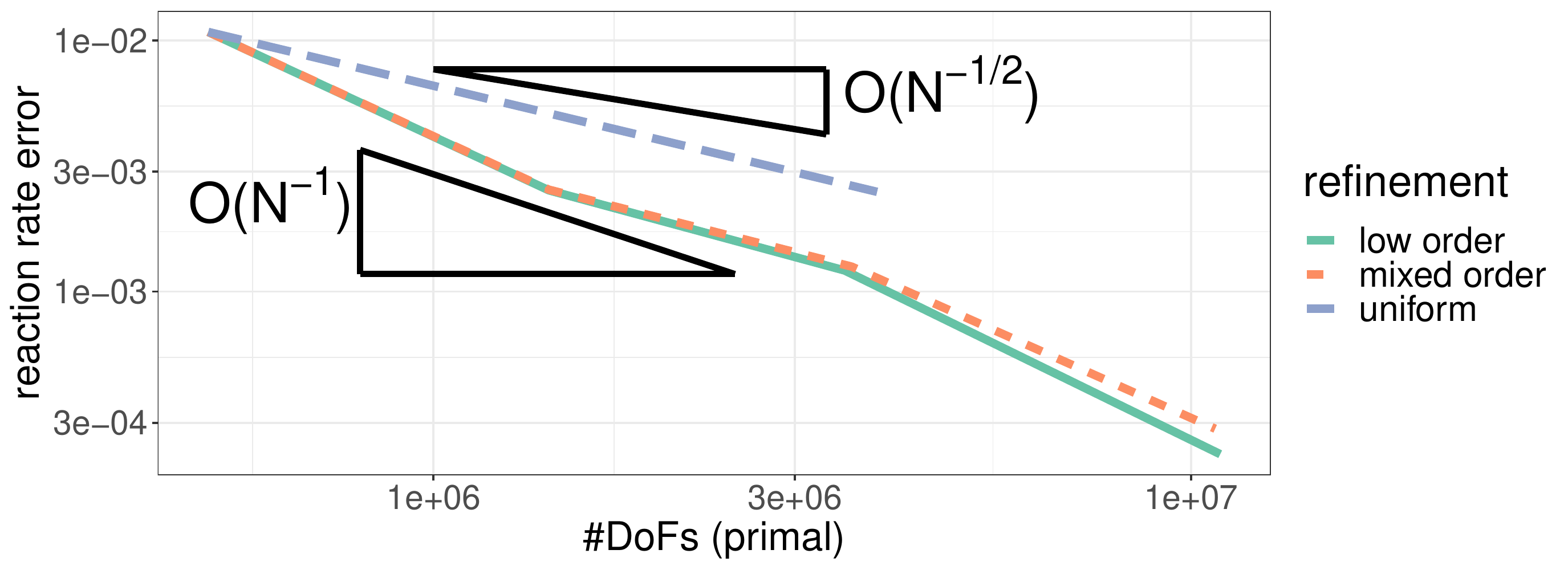}
 } 
 \subfloat{
  \includegraphics[width=0.48\textwidth]{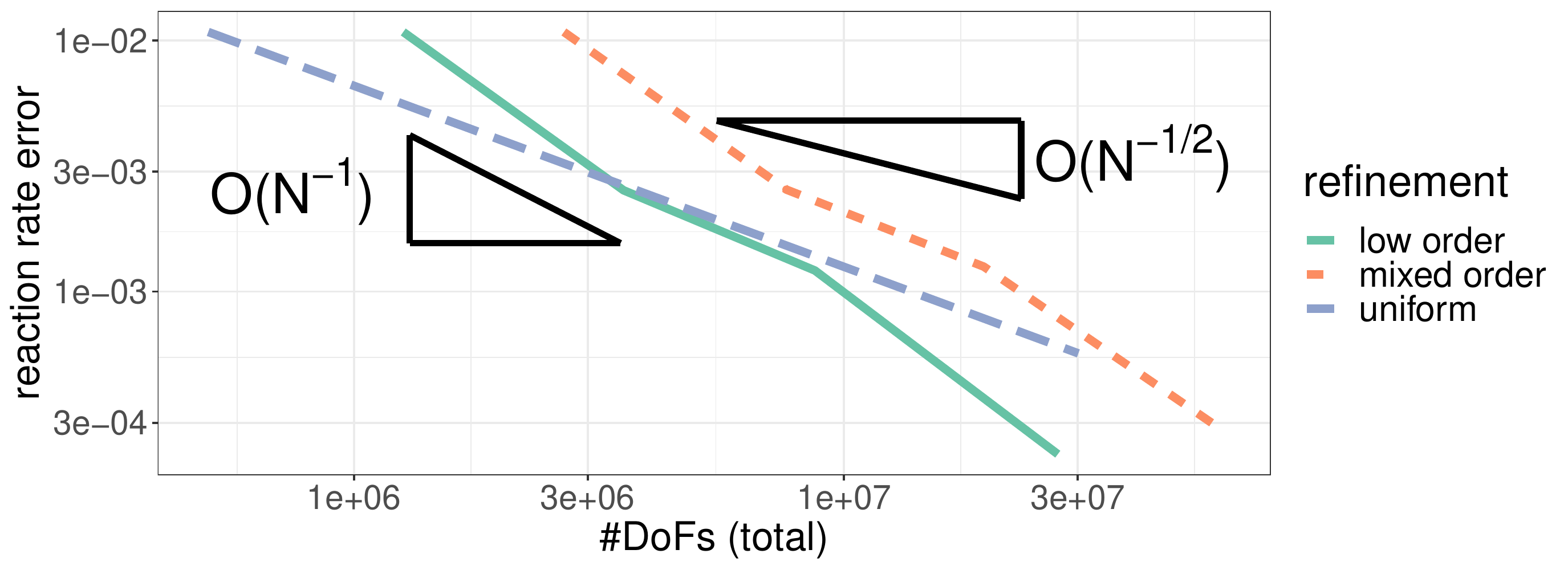}
 }
 \caption{Section \ref{sec:config_combustion}: Error convergence for the reaction rate functional. On the left only the number of unknowns for the primal problem
 and on the right all unknowns are taken into account.}
\label{figure:comb_omega_convergence}
\end{figure}

Figure \ref{figure:comb_omega_convergence} shows the error convergence under adaptive refinement when using only the primal estimator.
When only counting the primal unknowns for both estimators, our findings for equal low order and mixed order 
are comparable, and both perform better than global refinement. 
However, when both the number of Dofs of the adjoint and the PU are additionally taken into account
the low order approach clearly outperforms the mixed order approach. We can also see that the starting disadvantage 
of having the same error on the coarsest mesh with more DoFs is rectified by the first adaptive refinement step.

\begin{table}[H]
\centering
\caption{Section \ref{sec:config_combustion}: Performance of the primal error estimators under global refinement for $J_1$}
\label{table:comb_glob_prim_reaction}
\begin{tabular}{l|l|l|l|l|l}
\hline
 $M$ & $N$ & $J(u)-J(u_{kh})$ & $\eta^{1/1}_{kh}$ & $\eta^{1/2}_{kh}$ & $\eta^{2/2}_{kh}$ \\
\hline
  $256$ &    $896$ & $1.08154812e-02$ & $9.99898726e-04$ & $1.54141116e-03$ & $3.73046753e-03$ \\
  $512$ &   $3584$ & $2.49545713e-03$ & $5.00929317e-04$ & $1.13434333e-03$ & $1.33777161e-03$ \\
 $1024$ &  $14366$ & $5.67139257e-04$ & $2.17940341e-04$ & $5.88255119e-04$ & $5.98302490e-04$ \\
 $2048$ &  $57344$ & $1.11745245e-04$ & $8.34682775e-05$ & $3.46300197e-04$ & $3.45416860e-04$ \\
\hline
\end{tabular}
\end{table}

\begin{table}[H]
\centering
\caption{Section \ref{sec:config_combustion}: Performance of the adjoint error estimators under global refinement for $J_1$}
\label{table:comb_glob_adj_reaction}
\begin{tabular}{l|l|l|l|l|l}
\hline
 $M$ & $N$ & $J(u)-J(u_{kh})$ & $\eta^{1/1}_{kh}$ & $\eta^{1/2}_{kh}$ & $\eta^{2/2}_{kh}$ \\
\hline
  $256$ &    $896$ & $1.08154812e-02$ & $8.82511469e-04$ & $4.84743239e-01$ & $1.88545655e-01$ \\
  $512$ &   $3584$ & $2.49545713e-03$ & $4.48390977e-04$ & $1.99398985e-01$ & $1.28924222e-01$ \\
 $1024$ &  $14366$ & $5.67139257e-04$ & $2.11190606e-04$ & $6.62293691e-02$ & $5.78821980e-02$ \\
 $2048$ &  $57344$ & $1.11745245e-04$ & $8.29692632e-05$ & $2.08643837e-02$ & $1.99891703e-02$ \\
\hline

\end{tabular}
\end{table}

\begin{table}[H]
\centering
\caption{Section \ref{sec:config_combustion}: Performance of the full error estimators under global refinement for $J_1$}
\label{table:comb_glob_full_reaction}
\begin{tabular}{l|l|l|l|l|l}
\hline
 $M$ & $N$ & $J(u)-J(u_{kh})$ & $\eta^{1/1}_{kh}$ & $\eta^{1/2}_{kh}$ & $\eta^{2/2}_{kh}$ \\
\hline
  $256$ &    $896$ & $1.08154812e-02$ & $9.41205097e-04$ & $2.41600914e-01$ & $9.24075938e-02$ \\
  $512$ &   $3584$ & $2.49545713e-03$ & $4.74660147e-04$ & $9.91323210e-02$ & $6.37932253e-02$ \\
 $1024$ &  $14366$ & $5.67139257e-04$ & $2.14565473e-04$ & $3.28205570e-02$ & $2.86419478e-02$ \\
 $2048$ &  $57344$ & $1.11745245e-04$ & $8.32187704e-05$ & $1.02590417e-02$ & $9.82187670e-03$ \\
\hline

\end{tabular}
\end{table}

\begin{figure}[H]
 \centering
 \subfloat{
  \includegraphics[trim ={70 380 90 50},clip, width=0.48\textwidth]{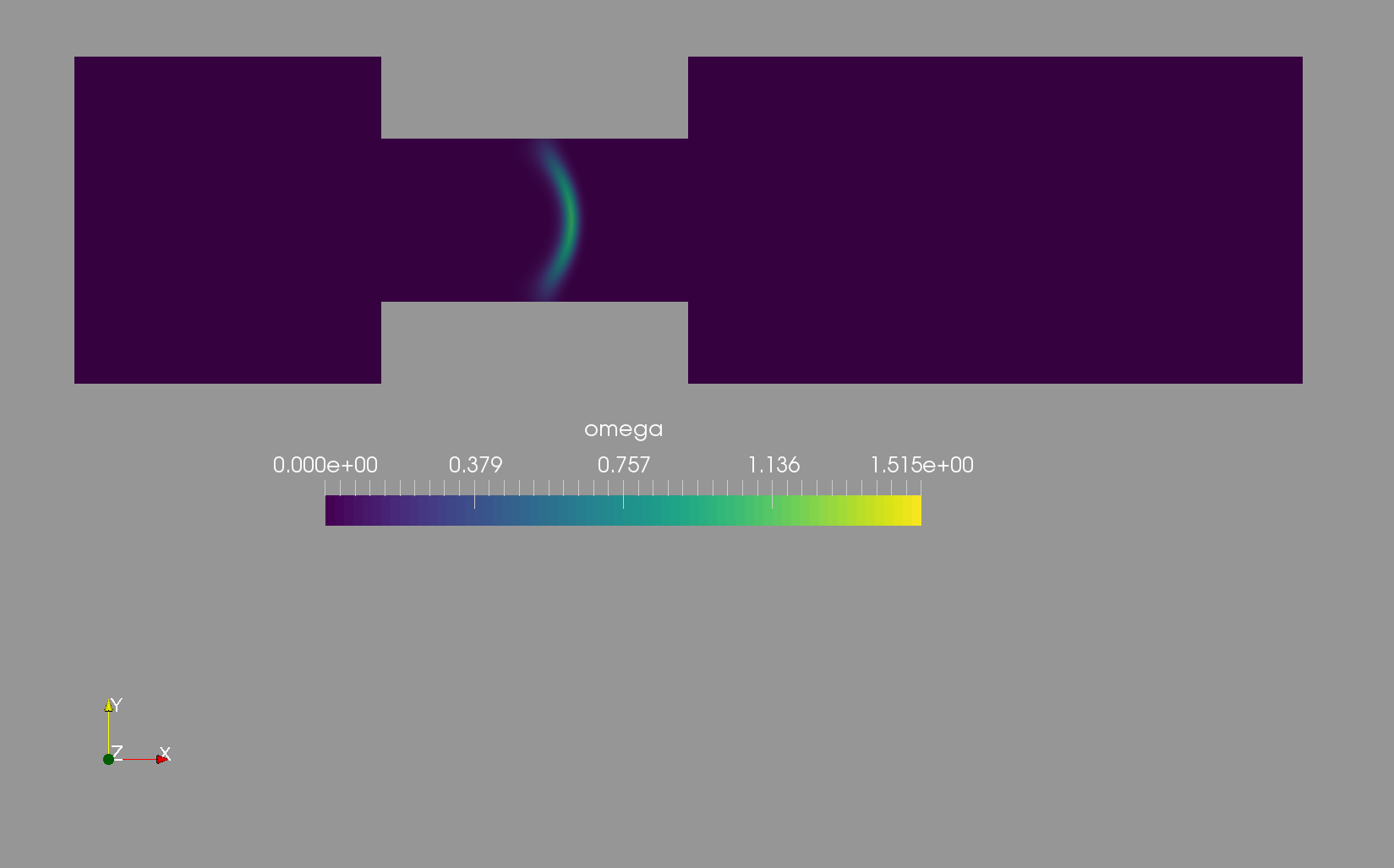}
 }
 \subfloat{
  \includegraphics[trim ={70 380 90 50},clip, width=0.48\textwidth]{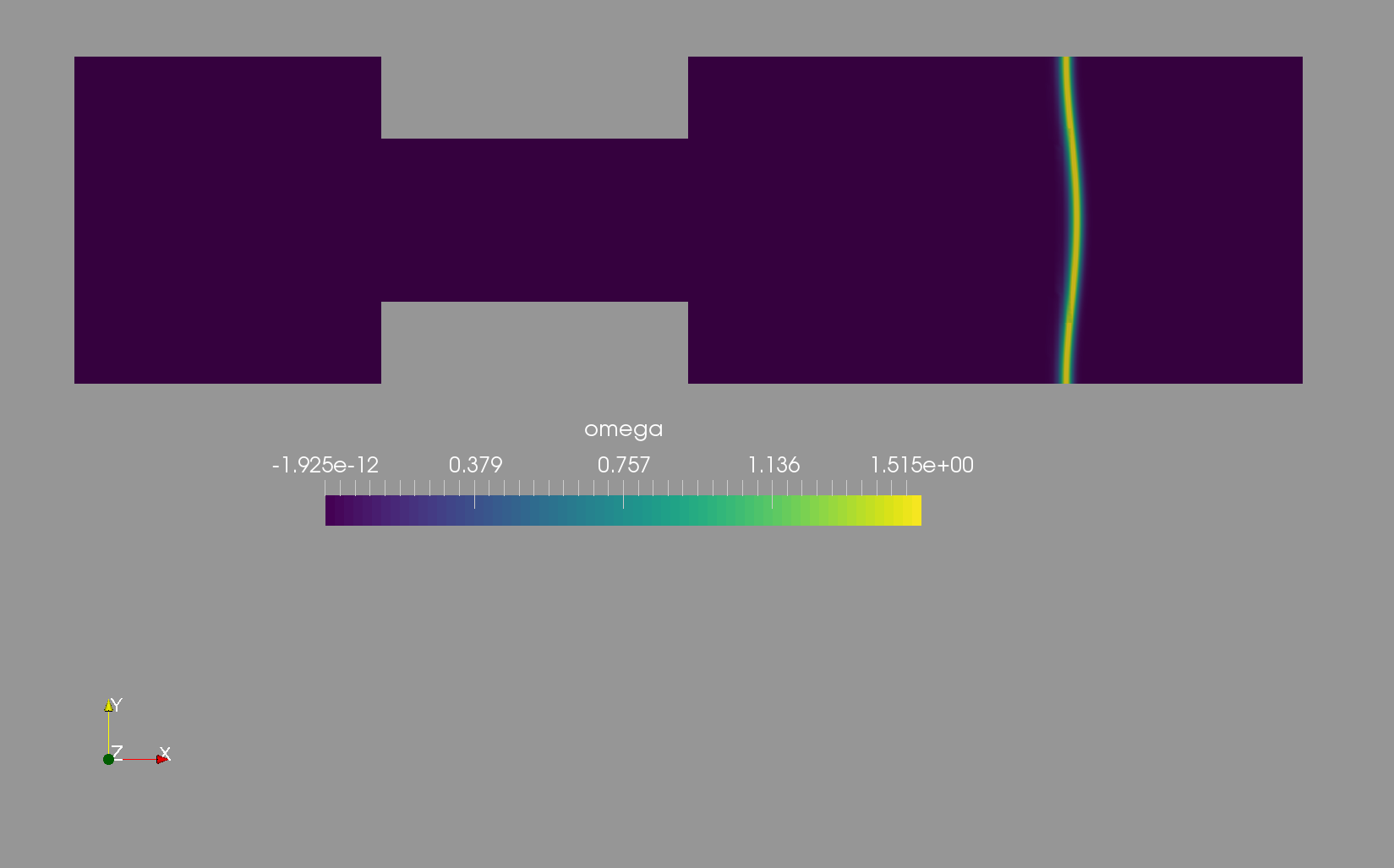}
 }\\
  \subfloat{
  \includegraphics[trim ={70 380 90 50},clip, width=0.48\textwidth]{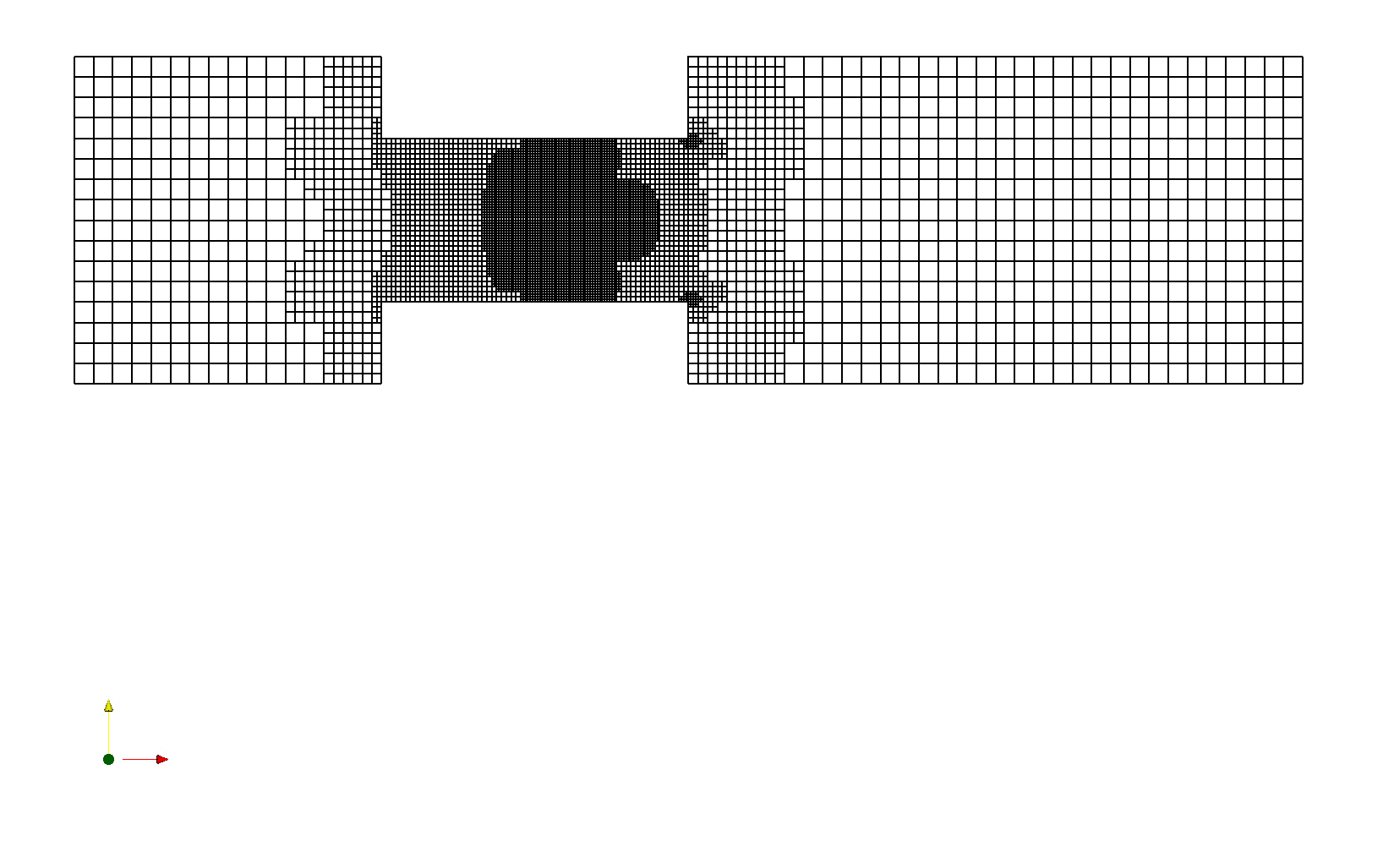}
 }
 \subfloat{
  \includegraphics[trim ={70 380 90 50},clip, width=0.48\textwidth]{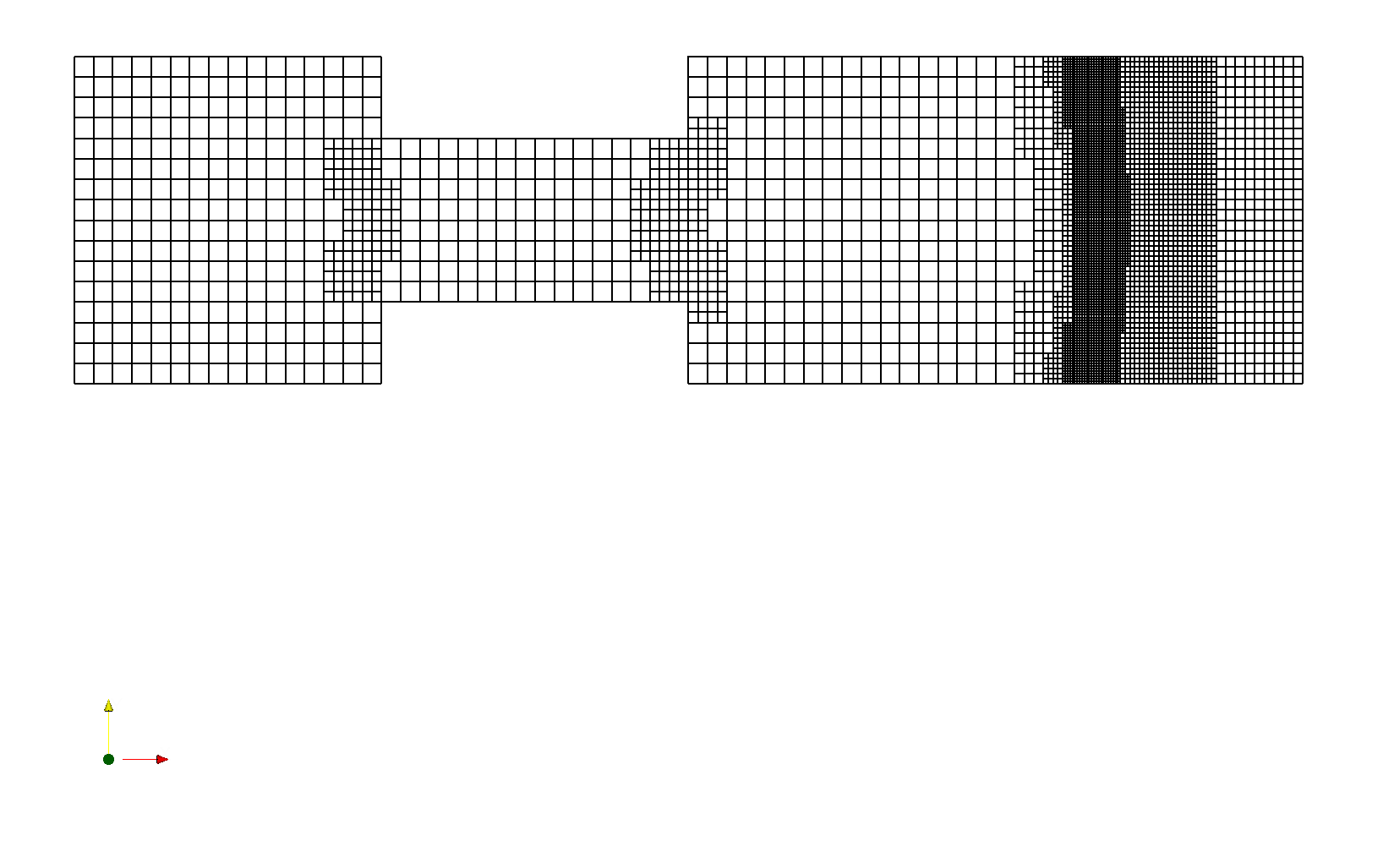}
 }
 \caption{Section \ref{sec:config_combustion}: reaction rate and grid at $t=20$ (left) and $t=60$ (right)}
 \label{figure:omega_flame}
\end{figure}

Figure \ref{figure:omega_flame} shows the reaction rate 
and the corresponding grids for two different time points. 
We can see that the grid evolves nicely and follows the combustion reaction. 
This shows that our localization works well in capturing the physics and refining accordingly.

\subsubsection{Discussion of findings for $J_2$}
\begin{figure}[H]
 \centering
 \subfloat{
  \includegraphics[width=0.48\textwidth]{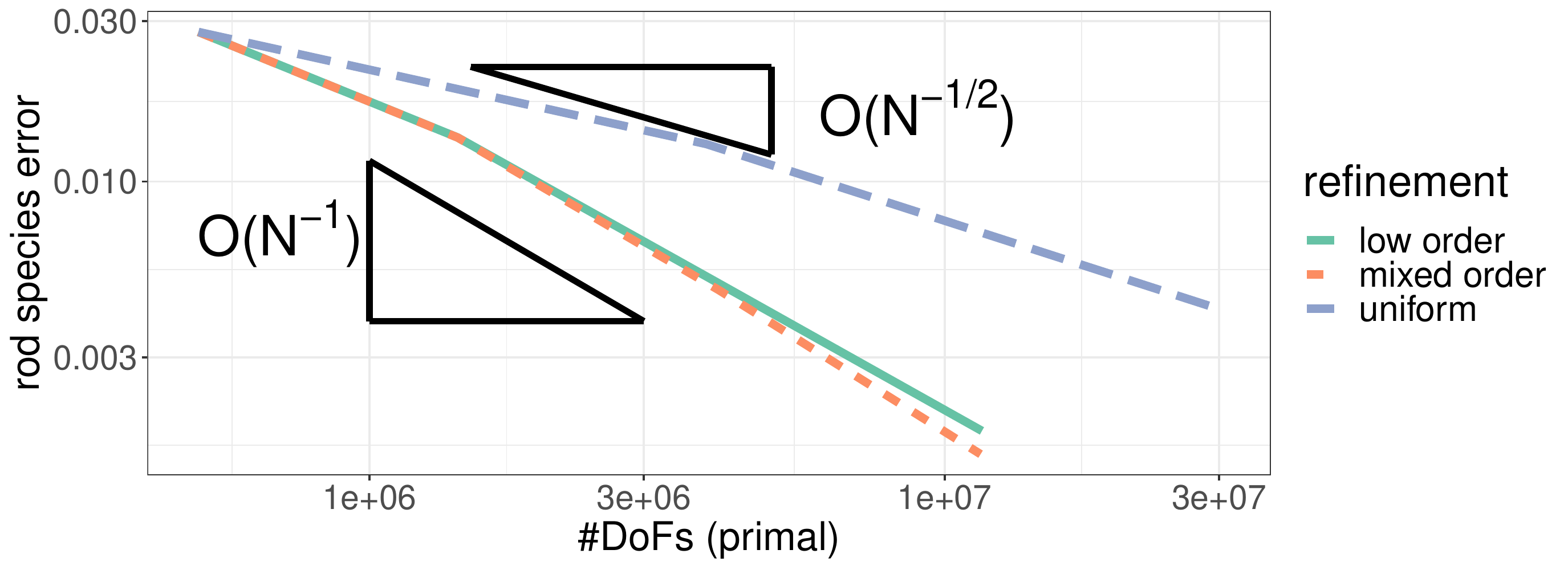}
 
 }
 \subfloat{
  \includegraphics[width=0.48\textwidth]{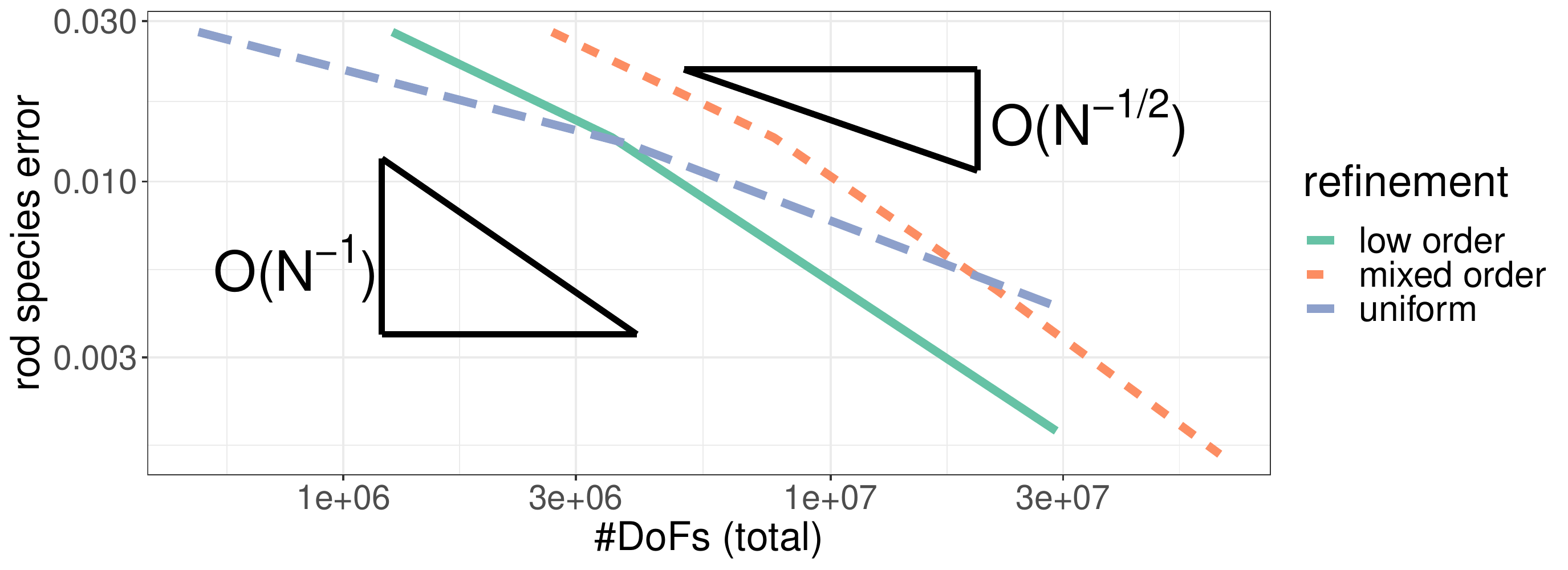}
 }
 \caption{Section \ref{sec:config_combustion}: Error convergence for the species concentration functional. 
 On the left only the number of unknowns for the primal problem
 and on the right all unknowns are taken into account.}
 \label{figure:comb_rod_convergence}
\end{figure}

Tables \ref{table:comb_glob_prim_species}, \ref{table:comb_glob_adj_species} 
and \ref{table:comb_glob_full_species} show the behaviour of the primal, adjoint and full estimators for $J_2$ respectively.
The equal low order approach shows a similar behaviour to the previous functional.
However, the overestimation of the adjoint estimators for the mixed order and high order approach is not as bad as for the nonlinear functional.
On the other hand their respective primal estimators are underestimating the error by quite a bit. 
In the full estimators these over- and underestimations are cancelling out quite nicely such that this estimator would be more useful here.
For this reason the simulations for Figure \ref{figure:comb_rod_convergence}
were done by using the full estimator for adaptivity 
in Algorithm \ref{algo_mark_refine_timestep} with the same marking 
strategy as before. As with the first functional we see that both approaches yield comparable results when only 
counting the memory cost for solving the primal problem. We also see that both approaches eventually outperform uniform refinement 
when taking the total cost into account. But the equal low order approach is again and still the best choice.

\begin{table}[H]
\centering
\caption{Section \ref{sec:config_combustion}: Performance of the primal error estimators under global refinement for $J_2$}
\label{table:comb_glob_prim_species}
\begin{tabular}{l|l|l|l|l|l}
\hline
 $M$ & $N$ & $J(u)-J(u_{kh})$ & $\eta^{1/1}_{kh}$ & $\eta^{1/2}_{kh}$ & $\eta^{2/2}_{kh}$ \\
\hline
  $256$ &    $896$ & $2.78670640e-02$ & $1.57134836e-02$ & $1.68459624e-03$ & $2.40117972e-03$\\
  $512$ &   $3584$ & $1.29349400e-02$ & $8.36950568e-03$ & $7.70832432e-04$ & $7.82054826e-04$\\
 $1024$ &  $14366$ & $4.17247100e-03$ & $3.88418655e-03$ & $3.57626628e-04$ & $3.52168282e-04$\\
 $2048$ &  $57344$ & $1.07518500e-03$ & $1.63053161e-03$ & $1.77990354e-04$ & $1.76929659e-04$\\
\hline

\end{tabular}
\end{table}

\begin{table}[H]
\centering
\caption{Section \ref{sec:config_combustion}: Performance of the adjoint error estimators under global refinement for $J_2$}
\label{table:comb_glob_adj_species}
\begin{tabular}{l|l|l|l|l|l}
\hline
 $M$ & $N$ & $J(u)-J(u_{kh})$ & $\eta^{1/1}_{kh}$ & $\eta^{1/2}_{kh}$ & $\eta^{2/2}_{kh}$ \\
\hline
  $256$ &    $896$ & $2.78670640e-02$ & $1.46726666e-02$ & $1.51762190e-01$ & $6.21329632e-02$ \\
  $512$ &   $3584$ & $1.29349400e-02$ & $7.65035908e-03$ & $4.68479946e-02$ & $3.34542288e-02$ \\
 $1024$ &  $14366$ & $4.17247100e-03$ & $3.79109953e-03$ & $1.26131290e-02$ & $1.13551533e-02$ \\
 $2048$ &  $57344$ & $1.07518500e-03$ & $1.62278719e-03$ & $3.29583742e-03$ & $3.19336781e-03$ \\
\hline

\end{tabular}
\end{table}

\begin{table}[H]
\centering
\caption{Section \ref{sec:config_combustion}: Performance of the full error estimators under global refinement for $J_2$}
\label{table:comb_glob_full_species}
\begin{tabular}{l|l|l|l|l|l}
\hline
 $M$ & $N$ & $J(u)-J(u_{kh})$ & $\eta^{1/1}_{kh}$ & $\eta^{1/2}_{kh}$ & $\eta^{2/2}_{kh}$ \\
\hline
  $256$ &    $896$ & $2.78670640e-02$ & $1.51930751e-02$ & $7.53851823e-02$ & $3.12572585e-02$ \\
  $512$ &   $3584$ & $1.29349400e-02$ & $8.00993238e-03$ & $2.31194532e-02$ & $1.64872919e-02$ \\
 $1024$ &  $14366$ & $4.17247100e-03$ & $3.83764304e-03$ & $6.13512088e-03$ & $5.51198299e-03$ \\
 $2048$ &  $57344$ & $1.07518500e-03$ & $1.62665940e-03$ & $1.55892353e-03$ & $1.50821907e-03$ \\
\hline

\end{tabular}
\end{table}

\begin{figure}[H]
 \centering
 \subfloat{
  \includegraphics[trim ={70 380 70 50},clip, width=0.48\textwidth]{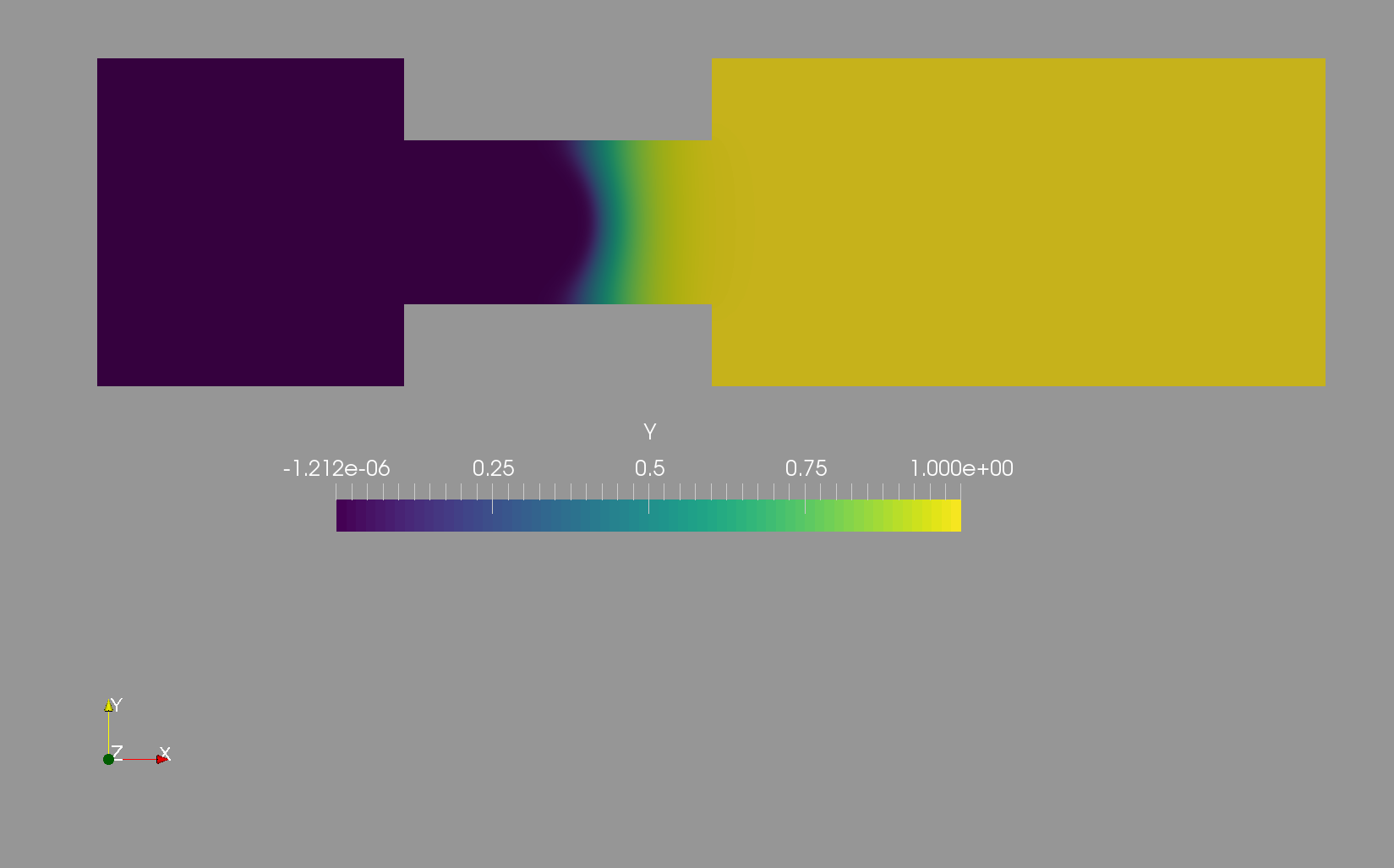}
 }
 \subfloat{
  \includegraphics[trim ={70 380 70 50},clip, width=0.48\textwidth]{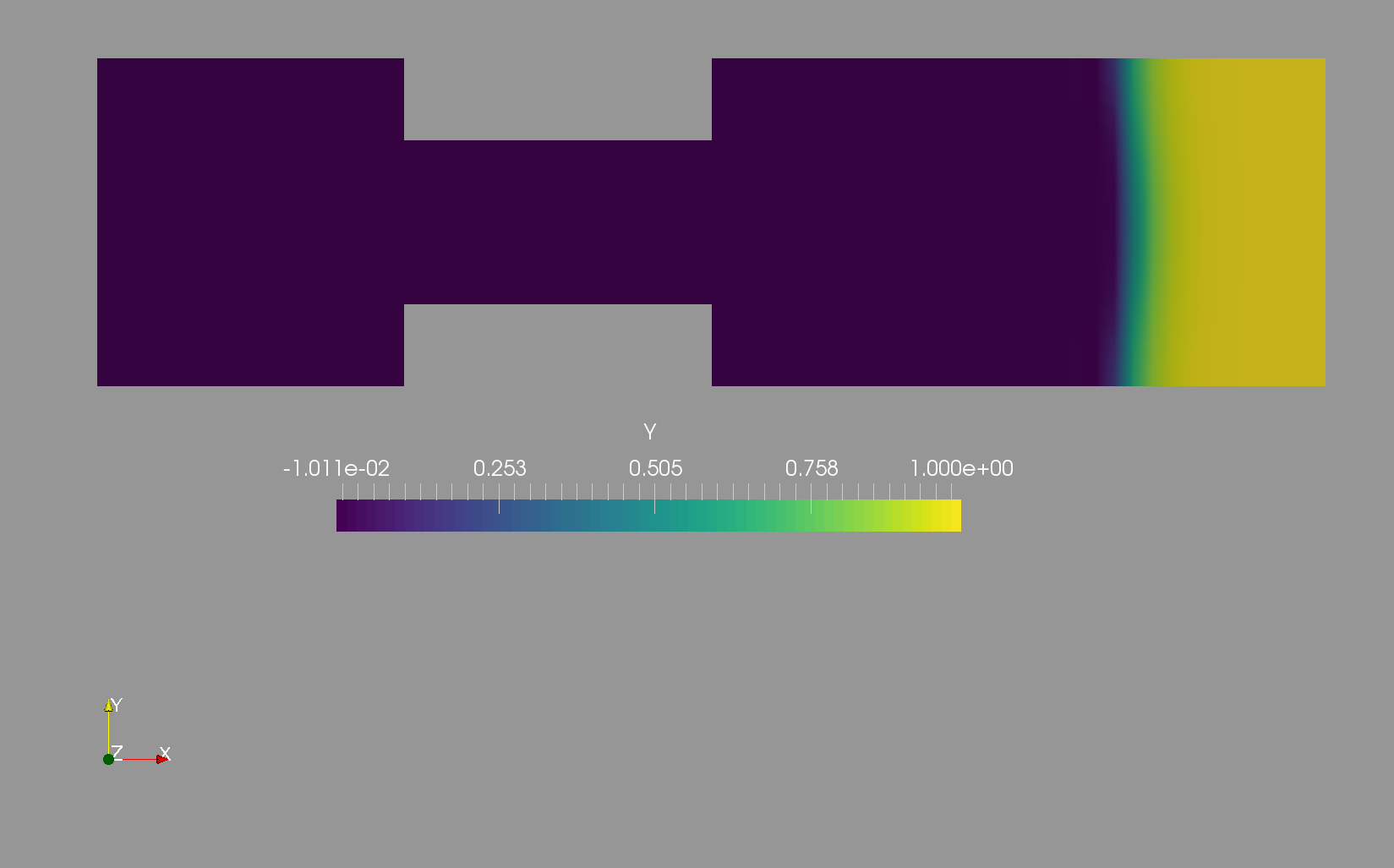}
 }\\
 \subfloat{
  \includegraphics[trim ={70 380 70 50},clip, width=0.48\textwidth]{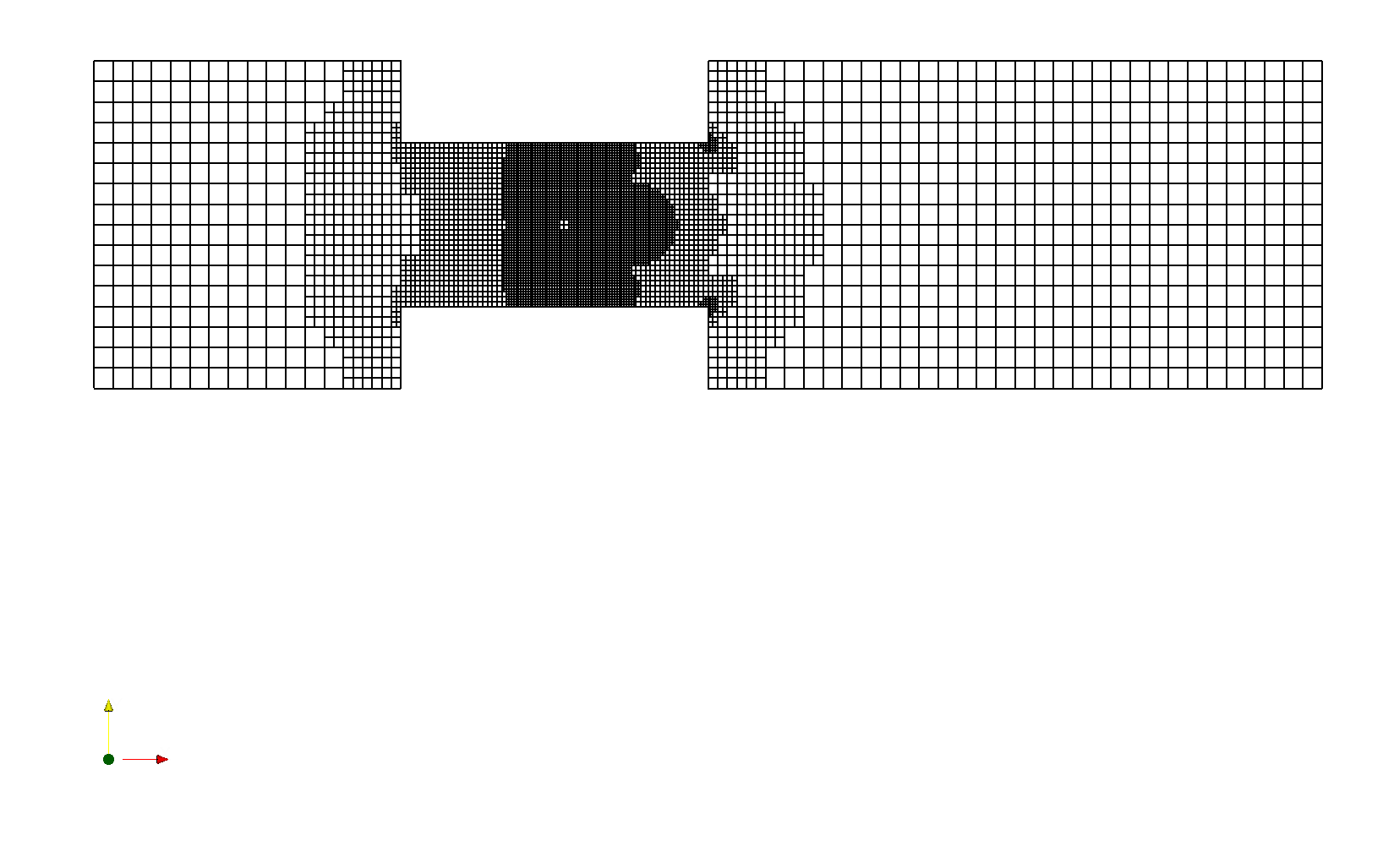}
 }
 \subfloat{
  \includegraphics[trim ={70 380 70 50},clip, width=0.48\textwidth]{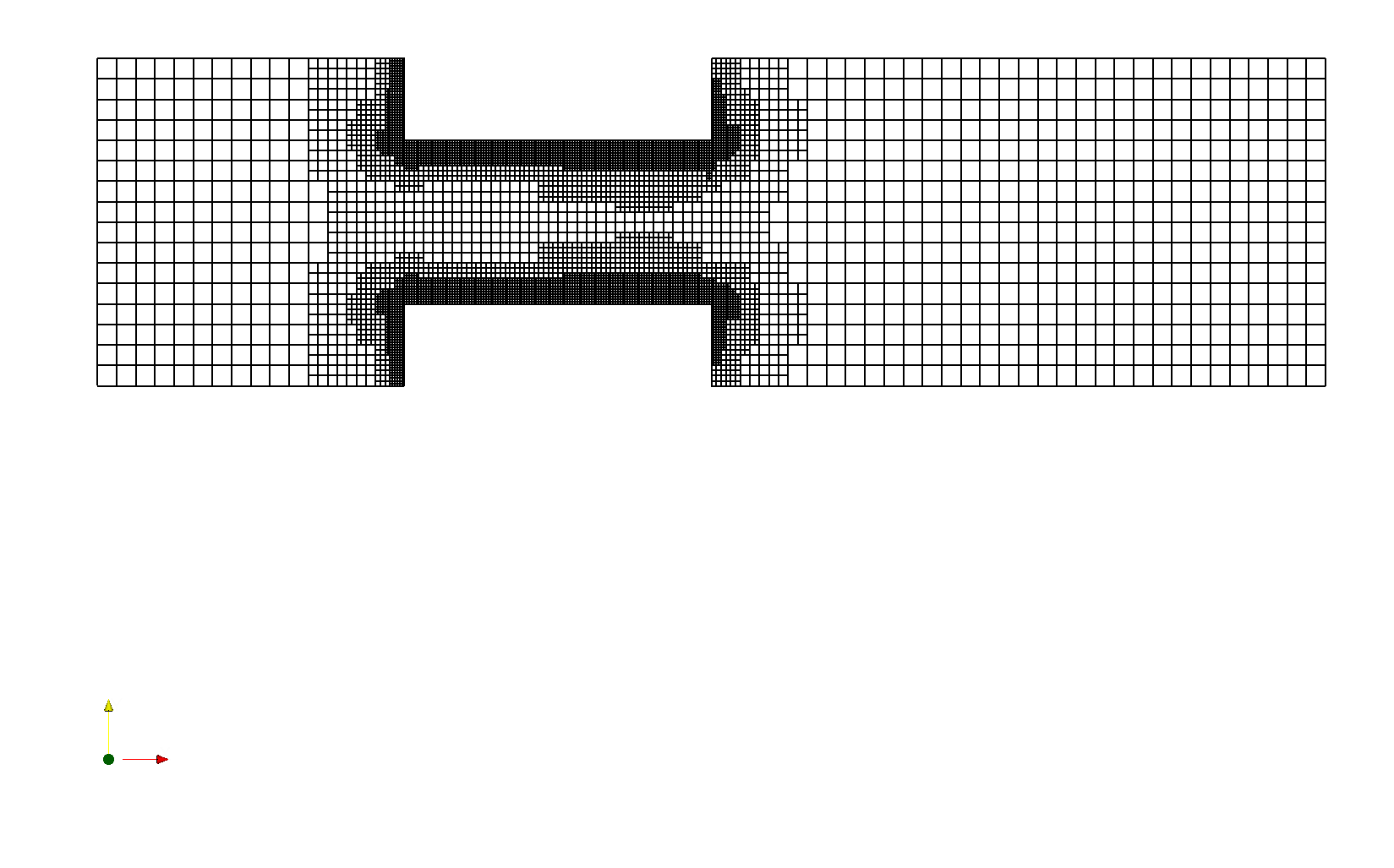}
 }
 \caption{Section \ref{sec:config_combustion}: Rod species concentration and corresponding grid at $t=20$ (left) and $t=60$ (right)}
 \label{figure:rod_flame}
\end{figure}

Figure \ref{figure:rod_flame}
shows the species concentration and the refined grids based on $J_2$ 
at different time steps. We see that along the cooled rods the grid again follows the combustion reaction. 
As that is the area where we have changes in the concentration this fits well. 
We also see that the mesh is refined around the cooled rod once the reaction moved past them. 
Together with the convergence behaviour we see again that the novel localization works well.

\section{Conclusions}
\label{Section: Conclusions}
In this work, we proposed partition-of-unity (PU) dual-weighted residual 
a posteriori error estimators and space-time adaptivity for linear and nonlinear 
partial differential equations. 
From the algorithmic side, the main novelties 
are the extension of the PU localization to space-time Galerkin finite element 
discretizations and the realization of split and joint error estimators.
From the implementation side, despite starting from pre-implementations 
in the DTM package \textit{dwr-diffusion} \cite{kocher.etal2019}
and deal.II \cite{arndt.etal2020}, extensive code developments and debugging
was necessary, which greatly exceed existing implementations,
specifically for the nonlinear features such as the nonlinear combustion PDE
as well as nonlinear goal functionals.
In three numerical examples, we studied in the detail the computational 
performance for the linear heat equation and also for 
a nonlinear low Mach number combustion problem. 
We also found that the equal low order approach yielded the best estimation and 
adaptive performance across the board and that a $cG(1)dG(0)$ PU is sufficient
for $cG(s)dG(r)$ discretizations of the primal problem.
Furthermore, an example of an immediate practical application of our framework can be found within the
excellence cluster PhoenixD\footnote{https://www.phoenixd.uni-hannover.de/} 
in which space-time methods and goal-oriented error estimation are of interest 
for the efficient solution of multiphysics problems and where the heat equation 
and the Navier-Stokes equations are needed.

\section*{Acknowledgments}
The authors thank Uwe K{\"o}cher (Hamburg), 
Julian Roth (IfAM), Johannes Lankeit (IfAM), and Bernhard Endtmayer (IfAM) for fruitful discussions. Moreover, we thank the anonymous referee 
for excellent questions that helped to improve the manuscript.

\section*{Declarations}
This work is funded by the Deutsche Forschungsgemeinschaft 
(DFG) under Germany’s Excellence Strategy within 
  the Cluster of Excellence PhoenixD (EXC 2122, Project ID 390833453).

\section*{Data availibility}
The algorithms and numerical results are implemented open source 
and can be found on \url{https://github.com/jpthiele/pu-dwr-diffusion} and 
\url{https://github.com/jpthiele/pu-dwr-combustion} respectively,
and follow good practices of sustainable research software developments \cite{Anztetal20}.


\bibliography{./Intro_nonstat_PU}

\begin{thebibliography}{10}

\bibitem{MUMPS}
P.~Amestoy, I.~S. Duff, J.~Koster, and J.-Y. L'Excellent.
\newblock A fully asynchronous multifrontal solver using distributed dynamic
  scheduling.
\newblock {\em SIAM Journal on Matrix Analysis and Applications}, 23(1):15--41,
  2001.

\bibitem{Anztetal20}
H.~Anzt, F.~Bach, S.~Druskat, F.~Löffler, A.~Loewe, B.~Renard, G.~Seemann,
  A.~Struck, E.~Achhammer, P.~Aggarwal, F.~Appel, M.~Bader, L.~Brusch,
  C.~Busse, G.~Chourdakis, P.~Dabrowski, P.~Ebert, B.~Flemisch, S.~Friedl,
  B.~Fritzsch, M.~Funk, V.~Gast, F.~Goth, J.~Grad, S.~Hermann, F.~Hohmann,
  S.~Janosch, D.~Kutra, J.~Linxweiler, T.~Muth, W.~Peters-Kottig, F.~Rack,
  F.~Raters, S.~Rave, G.~Reina, M.~Reißig, T.~Ropinski, J.~Schaarschmidt,
  H.~Seibold, J.~Thiele, B.~Uekerman, S.~Unger, and R.~Weeber.
\newblock An environment for sustainable research software in germany and
  beyond: current state, open challenges, and call for action [version 1; peer
  review: awaiting peer review].
\newblock {\em F1000Research}, 9(295), 2020.

\bibitem{arndt.etal2020}
D.~Arndt, W.~Bangerth, B.~Blais, T.~C. Clevenger, M.~Fehling, A.~V. Grayver,
  T.~Heister, L.~Heltai, M.~Kronbichler, M.~Maier, P.~Munch, J.-P. Pelteret,
  R.~Rastak, I.~Tomas, B.~Turcksin, Z.~Wang, and D.~Wells.
\newblock The deal.{II} library, {Version} 9.2.
\newblock {\em Journal of Numerical Mathematics}, 28(3):131--146, Sept. 2020.
\newblock Publisher: De Gruyter Section: Journal of Numerical Mathematics.

\bibitem{bangerth.etal2010}
W.~Bangerth, M.~Geiger, and R.~Rannacher.
\newblock Adaptive {Galerkin} {Finite} {Element} {Methods} for the {Wave}
  {Equation}.
\newblock {\em Computational Methods in Applied Mathematics}, 10(1):3--48, Jan.
  2010.
\newblock Publisher: De Gruyter Section: Computational Methods in Applied
  Mathematics.

\bibitem{bangerth.rannacher2003}
W.~Bangerth and R.~Rannacher.
\newblock {\em Adaptive finite element methods for differential equations}.
\newblock Lectures in mathematics {ETH} {Z\"{u}rich}. Birkh\"{a}user Verlag,
  Basel ; Boston, 2003.
\newblock OCLC: ocm51447498.

\bibitem{bause.etal2017}
M.~Bause, F.~A. Radu, and U.~K\"{o}cher.
\newblock Space--time finite element approximation of the {Biot} poroelasticity
  system with iterative coupling.
\newblock {\em Computer Methods in Applied Mechanics and Engineering},
  320:745--768, June 2017.

\bibitem{becker.etal2007}
R.~Becker, D.~Meidner, and B.~Vexler.
\newblock Efficient numerical solution of parabolic optimization problems by
  finite element methods.
\newblock {\em Optimization Methods and Software}, 22(5):813--833, Oct. 2007.
\newblock Publisher: Taylor \& Francis \_eprint:
  https://doi.org/10.1080/10556780701228532.

\bibitem{becker.rannacher1996}
R.~Becker and R.~Rannacher.
\newblock Weighted {A} {Posteriori} {Error} {Control} in {FE} {Methods}.
\newblock page~16, 1996.

\bibitem{becker.rannacher2001}
R.~Becker and R.~Rannacher.
\newblock An optimal control approach to \textit{a posteriori} error estimation
  in finite element methods.
\newblock {\em Acta Numerica}, 10:1--102, May 2001.

\bibitem{besier.rannacher2012}
M.~Besier and R.~Rannacher.
\newblock Goal-oriented space--time adaptivity in the finite element {Galerkin}
  method for the computation of nonstationary incompressible flow.
\newblock {\em International Journal for Numerical Methods in Fluids},
  70(9):1139--1166, 2012.
\newblock \_eprint: https://onlinelibrary.wiley.com/doi/pdf/10.1002/fld.2735.

\bibitem{braack.ern2003}
M.~Braack and A.~Ern.
\newblock A {Posteriori} {Control} of {Modeling} {Errors} and {Discretization}
  {Errors}.
\newblock {\em Multiscale Modeling \& Simulation}, 1(2):221--238, Jan. 2003.
\newblock Publisher: Society for Industrial and Applied Mathematics.

\bibitem{CaVe99}
C.~Carstensen and R.~Verf\"{u}rth.
\newblock Edge residuals dominate a posteriori error estimates for low order
  finite element methods.
\newblock {\em SIAM J. Numer. Anal.}, 36(5):1571--1587, 1999.

\bibitem{davis2004}
T.~A. Davis.
\newblock Algorithm 832: {UMFPACK} {V4}.3---an unsymmetric-pattern multifrontal
  method.
\newblock {\em ACM Transactions on Mathematical Software}, 30(2):196--199, June
  2004.

\bibitem{DiPietroErn2011}
D.~A. {Di Pietro} and A.~Ern.
\newblock {\em {Mathematical aspects of discontinuous Galerkin methods}},
  volume~69.
\newblock Springer Science \& Business Media, 2011.

\bibitem{DoerflerFindeisenWieners+2016+409+428}
W.~D\"{o}rfler, S.~Findeisen, and C.~Wieners.
\newblock Space-time discontinuous galerkin discretizations for linear
  first-order hyperbolic evolution systems.
\newblock {\em Computational Methods in Applied Mathematics}, 16(3):409--428,
  2016.

\bibitem{DoeWieZie21}
W.~D\"{o}rfler, C.~Wieners, and D.~Ziegler.
\newblock Parallel space-time solutions for the linear visco-acoustic and
  visco-elastic wave equation.
\newblock In W.~E. Nagel, D.~H. Kr\"{o}ner, and M.~M. Resch, editors, {\em High
  Performance Computing in Science and Engineering '19}, pages 589--599, Cham,
  2021. Springer International Publishing.

\bibitem{endtmayer2023goaloriented}
B.~Endtmayer, U.~Langer, and A.~Schafelner.
\newblock Goal-oriented adaptive space-time finite element methods for
  regularized parabolic p-laplace problems.
\newblock {\em arXiv:2306.07167}, 2023.

\bibitem{EndtLaWi18}
B.~Endtmayer, U.~Langer, and T.~Wick.
\newblock {Multigoal-Oriented Error Estimates for Non-linear Problems}.
\newblock {\em Journal of Numerical Mathematics}, 27(4):215--236, 2019.

\bibitem{EndtLaWi21_smart}
B.~Endtmayer, U.~Langer, and T.~Wick.
\newblock Reliability and efficiency of {DWR}-type a posteriori error estimates
  with smart sensitivity weight recovering.
\newblock {\em Computational Methods in Applied Mathematics}, 21(2), 2021.

\bibitem{ErEstHaJoh09}
K.~Eriksson, D.~Estep, P.~Hansbo, and C.~Johnson.
\newblock {\em Computational Differential Equations}.
\newblock Cambridge University Press, 2009.
\newblock http://www.csc.kth.se/~jjan/private/cde.pdf.

\bibitem{ErJo91}
K.~Eriksson and C.~Johnson.
\newblock {Adaptive Finite Element Methods for Parabolic Problems I: A Linear
  Model Problem}.
\newblock {\em SIAM Journal on Numerical Analysis}, 28(1):43--77, 1991.

\bibitem{ErJo95}
K.~Eriksson and C.~Johnson.
\newblock {Adaptive Finite Element Methods for Parabolic Problems II: Optimal
  Error Estimates in Linfty L2 and Linfty Linfty}.
\newblock {\em SIAM Journal on Numerical Analysis}, 32(3):706--740, 1995.

\bibitem{ErJoLo04}
K.~Eriksson, C.~Johnson, and A.~Logg.
\newblock {\em Adaptive Computational Methods for Parabolic Problems},
  chapter~24.
\newblock American Cancer Society, 2004.

\bibitem{Fai17}
L.~Failer.
\newblock {\em Optimal Control of Time-Dependent Nonlinear Fluid-Structure
  Interaction}.
\newblock PhD thesis, Technical University Munich, 2017.

\bibitem{failer.etal2016}
L.~Failer, D.~Meidner, and B.~Vexler.
\newblock Optimal {Control} of a {Linear} {Unsteady} {Fluid}--{Structure}
  {Interaction} {Problem}.
\newblock {\em Journal of Optimization Theory and Applications}, 170(1):1--27,
  July 2016.

\bibitem{failer.wick2018}
L.~Failer and T.~Wick.
\newblock Adaptive time-step control for nonlinear fluid--structure
  interaction.
\newblock {\em Journal of Computational Physics}, 366:448--477, Aug. 2018.

\bibitem{gander.neumuller2016}
M.~J. Gander and M.~Neum\"{u}ller.
\newblock Analysis of a {New} {Space}-{Time} {Parallel} {Multigrid} {Algorithm}
  for {Parabolic} {Problems}.
\newblock {\em SIAM Journal on Scientific Computing}, 38(4):A2173--A2208, Jan.
  2016.
\newblock Publisher: Society for Industrial and Applied Mathematics.

\bibitem{goll.etal2015}
C.~Goll, R.~Rannacher, and W.~Wollner.
\newblock The {Damped} {Crank}--{Nicolson} {Time}-{Marching} {Scheme} for the
  {Adaptive} {Solution} of the {Black}--{Scholes} {Equation}.
\newblock {SSRN} {Scholarly} {Paper} ID 2795622, Social Science Research
  Network, Rochester, NY, Apr. 2015.

\bibitem{hartmann1998}
R.~Hartmann.
\newblock A-posteriori {Fehlersch\"{a}tzung} und adaptive {Schrittweitein}- und
  {Ortsgittersteuerung} bei {Galerkin}-{Verfahren} f\"{u}r die
  {W\"{a}rmeleitungsgleichung}.
\newblock Diplomarbeit, University of Heidelberg, 1998.

\bibitem{hubner.etal2004}
B.~H\"{u}bner, E.~Walhorn, and D.~Dinkler.
\newblock A monolithic approach to fluid--structure interaction using
  space--time finite elements.
\newblock {\em Computer Methods in Applied Mechanics and Engineering},
  193(23):2087--2104, June 2004.

\bibitem{hughes.hulbert1988}
T.~J.~R. Hughes and G.~M. Hulbert.
\newblock Space-time finite element methods for elastodynamics: {Formulations}
  and error estimates.
\newblock {\em Computer Methods in Applied Mechanics and Engineering},
  66(3):339--363, Feb. 1988.

\bibitem{HuHu90}
G.~M. Hulbert and T.~J. Hughes.
\newblock Space-time finite element methods for second-order hyperbolic
  equations.
\newblock {\em Computer Methods in Applied Mechanics and Engineering},
  84(3):327--348, 1990.

\bibitem{KhiSteiWi22_JCP}
D.~Khimin, M.~Steinbach, and T.~Wick.
\newblock Space-time formulation, discretization, and computational performance
  studies for phase-field fracture optimal control problems.
\newblock {\em Journal of Computational Physics}, 470:111554, 2022.

\bibitem{Koecher2015}
U.~K\"{o}cher.
\newblock {\em Variational space-time methods for the elastic wave equation and
  the diffusion equation}.
\newblock PhD thesis, Helmut-Schmidt-Universit\"{a}t, 2015.

\bibitem{kocher.etal2019}
U.~K\"{o}cher, M.~P. Bruchh\"{a}user, and M.~Bause.
\newblock Efficient and scalable data structures and algorithms for
  goal-oriented adaptivity of space--time {FEM} codes.
\newblock {\em SoftwareX}, 10:100239, July 2019.

\bibitem{Lang2001}
J.~Lang, editor.
\newblock {\em Adaptive Multilevel Solution of Nonlinear Parabolic PDE
  Systems}.
\newblock Lecture Notes in Computational Science and Engineering book series
  (LNCSE, volume 16). Springer, 2001.

\bibitem{LaSchaf20}
U.~Langer and A.~Schafelner.
\newblock Adaptive space-time finite element methods for non-autonomous
  parabolic problems with distributional sources.
\newblock {\em Comput. Methods Appl. Math.}, 2020.

\bibitem{LaSchaf21b}
U.~Langer and A.~Schafelner.
\newblock Adaptive space-time finite element methods for parabolic optimal
  control problems.
\newblock {\em Journal of Numerical Mathematics}, 2021.

\bibitem{LaSchaf21}
U.~Langer and A.~Schafelner.
\newblock Space-time hexahedral finite element methods for parabolic evolution
  problems, 2021.
\newblock DK Report 2021-05.

\bibitem{LaStein19}
U.~Langer and O.~Steinbach, editors.
\newblock {\em Space-time methods: {A}pplication to Partial Differential
  Equations}.
\newblock volume 25 of Radon Series on Computational and Applied Mathematics,
  Berlin. de Gruyter, 2019.

\bibitem{langer.etal2021a}
U.~Langer, O.~Steinbach, F.~Tr\"{o}ltzsch, and H.~Yang.
\newblock Space-{Time} {Finite} {Element} {Discretization} of {Parabolic}
  {Optimal} {Control} {Problems} with {Energy} {Regularization}.
\newblock {\em SIAM Journal on Numerical Analysis}, 59(2):675--695, Jan. 2021.
\newblock Publisher: Society for Industrial and Applied Mathematics.

\bibitem{langer.etal2021}
U.~Langer, O.~Steinbach, F.~Tr\"{o}ltzsch, and H.~Yang.
\newblock Unstructured {Space}-{Time} {Finite} {Element} {Methods} for
  {Optimal} {Control} of {Parabolic} {Equations}.
\newblock {\em SIAM Journal on Scientific Computing}, 43(2):A744--A771, Jan.
  2021.
\newblock Publisher: Society for Industrial and Applied Mathematics.

\bibitem{meidner2007}
D.~Meidner.
\newblock {\em Adaptive {Space}-{Time} {Finite} {Element} {Methods} for
  {Optimization} {Problems} {Governed} by {Nonlinear} {Parabolic} {Systems}}.
\newblock Dissertation, 2007.

\bibitem{MeiRaVih109}
D.~Meidner, R.~Rannacher, and J.~Vihharev.
\newblock Goal-oriented error control of the iterative solution of finite
  element equations.
\newblock {\em Journal of Numerical Mathematics}, 17:143--172, 2009.

\bibitem{MeiRi14}
D.~Meidner and T.~Richter.
\newblock Goal-oriented error estimation for the fractional step theta scheme.
\newblock {\em Comput. Methods Appl. Math.}, 14(2):203--230, 2014.

\bibitem{meidner.richter2015}
D.~Meidner and T.~Richter.
\newblock A posteriori error estimation for the fractional step theta
  discretization of the incompressible {Navier}--{Stokes} equations.
\newblock {\em Computer Methods in Applied Mechanics and Engineering},
  288:45--59, May 2015.

\bibitem{meidner.vexler2007}
D.~Meidner and B.~Vexler.
\newblock Adaptive {Space}‐{Time} {Finite} {Element} {Methods} for
  {Parabolic} {Optimization} {Problems}.
\newblock {\em SIAM Journal on Control and Optimization}, 46(1):116--142, Jan.
  2007.
\newblock Publisher: Society for Industrial and Applied Mathematics.

\bibitem{NeiVe12}
I.~Neitzel and B.~Vexler.
\newblock A priori error estimates for space-time finite element discretization
  of semilinear parabolic optimal control problems.
\newblock {\em Numer. Math.}, 120(2):345--386, 2012.

\bibitem{neumuller2013}
M.~Neum\"{u}ller.
\newblock {\em Space-time methods fast solvers and applications}.
\newblock PhD thesis, TU Graz, June 2013.

\bibitem{Od69b}
J.~Oden.
\newblock {A general theory of finite elements II. Applications}.
\newblock {\em Internat. J. Numer. Methods Engrg.}, 1:247--259, 1969.

\bibitem{rannacher.suttmeier1999}
R.~Rannacher and F.-T. Suttmeier.
\newblock A posteriori error estimation and mesh adaptation for finite element
  models in elasto-plasticity.
\newblock {\em Computer Methods in Applied Mechanics and Engineering},
  176(1):333--361, July 1999.

\bibitem{RaVih13b}
R.~Rannacher and J.~Vihharev.
\newblock Adaptive finite element analysis of nonlinear problems: balancing of
  discretization and iteration errors.
\newblock {\em Journal of Numerical Mathematics}, 21(1):23--61, 2013.

\bibitem{richter.wick2015}
T.~Richter and T.~Wick.
\newblock Variational localizations of the dual weighted residual estimator.
\newblock {\em Journal of Computational and Applied Mathematics}, 279:192--208,
  May 2015.

\bibitem{RoThieKoeWi23}
J.~Roth, J.~P. Thiele, U.~K\"{o}cher, and T.~Wick.
\newblock Tensor-product space-time goal-oriented error control and adaptivity
  with partition-of-unity dual-weighted residuals for nonstationary flow
  problems.
\newblock {\em Computational Methods in Applied Mathematics}, 2023.

\bibitem{Schaf22}
A.~Schafelner.
\newblock {\em Space-time Finite Element Methods}.
\newblock PhD thesis, Johannes Kepler University Linz, 2022.

\bibitem{schmich2009}
M.~Schmich.
\newblock {\em Adaptive {Finite} {Element} {Methods} for {Computing}
  {Nonstationary} {Incompressible} {Flows}}.
\newblock PhD thesis, Dec. 2009.

\bibitem{schmich.vexler2008}
M.~Schmich and B.~Vexler.
\newblock Adaptivity with {Dynamic} {Meshes} for {Space}-{Time} {Finite}
  {Element} {Discretizations} of {Parabolic} {Equations}.
\newblock {\em SIAM Journal on Scientific Computing}, 30(1):369--393, Jan.
  2008.

\bibitem{SINGH2018893}
G.~Singh and M.~F. Wheeler.
\newblock A space--time domain decomposition approach using enhanced velocity
  mixed finite element method.
\newblock {\em Journal of Computational Physics}, 374:893--911, 2018.

\bibitem{SteiYa18}
O.~Steinbach and H.~Yang.
\newblock Comparison of algebraic multigrid methods for an adaptive space--time
  finite-element discretization of the heat equation in {3D} and {4D}.
\newblock {\em Numerical Linear Algebra with Applications}, 25(3):e2143, 2018.
\newblock e2143 nla.2143.

\bibitem{TaTe11}
K.~Takizawa and T.~Tezduyar.
\newblock Multiscale space--time fluid--structure interaction techniques.
\newblock {\em Computational Mechanics}, 48:247--267, 2011.

\bibitem{tezduyar.etal1992}
T.~E. Tezduyar, M.~Behr, S.~Mittal, and J.~Liou.
\newblock A new strategy for finite element computations involving moving
  boundaries and interfaces---{The} deforming-spatial-domain/space-time
  procedure: {II}. {Computation} of free-surface flows, two-liquid flows, and
  flows with drifting cylinders.
\newblock {\em Computer Methods in Applied Mechanics and Engineering},
  94(3):353--371, Feb. 1992.

\bibitem{TeSa07}
T.~E. Tezduyar and S.~Sathe.
\newblock Modelling of fluid--structure interactions with the space--time
  finite elements: Solution techniques.
\newblock {\em International Journal for Numerical Methods in Fluids},
  54(6‐8):855--900, 2007.

\bibitem{tezduyar.etal2006}
T.~E. Tezduyar, S.~Sathe, R.~Keedy, and K.~Stein.
\newblock Space--time finite element techniques for computation of
  fluid--structure interactions.
\newblock {\em Computer Methods in Applied Mechanics and Engineering},
  195(17):2002--2027, Mar. 2006.

\bibitem{tezduyar.etal2006a}
T.~E. Tezduyar, S.~Sathe, and K.~Stein.
\newblock Solution techniques for the fully discretized equations in
  computation of fluid--structure interactions with the space--time
  formulations.
\newblock {\em Computer Methods in Applied Mechanics and Engineering},
  195(41):5743--5753, Aug. 2006.

\bibitem{Thi23_phd}
J.~Thiele.
\newblock {\em Error-controlled space-time finite elements, algorithms, and
  implementations for nonstationary problems}.
\newblock PhD thesis, Leibniz University Hannover, in preparation, 2023.

\bibitem{ThiWi21_PAMM}
J.~P. Thiele and T.~Wick.
\newblock Space-time pu-dwr error control and adaptivity for the heat equation.
\newblock {\em PAMM}, 21(1):e202100174, 2021.

\bibitem{ThiWi23_ICO}
J.~P. Thiele and T.~Wick.
\newblock Space-time error control using a partition-of-unity dual-weighted
  residual method applied to low mach number combustion.
\newblock In J.~M. Melenk, I.~Perugia, J.~Sch\"{o}berl, and C.~Schwab, editors,
  {\em Spectral and High Order Methods for Partial Differential Equations
  ICOSAHOM 2020+1}, pages 509--520, Cham, 2023. Springer International
  Publishing.

\bibitem{Ve03}
R.~Verf\"{u}rth.
\newblock A posteriori error estimates for finite element discretizations of
  the heat equation.
\newblock 40:195--212, 2003.

\bibitem{Wi21_WCCM}
T.~Wick.
\newblock {On the Adjoint Equation in Fluid-Structure Interaction}.
\newblock WCCM-ECCOMAS2020, 2021.

\bibitem{Wi22_num_pde}
T.~Wick.
\newblock Numerical methods for partial differential equations.
\newblock Hannover : Institutionelles Repositorium der Leibniz Universit\"{a}t
  Hannover, DOI: https://doi.org/10.15488/11709, January 2022.

\end{thebibliography}
\bibliographystyle{abbrv}

\end{document}